\documentclass[11pt]{article}
\usepackage{amssymb}
\usepackage{amsmath}
\usepackage{color}
\usepackage{graphics}
\usepackage{hyperref}
\numberwithin{equation}{section}
\usepackage{mathrsfs}
\usepackage{ifpdf}
\ifpdf
\usepackage{epstopdf}
\fi
\def\tr{\mathrm{tr}}
\def\div{\rm div}
\def\Hdiv{\mathbb{H}(\mathbf{div};\Omega)}
\def\btau{\boldsymbol \tau}
\def\bu{\boldsymbol u}
\def\bv{\boldsymbol v}
\def\bw{\boldsymbol w}
\def\bt{\boldsymbol t}
\def\br{\boldsymbol r}

\def\bsig{\boldsymbol \sigma}

\def\bA{\mathbf A}
\def\bB{\mathbf B}
\def\bC{\mathbf C}
\def\GD{{\Gamma_\mathrm{D}}}
\def\GN{{\Gamma_\mathrm{N}}}
\def\D{{\mathrm{D}}}
\DeclareMathAlphabet{\mathpzc}{OT1}{pzc}{m}{it} 
\def\A{\text{\Large $\mathpzc{A}$}} 
\def\B{\text{\Large $\mathpzc{B}$}} 
\def\L{\text{\large $\mathscr{L}$}}
\def\CA{C_{\!\text{\footnotesize $\mathpzc{A}$}}}
\def\CB{C_{\!\text{\footnotesize $\mathpzc{B}$}}}
\def\LA{L_{\!\text{\footnotesize $\mathpzc{A}$}}}
\def\LB{L_{\!\text{\footnotesize $\mathpzc{B}$}}}
\def\CBt{\widetilde{C}_{\!\text{\footnotesize $\mathpzc{B}$}}}
\def\LAt{\widetilde{L}_{\!\text{\footnotesize $\mathpzc{A}$}}}
\def\ww{\boldsymbol{\omega}}
\newcommand{\T}{\mathcal{T}_{h}}
\newtheorem{theorem}{Theorem}[section]
\newtheorem{rem}{Remark}[section]
\newtheorem{lem}[theorem]{Lemma}

\newenvironment{proof}{\noindent{\it Proof.}}{\hfill$\square$}

\def\bphi{\boldsymbol{\phi}}
\def\bchi{\boldsymbol{\chi}}
\def\bvarphi{\boldsymbol{\varphi}}
\def\bpsi{\boldsymbol{\psi}}
\def\bphi{\boldsymbol{\phi}}
\def\alp{\boldsymbol{\alpha}}
\def\vbx{\vec{\boldsymbol{\xi}}}
\def\vby{\vec{\boldsymbol{\eta}}}
\def\vbz{\vec{\boldsymbol{\zeta}}}
\def\Div{\mathbf{div}}
\def\Cpp{C{\scriptsize $^{++}$}\,}

\newcommand\bsi{{\boldsymbol \sigma}}

\newcommand{\bn}{{\boldsymbol n}}

\newcommand\bg{{\boldsymbol g}}
\newcommand\bx{{\boldsymbol x}}

\newcommand\bdiv{{\mathbf{div}}}
\renewcommand\div{{\mathrm{div}}}
\newcommand\qan{{\quad\hbox{and}\quad}}
\newcommand\qin{{\quad\hbox{in}\quad}}
\newcommand\qon{{\quad\hbox{on}\quad}}
\newcommand\disp{\displaystyle}

\def\thanksMario{Secci\'on de Matem\'atica, Sede de Occidente, Universidad de Costa Rica,
                 San Ram\'on, Alajuela, Costa Rica, email: {\tt mario.alvarezguadamuz@ucr.ac.cr}.}
\def\thanksEligio{Departamento de Ciencias B\'asicas, Facultad de Ciencias, Universidad del
                  B\'io-B\'io, Campus Fernando May, Chill\'an, Chile, email: {\tt ecolmenares@ubiobio.cl}.}
\def\thanksFila{Escuela de Matem\'atica, Universidad Nacional, Campus Omar Dengo,
				Heredia, Costa Rica, email: {\tt filander.sequeira@una.cr}.}

\begin{document}
\title{Analysis of a  semi-augmented mixed finite element method for double-diffusive natural convection in porous media
\thanks{This work was partially supported by  Vicerrector\'ia de Investigaci\'on, through the
project C00-89,  Sede de Occidente, Universidad de Costa Rica; by 
CONICYT-Chile through the project 1190241; and by Universidad Nacional, Costa Rica,  through the project 0140-20.}}

\author{{\sc Mario Alvarez}\thanks{\thanksMario} \quad
        {\sc Eligio Colmenares}\thanks{\thanksEligio} \quad
        {\sc Fil\'ander A. Sequeira}\thanks{\thanksFila}}
	
\date{ }

\maketitle

\begin{abstract}
	\noindent In this paper we study a stationary double-diffusive natural convection problem in porous media given by a Navier-Stokes/Darcy type system,  for describing the velocity and the pressure, coupled to a vector advection-diffusion equation  describing the heat and substance concentration,   of a viscous fluid in a porous media with physical boundary conditions. The model problem is rewritten in terms of a first-order system, without the pressure,  based on the introduction of the strain tensor and  a nonlinear pseudo-stress tensor in the fluid equations. After a variational approach, the  resulting weak model is then augmented using appropriate redundant penalization terms for the fluid equations along with a standard primal formulation for the heat and substance concentration. Then, it is rewritten  as an equivalent  fixed-point problem. Well-posedness and uniqueness results for both the continuous and the discrete schemes are stated,  as well as the respective convergence result under certain regularity assumptions combined with the Lax-Milgram theorem, and the Banach and Brouwer fixed-point theorems. In particular, Raviart-Thomas elements of order $k$ are used for approximating the pseudo-stress tensor, piecewise polynomials of degree  $ \leq k$ and $\leq k+1$ are utilized for approximating the strain tensor and the velocity, respectively, and the heat and substance concentration are approximated by means of Lagrange finite elements of order $\leq k+1$.  Optimal a priori error estimates are derived and confirmed through some numerical examples that illustrate the performance of the proposed semi-augmented mixed-primal scheme.
	
\end{abstract}



\section{Introduction}

\noindent In nature and several technological applications, transport phenomena widely occur as a result  of a combined heat and mass transfer that are driven by buoyancy effects due to both temperature and concentration variations (see, e.g., \cite{Cheng79,Mojtabi2005,Nield99}). Such processes, also known as thermosolutal or double-diffusive natural convection, involving fluid circulation in a porous media, are frequently found in astrophysics, oceanology, metallurgy, electrophysics and geophysics, but also appear in several engineering applications such as filtration processes, geothermal energy exploitation, spreading on porous substrates, bio-film growth, gasification of biomass,  to name a few. 

\noindent From the mathematical point of view the Darcy-Oberbeck-Boussinesq model allows to adequately describe and quantify this complex flow by means of a nonlinear partial differential equations system. More precisely, the momentum  and conservation of fluid mass give rise to a  Navier-Stokes/Darcy type system for describing the fluid flow in the porous media which, in turn,  is coupled via buoyancy forces and  convective mass and heat transfer to a vector advection-diffusion equation for describing the substance concentration and the temperature, as a result of an energy and mass transfer balance (see, e.g., \cite{Mojtabi2005,Nield99}). 

\noindent Many computational techniques  have been developed so far in order  to  numerically solve and simulate this problem and related ones (see \cite{abedi,dimitri,bmr19,cibik,CGMR-2020,ob-2,ob-3,ob-1,moraga,ob-5}, and the references therein). Particularly,  the contributions  \cite{abedi,dimitri,moraga,ob-2,ob-3,ob-1,ob-5} deal with double-diffusive convection in a cavity, whereas  in  \cite{bmr19,cibik,CGMR-2020} the authors consider the phenomenon in a porous media.   

\noindent In \cite{moraga}, a finite volume method is proposed and applied to solve agro-food processes, whereas some methods based on finite elements for this problem are \cite{abedi,dimitri}. In \cite{abedi}, the authors proposed stabilized finite element formulations based on the SUPG (Streamline-Upwind/Petrov-Galerkin) and PSPG (Pressure-Stabilization/Petrov-Galerkin) methods to solve the problem in unsteady state. Numerical simulations in two and three dimensions illustrate the accuracy and performance of this technique. However, the theoretical analyses of the associate continuous and the discrete variational problems as well as the convergence of the method are not carried out there, and the method only allows to carry out low-order approximations of the main unknowns.  On the other hand, in \cite{dimitri} the problem is considered in steady state and analyzed by using the boundary control theory. The authors formulate and prove solvability results for the corresponding boundary control problem, state local uniqueness and stability of optimal solutions.

\noindent Focusing in double-diffusive viscous flow in porous media,  \cite{cibik} proposes a technique consisting of a projection-based stabilization method in the unsteady state. There, the convergence of the velocity, temperature and concentration  in the semi-discrete case are derived. In addition, some numerical experiments are reported to confirm optimal order error estimates  and to compare their results with previous ones. In \cite{bmr19} the authors construct  a divergence-conforming primal scheme and establish the  existence and uniqueness results for the continuous problem and the discrete scheme as well as convergence properties. On the other hand, in \cite{CGMR-2020} the authors propose a high-order fully-mixed method based on the introduction of  the velocity gradient, the temperature gradient and the concentration gradient  as new unknowns into the problem. The resulting formulation  has a saddle-point structure on reflexive Banach spaces for both the Navier-Stokes/Darcy and the thermal energy conservation equations. There, it is particularly shown that the discrete scheme is well-posed  and an  a priori error estimate for the Galerkin scheme is also derived under sufficiently small data. However, feasible  finite element subspaces must be constructed over meshes with a macro-element structure in order to satisfy an inf-sup compatibility condition, which in turn significantly increased the computational cost, especially in the three-dimensional case. 

\noindent According to the above discussion and in order to contribute with the design and analysis  of new mixed finite schemes for simulating double-diffusive convection in porous media.  The main goal of this work is precisely  to propose a new semi-augmented mixed finite element method in which the strain tensor and a pseudo-stress tensor are introduced as primary unknowns of interest in the fluid equations and the pressure is eliminated from the system by its own definition. To avoid any inf-sup restriction, to guarantee greater flexibility regarding finite element spaces and lower computational cost than \cite{CGMR-2020},  the respective variational formulation is  augmented by using appropriate redundant penalization Galerkin type-terms based on the constitutive and equilibrium equations combined with a primal formulation for the heat and substance concentration equations in standard Hilbert spaces. In this way, the aforementioned  strain tensor, the nonlinear pseudostress, the velocity and a vector field whose components are the temperature of the fluid and the substance concentration are the main unknowns of the resulting coupled system. Moreover, physical boundary conditions are considered. Indeed, a no-slip condition (that is, homogeneous Dirichlet condition) for the fluid velocity, a prescribed temperature and substance concentration on a Dirichlet boundary and no heat/mass flow across an isolated surface/homogeneous Neumann condition.

\noindent Concerning the solvability analysis, we proceed similarly to \cite{almonacid-2018,colm-2016} by introducing an equivalent fixed-point setting. According to this, combining  the Lax-Milgram theorem with the classical Banach and Brouwer fixed-point theorems, we establish the respective solvability of the continuous problem and the associated Galerkin scheme, under suitable regularity assumptions (see \cite{alvarez-2015}, for further details), a feasible choice of parameters and, in the discrete case, for any family of finite element subspaces. To handle the non-homogeneous Dirichlet condition for the temperature and concentration, we carry out a rigorous treatment of the boundary data throughout the analysis by means of appropriate extensions involving the Scott-Zhang interpolator (in the discrete case), which allows us to establish the well-posedness of our scheme, along with its convergence result and the respective a priori error bounds.
Up to the best of our knowledge, because of the difficulties that can arise in the analysis, the physically relevant non-homogeneous Dirichlet condition case is usually either omitted or not considered; this is what motivates us to contribute in this direction as well.  

\noindent A Strang-type lemma,  enables us to derive the corresponding C\'ea estimate and to provide optimal a priori error bounds for the Galerkin solution. In turn, the pressure can be recovered by a  post-processed  of the discrete solutions, preserving the same rate of convergence. Finally, some  numerical experiments are presented to illustrate the performance of the technique, including the unsteady case with unknown solution,  to confirm the expected orders. 

The contents of this paper is presented  as follows. At the end of this section, we introduce some standard notations to be used throughout the manuscript. In Section \ref{section2}, we introduce the model problem and the data. We also  derive  an  equivalent first-order equations in terms of additional variables. Then, in Section \ref{section3}, we derive the semi-augmented mixed-primal variational formulation and establish its well-posedness. The associated Galerkin scheme is introduced and analyzed in Section \ref{FEM-specific-spaces}.  In Section \ref{section5}, we derive the corresponding C\'ea estimate and, finally, in Section \ref{sec:Nu-Ex} we present some numerical examples illustrating the performance of our semi-augmented mixed-primal finite element method.  

\subsection{Notations}\label{Notations}

Let us denote by $\Omega\subseteq\mathrm{R}^n$, with $n\in\{2,3\}$,  a given bounded domain with polygonal/polyhedral boundary $\Gamma$  with outward unit normal vector  $\bn$ and
let $\Gamma_{\mathrm{D}}\,,\Gamma_{\mathrm{N}}\,\subseteq\,\Gamma$ be
such that $\Gamma_{\mathrm{D}}\,\cap\,\Gamma_{\mathrm{N}}\,=\,\varnothing$,
$|\Gamma_{\mathrm{D}}|\,>\,0 $ and $\Gamma\,=\,\overline{\Gamma}_{\mathrm{D}}\,\cup\,\overline{\Gamma}_{\mathrm{N}}$. Standard notation will be adopted for Lebesgue spaces 
 $\mathrm{L}^p(\Omega)$ with norm $\|\cdot\|_{0, p,\Omega}$ or  $\|\cdot\|_{0,\Omega}$ if $p=2$,  and Sobolev spaces  $\mathrm{H}^s(\Omega)$ with norm  $\|\cdot\|_{s,\Omega}$, and semi-norm $|\cdot|_{s,\Omega}$. 
 In particular, when $\mathrm{A}$ denotes a generic scalar functional space, then  we  will denote its vectorial and tensor counterparts by $\mathbf{A}$ and $\mathbb{A}$, respectively.  

 For vector fields $\bv=(v_i)_{1\leq i\leq n}$ and $\bw=(w_i)_{1\leq i\leq n}$,  we set  the gradient, divergence and dyadic product operators, as
\[\nabla\bv:=\left(\frac{\partial v_i}{\partial x_j}\right)_{1\leq i,j\leq n},\quad\mathrm{div}\,\bv:=\sum_{j=1}^n\frac{\partial v_j}{\partial x_j}\,,\quad\mathrm{and}\quad \bv\otimes\bw:=(v_iw_j)_{1\leq i,j\leq n}\,,\]
respectively. Furthermore, given the tensor fields $\btau=(\tau_{ij})_{1\leq i,j\leq n}$ and $\boldsymbol{\zeta}=(\zeta_{ij})_{1\leq i,j\leq n}$, 
we let $\mathbf{div}\,\btau$ be the divergence operator $\mathrm{div}$ acting along the rows of $\btau$, and define 
the transpose, the trace, the tensor inner product, and the deviatoric tensor, respectively, as
\[\btau^{\mathrm{t}}\ :=\ (\tau_{ji})_{1\leq i,j\leq n},\quad \tr(\btau):=\sum_{j=1}^n \tau_{ii},\quad\btau:\boldsymbol{\zeta}:=\sum_{j=1}^n\tau_{ij}\zeta_{ij}\,,\quad\mathrm{and}\quad\btau^{\mathrm{d}}:=\btau-\frac{1}{n}\tr(\btau)\mathbb{I},\]
where $\mathbb{I}$ stands for the identity tensor in $\mathrm{R}^{n\times n}$. 
We recall that the space 
$$\Hdiv\ :=\ \Big\{\btau\in\mathbb{L}^2(\Omega):\quad\mathbf{div}\,\btau\in\mathbf{L}^2(\Omega)\Big\}\,,$$
equipped with the usual norm $$\|\btau\|^2_{\mathbf{div},\Omega}:=\|\btau\|^2_{0,\Omega}+\|\mathbf{div}\,\btau\|^2_{0,\Omega},$$
is a Hilbert space. Finally, we employ $\boldsymbol{0}$ to denote a generic null vector and use $c$ or $C$, with or  without  subscripts,  
bars,  tildes  or  hats,  to  mean  generic  positive  constants  independent  of  the discretization parameters, which may take different values at different places.

\section{The model problem}\label{section2}

This section introduces the mathematical model we address in the present work as well as the auxiliary unknowns  that are introduced and considered in the subsequent variational formulation. Under the Oberbeck-Boussinesq approximation framework,  double-diffusive natural convection phenomenon  in a porous media is described in terms of a Navier-Stokes-Brinkman model coupled to a system of diffusion-advection equations. In the stationary state, the problem consists of finding the velocity  $\bu$, the pressure $p$, and the vector $\bvarphi=(\varphi_1,\varphi_2)^{\mathrm{t}}$ when $n=2$ and $\bvarphi=(\varphi_1,\varphi_2, 0)^{\mathrm{t}}$ when $n=3$,  where $\varphi_1$ and $\varphi_2$ are the temperature and the concentration fields, respectively, of a fluid in a confined porous  enclosure $\Omega$, satisfying the set of  equations:  
\begin{equation}\label{eqn:model}
\begin{array}{rcl}
\gamma\,\bu-2\,\mathbf{div}\left(\nu(\bvarphi)\mathbf{e}(\bu)\right)+\left(\bu\cdot \nabla \right)\bu+\nabla \,p & = & \left(\alp\cdot\bvarphi\right)\bg \qin \Omega\,,\\[1ex]
\mathrm{div}\,\bu & = & 0 \qin \Omega\,,\\[1ex]
-\bdiv(\mathbb{K}\nabla \bvarphi)\,+\,(\nabla\bvarphi)\bu & = & \mathbf{0} \qin \Omega\,,
\end{array}
\end{equation}
where $\mathbf{e}(\bu)$ stands for the symmetric part of the velocity gradient, i.e., $\mathbf{e}(\bu) := \frac{1}{2}(\nabla \bu+(\nabla \bu)^{\mathrm{t}})$\,. The data are the gravitational force $\boldsymbol{g}$, the positive constant $\gamma$ corresponding to the reciprocal of the Darcy number, the  (thermal and solute) expansion coefficients vector $\alp$,  and the diagonal  matrix of  (thermal and mass)  diffusion constants $\mathbb{K}=\mathrm{diag}(k_i)_{1\leq i \leq n}\in \mathbb{L}^{\infty}(\Omega)$, with $k_i=0$ when $i=3$,  which is assumed to be a uniformly positive definite tensor, which means that there exists a positive constant $k_0$ such that
\begin{equation}\label{eq:def-pos-k}
(\mathbb{K}\, \boldsymbol{x})\cdot \boldsymbol{x}\, \geq\, k_0\, |\boldsymbol{x}|^2\,\qquad  \forall \,\boldsymbol{x}\in\mathrm{R}^n\,, \quad \mbox{with}\quad k_0=\min\{k_1,k_2,\dots, k_n\}
\,.
\end{equation}
In turn, the kinematic viscosity $\nu: \mathrm R^+ \times \mathrm R^+ \longrightarrow \mathrm R^+$ is assumed to be a bounded and Lipschitz continuous function that might depend on both the temperature and the mass concentration. That is, we assume  the existence of positive constants $\nu_1, \,\nu_2,$ and $ L_\nu$ that satisfy
\begin{equation}\label{eqn:nu-properties}
\nu_1\,\leq \,\nu(\boldsymbol \phi)\,\leq\, \nu_2 
\qan
|\nu(\boldsymbol \phi)-\nu(\boldsymbol \psi)|\,\leq\, L_{\nu}\,|\boldsymbol \phi - \boldsymbol \psi|
\qquad \forall\, \boldsymbol \phi, \,  \boldsymbol \psi \in \mathrm R \times \mathrm R^+ \,.
\end{equation}

\noindent The system \eqref{eqn:model} is finally completed with a non-slip condition on the whole boundary  for the velocity and physical boundary conditions for both the temperature and the concentration fields, that is
\begin{equation}\label{eqn:bc}
\bu\,=\,\mathbf{0}\,\hspace{0.2cm}\mbox{on $\Gamma$\,,}
\qquad 
\disp
\bvarphi\,=\,\bvarphi_{\D}\,
\,\, 
\mbox{on $\GD$}
\qquad \mbox{and}\qquad
\disp
(\nabla\bvarphi)\,\bn =\mathbf{0}\,
\,\, \mbox{on}\,\, \GN\,,
\end{equation}
where $\bvarphi_\D\in \mathbf{H}^{1/2}(\Gamma_\D)$  is a known trace function on  $\GD.$ 

Next, proceeding similarly to \cite{caucao-2017}, we introduce the strain tensor $\bt:=\mathbf{e}(\bu)$ as an auxiliary variable and then define 
\begin{equation}\label{eqn:sigma}
\bsig\ :=\ 2\nu(\bvarphi)\bt-(\bu\otimes\bu)-\big(p+c_{\bu}\big)\mathbb{I}\quad\mathrm{in}\quad\Omega\,,
\end{equation}
as an additional tensorial unknown,  where $c_{\bu}$ is a constant to be suitably defined next (see equation \eqref{eqn:def-c0} below). Thus, noting that the incompressibility condition of the fluid  $(\div\, \bu=0)$ implies that $\bdiv (\bu\otimes\bu)=(\bu\cdot \nabla)\bu$, we get from the first relation of \eqref{eqn:model} the equilibrium equation 
\begin{equation*}\label{eqn}
\gamma\bu-\bdiv\,\bsig\ =\ (\alp\cdot\bvarphi)\,\bg\quad\mathrm{in}\quad\Omega\,,
\end{equation*}
and taking deviatoric part in \eqref{eqn:sigma}, we find that  the  constitutive relation defined by \eqref{eqn:sigma} can be written as
\begin{equation*}\label{eqn:sigma-d}
\bsig^{\mathrm{d}}\ =\ 2\nu(\bvarphi)\bt-(\bu\otimes\bu)^{\mathrm{d}}\,\qin\,\,\Omega\,
\end{equation*}
where $\bt^{\mathrm d}=\bt$ since $\div \,\bu=0.$ Thus, the pressure is eliminated from the system, however taking trace in \eqref{eqn:sigma}, we find that it can be easily recovered according to the formula
\begin{equation}\label{eqn:pressure}
p\,=\,-\dfrac{1}{n}\Big[\tr(\bsig+(\bu\otimes\bu))\Big]\,-\,c_{\bu}\,,\qquad \mbox{in} \qquad \Omega\,.
\end{equation} 
Now, since   $p\in\mathrm{L}^2_0(\Omega)$ (which as usual  is clearly required for uniqueness of an eventual pressure solution to \eqref{eqn:model}), the   equation \eqref{eqn:pressure} suggests to define
\begin{equation}\label{eqn:def-c0}
c_{\bu}\,=\,-\dfrac{1}{n |\Omega|}\int_\Omega  \tr(\bu\otimes\bu)\,=\,-\dfrac{1}{n |\Omega|}\|\bu\|_{0,\Omega}^2\,,
\end{equation} 
so as to get the equivalence
\begin{equation}\label{eqn:mean-trace-sig}
\int_\Omega p \,=\, 0\qquad \mbox{if, and only if,} \qquad \int_\Omega  \tr\,\bsi\,=0\,.
\end{equation} 

\noindent According to the above discusion, the system \eqref{eqn:model} and \eqref{eqn:bc}  equivalently reads: Find $\bt$, $\bsi$ and $\bvarphi$ such that
\begin{equation}\label{eqn:model-1}
\begin{array}{c}
\displaystyle 2\nu(\bvarphi)\bt - \bsig^{\mathrm{d}} - (\bu\otimes\bu)^{\mathrm{d}}\ =\ \mathbf{0}\qin \Omega\,,\quad
\displaystyle \gamma\bu-\bdiv\,\bsig\ =\ (\alp\cdot\bvarphi)\,\bg\qin \Omega\,,\\[2ex]
\bt+\boldsymbol{\omega}(\bu)\ =\ \nabla\bu\qin \Omega\,,\quad
-\bdiv(\mathbb{K}\nabla \bvarphi)\,+\,(\nabla\bvarphi)\bu\ =\ \mathbf{0}\qin \Omega\,,\\[2ex]
\displaystyle \bu\ =\ \boldsymbol{0}\qon \Gamma\,, \quad
\bvarphi\ =\ \bvarphi_{\D} \qon  \Gamma_{\mathrm{D}}\,,\quad
(\nabla\bvarphi)\,\bn\ =\  \mathbf{0} \qon  \Gamma_{\mathrm{N}}\,,\qan
\int_{\Omega} \tr\,\bsi\ =\ 0\,,
\end{array}
\end{equation}
where $\boldsymbol{\omega}(\bu):=\frac{1}{2}(\nabla \bu-(\nabla \bu)^{\mathrm{t}})$ is the skew-symmetric part of the velocity gradient. We emphasize that the introduction of the variables $\bt$ and $\bsi$ as new unknowns in the system allows us to equivalently rewrite the Navier-Stokes-Brinkman model (first row of \eqref{eqn:model}) in terms of a first-order set of equations. Also, observe that  according to \eqref{eqn:pressure}-\eqref{eqn:mean-trace-sig}, the zero mean value restriction of the pressure on the domain  is imposed via  the respective  restriction on $\tr\,\bsi$ in the last equation of \eqref{eqn:model-1}. Note further that the incompressibility condition of the fluid is implicitly present via the equilibrium relation defined by $\bsi$ according to the second equation in the first row of \eqref{eqn:model-1}. 

\section{The continuous formulation}\label{section3}

In this section we introduce and analyze the weak formulation proposed for  the system described by \eqref{eqn:model-1}. In Section  \ref{sec:augmented} we firstly deduce an augmented mixed variational formulation of \eqref{eqn:model-1} and then we rewrite it  as a fixed-point problem in Section \ref{sec:FPA}, whose analysis  is addressed through Sections \ref{sec:WD-FPA} and \ref{sec:WD-FPA2}.

\subsection{The semi-augmented mixed-primal variational formulation}\label{sec:augmented}

To begin with,  the fact that the trace of the  tensor  solution $\bsi$ of system  \eqref{eqn:model-1} has zero mean value in $\Omega$ (see last equation of \eqref{eqn:model-1}) suggests to introduce  the space
$$\mathbb{H}_0(\mathbf{div};\Omega)\ :=\ \left\{\btau\in \Hdiv:\quad\int_{\Omega}\tr\,\btau=0\right\}.$$
\noindent In  addition, by their own definitions, we introduce the following space for the strain tensor $\bt$, as 
$$\mathbb{L}_{\tr}^2(\Omega)\ :=\ \Big\{\br\in\mathbb{L}^2(\Omega):\quad\br^{\mathrm{t}}=\br\quad\mathrm{and}\quad \tr\,\br=0 \Big\}.$$

%
\noindent In turn, because of the boundary conditions for the temperature and concentrations (see second and third equations of last row in \eqref{eqn:model-1}) we consider the closed subspace 
$$\mathbf{H}^{1}_{\GD}(\Omega)\ :=\ \Big\{\bpsi\in\mathbf{H}^1(\Omega):\quad\bpsi\big|_{\GD}=\mathbf{0}\,\Big\}\,, $$
for which, from the generalized Poincar\'e inequality,  we know that there exists $c_{\text{gp}}>0$, depending only on $\Omega$ and $\GD$, such that
\begin{equation}\label{Poincare-inequality}
\|\bpsi\|_{1, \Omega}\;\leq \; c_{\text{gp}} |\bpsi|_{1, \Omega} \qquad \forall\ \bpsi\in \mathbf{H}^{1}_{\GD}(\Omega)\,.
\end{equation}

\noindent Now, testing the first equation from first row in \eqref{eqn:model-1} with $\br\in \mathbb{L}_{\tr}^2(\Omega)$, we obtain
\begin{equation}\label{trilinear term}
2\int_{\Omega}\nu(\bvarphi)\bt:\br-\int_{\Omega} \bsig^{\mathrm{d}}:\br-\int_{\Omega}(\bu\otimes\bu)^{\mathrm{d}}:\br\ =\ 0\qquad \forall\ \br\in \mathbb{L}_{\tr}^2(\Omega)\,.
\end{equation}

\noindent At this point, we readily note that in order to bound the third term on the left hand side of \eqref{trilinear term}, and thanks to the continuous injection $\boldsymbol{i}:\mathbf{H}^1(\Omega)\rightarrow \mathbf{L}^4(\Omega)$ (see e.g. \cite{adams-2003} or \cite{quartero-1994}), we requiere the unknown $\bu$ to live in $\mathbf{H}^{1}_{0}(\Omega):=\{\bv\in\mathbf{H}^{1}(\Omega):\; \bv|_{\Gamma}=\bf{0}\}$. Indeed, by applying Cauchy-Schwarz and H\"older inequalities, we deduce that there exists a  positive constant $c_1(\Omega):=\|\boldsymbol{i}\|^2$, such that
\begin{equation}\label{eqn:c2}
\left|\int_{\Omega}(\bu\otimes\boldsymbol{w})^{\mathrm{d}}:\br\right|\,\leq c_1(\Omega)\|\bu\|_{1,\Omega}\|\boldsymbol{w}\|_{1,\Omega}\|\br\|_{0,\Omega}\,\quad\forall\,\bu,\boldsymbol{w}\in\mathbf{H}^1_0(\Omega),\qquad  \forall\,\br\in\mathbb{L}^2(\Omega)\,.
\end{equation}

\noindent Next, multiplying the first equation from second row in \eqref{eqn:model-1} by a test function $\btau\in \mathbb{H}_0(\mathbf{div};\Omega)$, integrating by parts, using the Dirichlet condition for $\bu$, and the identity $\bt:\btau=\bt:\btau^{\mathrm{d}}$ (since $\bt$ is trace-free), we get 
$$\int_{\Omega}\bt:\btau^{\mathrm{d}}+\int_{\Omega}\boldsymbol{\omega}(\bu):\btau+\int_{\Omega}\bu\cdot\mathbf{div}\,\btau\ =\ 0\qquad \forall\,\btau\in \mathbb{H}_0(\mathbf{div};\Omega).$$

\noindent Likewise, the equilibrium relation associated to $\bsi$ (second equation from first row in \eqref{eqn:model-1}) is written as
$$\gamma\int_{\Omega}\bu\cdot\bv-\int_{\Omega}\bv\cdot\mathbf{div}\,\bsig\ =\ \int_{\Omega}(\alp\cdot\bvarphi)\bg\cdot\bv\qquad \forall\ \bv\in \mathbf{H}^1_0(\Omega)\,,$$
whereas, the symmetry of the pseudo-stress tensor is imposed in an ultra-weak sense (see e.g \cite{alvarez-2019}) through the identity
$$\displaystyle-\int_{\Omega}\bsig:\boldsymbol{\omega}(\bv)\ =\ 0\,\qquad \forall\,\bv\in\mathbf{H}^{1}_0(\Omega)\,.$$

As for the equation associated to the temperature  and concentration (second equation from second row in \eqref{eqn:model-1}), we simply multiply it  by a test function $\bpsi \in \mathbf{H}^1_{\GD}(\Omega)$ and, after integrating by parts, and employing the Neumann boundary contidions on $\GN$, we find
%
%
$$\int_{\Omega}\mathbb{K}\nabla\bvarphi: \nabla \bpsi\ +\ \int_{\Omega}(\nabla\bvarphi)\bu\cdot \bpsi\ =\ 0\qquad \forall\,\bpsi \in \mathbf{H}^{1}_{\GD}(\Omega)\,.$$

In this way, a preliminary weak formulation for the coupled problem (\ref{eqn:model-1}) takes the form: Find  $(\bt, \bsig, \bu, \bvarphi)\in\mathbb{L}_{\mathrm{tr}}^2(\Omega)\times\mathbb{H}_0(\mathbf{div};\Omega)\times \mathbf{H}^1_{0}(\Omega)\times \mathbf{H}^{1}(\Omega)$, with $\bvarphi\big|_{\GD}=\bvarphi_\D$, such that 
\begin{equation}\label{eqn:weak}
\begin{array}{rcll}
\displaystyle{2\int_{\Omega}\nu(\bvarphi)\bt:\br-\int_{\Omega} \bsig^{\mathrm{d}}:\br-\int_{\Omega}(\bu\otimes\bu)^{\mathrm{d}}:\br} & = & 0 &\forall\,\br\in \mathbb{L}_{\tr}^2(\Omega)\,,\\[3ex]
\displaystyle{\int_{\Omega}\bt:\btau^{\mathrm{d}}+\int_{\Omega}\boldsymbol{\omega}(\bu):\btau+\int_{\Omega}\bu\cdot\mathbf{div}\,\btau}& = & 0&\forall\,\btau\in \mathbb{H}_0(\mathbf{div};\Omega)\,,\\[3ex]
\displaystyle-\int_{\Omega}\bv\cdot\mathbf{div}\,\bsig-\int_{\Omega}\boldsymbol{\omega}(\bv): \bsig+\gamma\int_{\Omega}\bu\cdot\bv & = & \displaystyle\int_{\Omega}(\alp\cdot\bvarphi)\bg\cdot\bv\quad & \forall\,\bv\in \mathbf{H}^1_{0}(\Omega)\,,\\[3ex]
\displaystyle\int_{\Omega}\mathbb{K}\nabla\bvarphi: \nabla \bpsi + \int_{\Omega}(\nabla\bvarphi)\bu\cdot \bpsi & = & 0 & \forall\,\bpsi \in \mathbf{H}^{1}_{\GD}(\Omega)\,.
\end{array}
\end{equation}
%
%
In order to analyze  the variational formulation \eqref{eqn:weak}, and similarly as in \cite[Section 2]{caucao-2017} (see also \cite{almonacid-2018, alvarez-2019, colm-2017}), we additionally  augment \eqref{eqn:weak} by incorporating  the following residual Galerkin type-terms 
coming from the constitutive and equilibrium equations for the fluid, 
\begin{equation}\label{eqn:redundant-terms-1}
\begin{array}{rcll}
\disp
\kappa_1\int_{\Omega}(\mathbf{e}(\bu)-\bt):\mathbf{e}(\bv) & = & 0 & \forall\,\bv\in\mathbf{H}^1_0(\Omega)\,,\\[4ex]
\disp
-\kappa_2\gamma\int_{\Omega}\bu\cdot\mathbf{div}\,\btau+\kappa_2\int_{\Omega}\mathbf{div}\,\bsig\cdot\mathbf{div}\,\btau & = & \displaystyle -\kappa_2\int_{\Omega}(\alp \cdot \bvarphi)\bg\cdot\mathbf{div}\,\btau\quad & \forall\,\btau\in\mathbb{H}_0(\mathbf{div};\Omega)\,,\\[4ex]
\disp
\quad -\kappa_3\int_{\Omega}\left\{2\nu(\bvarphi)\bt - \bsig^{\mathrm{d}} - (\bu\otimes\bu)^{\mathrm{d}}\right\}:\btau^{\mathrm{d}}& = & 0 &\forall\,\btau\in\mathbb{H}_0(\mathbf{div};\Omega)\,,
\end{array}
\end{equation}
where $(\kappa_1,\kappa_2,\kappa_3)$ is a vector of positive parameters to be specified later in Section \ref{sec:WD-FPA}.  

Hence, letting 
$$\vbx\ :=\ (\bt,\bsig,\bu)\in \mathbb{H}:=\mathbb{L}_{\mathrm{tr}}^2(\Omega)\times \mathbb{H}_0(\mathbf{div};\Omega)\times\mathbf{H}^{1}_0(\Omega),$$
where $\mathbb{H}$ is endowed with the natural norm
$$\|\vby\|_{\mathbb{H}}:=\Big\{\|\br\|_{0,\Omega}^2+\|\btau\|_{\mathbf{div},\Omega}^2+\|\bv\|_{1,\Omega}^2\Big\}^{1/2},\quad\forall\,\vby:=(\br,\btau, \bv)\in \mathbb{H}\,,$$
and adding up \eqref{eqn:weak} with \eqref{eqn:redundant-terms-1}, we arrive at the following semi-augmented mixed-primal formulation for the double-diffusive natural convection problem in porous media: Find $(\vbx,\bvarphi)\in\mathbb{H}\times\mathbf{H}^1(\Omega)$ , with $\bvarphi\big|_{\GD}=\bvarphi_\D$,   such that
\begin{subequations}\label{augmented-mixed-primal}
\begin{eqnarray}
\label{augmented-mixed-primal-1}\bA_{\bvarphi}(\vbx,\vby)+\bB_{\bu}(\vbx,\vby) & = & F_{\bvarphi}(\vby)\qquad\quad \forall\,\vby\in\mathbb{H}, \\[2ex]
\label{augmented-mixed-primal-2}\bC(\bvarphi,\bpsi) & = & G_{\bu,\bvarphi}(\bpsi)\qquad \forall\,\bpsi\in\mathbf{H}^1_{\GD}(\Omega),
\end{eqnarray}
\end{subequations}
where, given $\bphi \in \mathbf{H}^{1}(\Omega)$ and $\bw\in\mathbf{H}_0^1(\Omega)$, $\bA_{\bphi}$, $\bB_{\bw}$ and $\bC$ are the bilinear forms defined, respectively, as
\begin{equation}\label{A}
\begin{array}{rcl}
\bA_{\bphi}(\vbx, \vby) & := & \displaystyle 2\int_{\Omega}\nu(\bphi)\bt:\big\{\,\br-\kappa_3\btau^{\mathrm{d}}\,\big\}\, + \,\int_{\Omega}\bt:\{\btau^{\mathrm{d}}-\kappa_1\mathbf{e}(\bv)\}\, - \,\int_{\Omega} \bsig^{\mathrm{d}}:\big\{\br-\kappa_3\btau^{\mathrm{d}}\big\}\\[3ex]
&&\displaystyle  +\ (1-\kappa_2\gamma)\int_\Omega\bu\cdot \bdiv\,\btau\, - \,\int_{\Omega}\bv\cdot\mathbf{div}\,\bsig\, +\, \int_{\Omega}\ww(\bu):\btau-\int_{\Omega} \bsig:\ww(\bv)\\[3ex]
&&\displaystyle  +\ \gamma \int_{\Omega}\bu\cdot\bv\, + \,\kappa_2\int_{\Omega}\mathbf{div}\,\bsig\cdot \mathbf{div}\,\btau\, + \,\kappa_1\int_{\Omega}\mathbf{e}(\bu):\mathbf{e}(\bv)\,,
\end{array}
\end{equation}
and
\begin{equation}\label{B}
\bB_{\bw}(\vbx, \vby)\ :=\ \int_{\Omega}(\bu\otimes\bw)^{\mathrm{d}}:\big\{\kappa_3\btau^{\mathrm{d}}-\br\big\}\,,
\end{equation}
for all $\vbx, \vby\in\mathbb{H}$ and
\begin{equation}\label{C}
\bC(\bchi,\bpsi)\ :\ =\,\int_{\Omega}\mathbb{K}\nabla \bchi : \nabla\bpsi\qquad \forall\,\bchi,\bpsi\in\mathbf{H}^1(\Omega)\,.
\end{equation}
In turn, $F_{\bphi}$ and $G_{\bw, \bphi}$ are linear functionals defined as
\begin{equation}\label{F}
F_{\bphi}(\vby)\ :=\ \int_{\Omega}(\alp\cdot \bphi)\bg \cdot\{\bv-\kappa_2\mathbf{div}\,\btau\big\}\,\qquad \forall\,\vby\in\mathbb{H}\,,
\end{equation}
and
\begin{equation}\label{G}
G_{\bw,\bphi}(\bpsi)\ :=\ -\int_{\Omega}(\nabla\bphi)\bw\cdot\bpsi \,\qquad \forall\,\bpsi\in\mathbf{H}^1_{\GD}(\Omega)\,.
\end{equation}

\subsection{The fixed point approach}\label{sec:FPA}

In this Section we describe a fixed point strategy that allow us to solve the coupled problem given by \eqref{augmented-mixed-primal}. We start by denoting $\mathbf{H}:=\mathbf{H}_{0}^{1}(\Omega)\times \mathbf{H}^{1}(\Omega)$, and defining the operator $\A:\mathbf{H} \rightarrow \mathbb{H}$ by
$$\A(\bw,\bphi)\,=\,(\A_{1}(\bw,\bphi),\, \A_{2}(\bw,\bphi),\, \A_{3}(\bw, \bphi))\,:=\,\vbx\quad \quad \forall\, (\bw, \bphi)\in \mathbf{H}$$
where $\vbx=(\bt, \bsig, \bu)\in \mathbb{H}$  is the unique solution of the augmented mixed formulation  given by \eqref{augmented-mixed-primal-1}, with $(\bw,\bphi)$ instead $(\bu,\bvarphi)$, that is: 
\begin{equation}\label{Defining_operator_A}
\begin{array}{rlll}
\bA_{\bphi}(\vbx, \vby)\, +\, \bB_{\bw}(\vbx, \vby) &=& F_{\bphi}(\vby)\, &\quad \forall\,\vby\in\mathbb{H}\,,
\end{array}
\end{equation}
where the bilinear forms $\bA_{\bphi}$, $\bB_{\bw}$ and the functional $F_{\bphi}$  are defined exactly as in \eqref{A}, \eqref{B}, and \eqref{F}, respectively.
In addition, we also introduce the operator $\B:\mathbf{H} \rightarrow \mathbf{H}^{1}(\Omega)$ defined as
$$\B(\bw,\bphi)\,:=\,\bvarphi \quad \quad \forall \,(\bw,\bphi)\in \mathbf{H},$$
where $\bvarphi$ is the unique solution of the problem \eqref{augmented-mixed-primal-2},
with $(\bw,\bphi)$ instead $(\bu,\bvarphi)$, that is: 
\begin{equation}\label{Defining_operator_B}
\bC(\bvarphi,\bpsi)\ =\ G_{\bw,\bphi}(\bpsi)\qquad  \forall\,\bpsi\in\mathbf{H}^1_{\GD}(\Omega),
\end{equation}
where the bilinear form $\bC$, and the functional $G_{\bw,\bphi}$ are defined by
\eqref{C} and \eqref{G}, respectively.
\medskip

In this way, by introducing the operator $\L: \mathbf{H} \rightarrow \mathbf{H}$ as
\begin{equation}\label{Defining_operator_L}
\L(\bw,\bphi)\,:=\,(\A_{3}(\bw, \bphi),\, \B(\A_{3}(\bw,\bphi), \bphi))\quad \quad \forall\, (\bw, \bphi)\in \mathbf{H}
\end{equation}
we readily realize that solving \eqref{augmented-mixed-primal} is equivalent to seeking a fixed point of $\L$, that is: Find
$(\bu,\bvarphi)\in \mathbf{H}$ such that
\begin{equation}\label{Fixed-point-equation}
\L(\bu,\bvarphi)\, =\, (\bu,\bvarphi)\,.
\end{equation}
The following two sections establish the well-posedness of \eqref{Fixed-point-equation}.

\subsection{Well-posedness of the uncoupled problems}\label{sec:WD-FPA}

 We begin by recalling the following lemmas which are useful to prove the ellipticity of the bilinear form $\bA_{\bphi}+\bB_{\bw}$.
\begin{lem}\label{decomp-1}
There exists $c_2(\Omega)>0$ such that
\begin{equation*}
c_2(\Omega)\|\btau_0\|^2_{0,\Omega}\ \leq\ \|\btau^{\mathrm{d}}\|^2_{0,\Omega}\, +\, \|\Div(\btau)\|^2_{0,\Omega}
\quad \quad \forall \,\btau = \btau_0 + c\,\mathbb{I}\in \mathbb{H}(\Div;\Omega)\,,
\end{equation*}
with $\btau_0\in\mathbb{H}_0(\Div;\Omega)$ and $c\in \mathrm R$.
\end{lem}

\begin{proof}
See \cite[Proposition 3.1]{brezzi-1991}.
\end{proof}

\begin{lem}\label{korn-estimate}
There holds
\begin{equation*}
\frac{1}{2}|\bv|_{1,\Omega}^{2}\;\leq \; \|\mathbf{e}(\bv)\|_{0,\Omega}^{2}\quad \quad \forall\, \bv\in \mathbf{H}^{1}_{0}(\Omega)\,.
\end{equation*}
\end{lem}

\begin{proof}
See \cite[Theorem 10.1]{McLean_2000}.
\end{proof}

We now provide sufficient conditions under which the uncoupled problems \eqref{Defining_operator_A} and \eqref{Defining_operator_B} are indeed uniquely solvable. 
\begin{lem}\label{Lemma_Flow}
Assume that $\kappa_1\in\big(0, 2\delta_3\big(2\nu_1-\frac{\kappa_3\nu_2}{\delta_1}\big)\big)$, $\kappa_2\in(0,2\delta_2)$, $\kappa_3\in\big(0,\frac{2\nu_1\delta_1}{\nu_2}\big)$ with $\delta_1\in\big(0,\frac{1}{\nu_2}\big)$, $\delta_2\in\big(0, \frac{2}{\gamma}\big)$, $\delta_3\in (0,1)$. Then, there exists $r_0>0$ such that for each $r\in (0,r_0)$, problem \eqref{Defining_operator_A} has a unique solution $\A(\bw, \bphi):=\vbx\in \mathbb{H}$, for each $(\bw, \bphi)\in \mathbf{H}$ with $\|\bw\|_{1, \Omega}\leq r$. Moreover, there exists  $\CA$, independent of $(\bw, \bphi)$, such that
\begin{equation}\label{Estimate-operatorA}
\|\A(\bw, \bphi)\|_{\mathbb{H}}\;=\;\|\vbx\|_{\mathbb{H}}\;\leq\; \CA\;|\alp|\,|\bg|\,\|\bphi\|_{0,\Omega}\quad \forall\, (\bw, \bphi)\in \mathbf{H}\,.
\end{equation}
\end{lem}

\begin{proof}
We begin by deriving the continuity of the bilinear forms $\bA_{\bphi}$ and $\bB_{\bw}$
(sf. \eqref{A} and \eqref{B}, respectively). Indeed, applying Cauchy-Schwarz's inequality,
the assumptions \eqref{eqn:nu-properties}, and the fact that $\|\mathbf{e}(\bv)\|_{0,\Omega} \leq |\bv|_{1,\Omega}$
and $\|\ww(\bv)\|_{0,\Omega} \leq |\bv|_{1,\Omega}\;\;\forall\, 
\bv\in\mathbf{H}^{1}(\Omega)$, we deduce that there exists a positive constant $C_{\bA_{\bphi}}$, depending on $\kappa_1, \kappa_2, \kappa_3, \nu_2, \gamma$, such that
\begin{equation}\label{Bouded-A_phi}
|\bA_{\bphi}(\vbx, \vby)|\;\leq\; C_{\bA_{\bphi}}\|\vbx\|_{\mathbb{H}}\|\vby\|_{\mathbb{H}}\qquad \forall\, \vbx, \vby\in \mathbb{H}\,,
\end{equation}
where $C_{\bA_{\bphi}} := 3\max\{2,\gamma+\kappa_1,\kappa_2+\kappa_3,2\nu_2,1+2\nu_2\kappa_3,1+|1-\gamma\kappa_2|\}$.
In turn, by applying H\"{o}lder's inequality and the continuous injection $\mathbf{i}:\mathbf{H}^{1}(\Omega)\to \mathbf{L}^{4}(\Omega)$, with constant $c_1(\Omega)$, we obtain that
\begin{equation}\label{Bounded-B_w}
|\bB_{\bw}(\vbx, \vby)|\;\leq\;c_1(\Omega)(2+\kappa_3^2)^{1/2}\|\bw\|_{1, \Omega}\,\|\vbx\|_{\mathbb{H}}\,\|\vby\|_{\mathbb{H}}\qquad \forall\, \vbx, \vby \in \mathbb{H}\,.
\end{equation}
Then, it follows from \eqref{Bouded-A_phi} and \eqref{Bounded-B_w}, that there exists a positive constant denote by $\|\bA_{\bvarphi}+\bB_{\bw}\|$, that depends on $\kappa_1, \kappa_2, \kappa_3, \nu_2, \gamma, c_1(\Omega)$ and $\|\bw\|_{1,\Omega}$, such that
\begin{equation}\label{Bounded-A_phi+B_w}
|(\bA_{\bphi}+\bB_{\bw})(\vbx,\vby)|\;\leq\;\|\bA_{\bphi}+\bB_{\bw}\|\,\|\vbx\|_{\mathbb{H}}\,\|\vby\|_{\mathbb{H}}\qquad \forall\, \vbx, \vby \in \mathbb{H}\,.
\end{equation}
We now address the proof of the ellipticity for the bilinear form $\bA_{\bphi}+\bB_{\bw}$.
To this end, we first derive this property for the bilinear form $\bA_{\bphi}$. In fact,
from \eqref{A}, by applying \eqref{eqn:nu-properties} and the Cauchy-Schwarz inequality,
together with Lemma \ref{korn-estimate}, it followss that
\begin{eqnarray*}
\bA_{\bphi}(\vbx, \vbx) & = & 2\int_{\Omega}\nu(\bphi)\bt:\bt\, -\, 2\kappa_3\int_{\Omega}\nu(\bphi)\bt:\bsig^{\mathrm{d}}\, -\, \kappa_1\int_{\Omega}\bt:\mathbf{e}(\bu)\, +\, \kappa_3\|\bsig^{\mathrm{d}}\|_{0,\Omega}^2\\
&& -\ \kappa_2\gamma\int_\Omega\bu\cdot \bdiv\,\bsig\, +\, \gamma\|\bu\|_{0,\Omega}^2\, + \,\kappa_2\|\mathbf{div}\,\bsig\|_{0,\Omega}^2\, +\, \kappa_1\|\mathbf{e}(\bu)\|_{0,\Omega}^2\\[1ex]
& \geq & 2\nu_1\|\bt\|_{0,\Omega}^2\, -\, 2\kappa_3\nu_2\|\bt\|_{0,\Omega}\|\bsig^{\mathrm{d}}\|_{0,\Omega}\, -\, \kappa_1\|\bt\|_{0,\Omega}\|\mathbf{e}(\bu)\|_{0,\Omega}\, +\, \kappa_3\|\bsig^{\mathrm{d}}\|_{0,\Omega}^2\\
&& -\ \kappa_2\gamma\|\bu\|_{0,\Omega}\|\bdiv\,\bsig\|_{0,\Omega}\, +\, \gamma\|\bu\|_{0,\Omega}^2\, + \,\kappa_2\|\mathbf{div}\,\bsig\|_{0,\Omega}^2\, +\, \kappa_1\|\mathbf{e}(\bu)\|_{0,\Omega}^2\\[1ex]
& \geq & 2\nu_1\|\bt\|_{0,\Omega}^2\, -\, 2\kappa_3\nu_2\|\bt\|_{0,\Omega}\|\bsig^{\mathrm{d}}\|_{0,\Omega}\, -\, \kappa_1\|\bt\|_{0,\Omega}|\bu|_{1,\Omega}\, +\, \kappa_3\|\bsig^{\mathrm{d}}\|_{0,\Omega}^2\\
&& -\ \kappa_2\gamma\|\bu\|_{0,\Omega}\|\bdiv\,\bsig\|_{0,\Omega}\, +\, \gamma\|\bu\|_{0,\Omega}^2\, + \,\kappa_2\|\mathbf{div}\,\bsig\|_{0,\Omega}^2\, +\, \frac{\kappa_1}{2}|\bu|_{1,\Omega}^2
\end{eqnarray*}
Next, employing the Young inequality and gathering similar terms, we obtain
\begin{eqnarray*}
\bA_{\bphi}(\vbx, \vbx)
& \geq & \left(2\nu_1 -\frac{\kappa_3\nu_2}{\delta_1} - \frac{\kappa_1}{2\delta_3}\right)\|\bt\|_{0,\Omega}^2\, +\,
\kappa_3\big(1-\nu_2\delta_1\big)\|\bsig^{\mathrm{d}}\|_{0,\Omega}^2\, +\,  
\kappa_2\left(1 - \frac{\gamma\delta_2}{2}\right)\|\bdiv\,\bsig\|_{0,\Omega}^2\\
& & +\ \gamma\left(1 -\frac{\kappa_2}{2\delta_2}\right)\|\bu\|_{0,\Omega}^2\, +\, 
\frac{\kappa_1}{2}\big(1 - \delta_3\big)|\bu|_{1,\Omega}^2
\end{eqnarray*}
and then, choosing $\kappa_1$, $\kappa_2$, $\kappa_3$, $\delta_1$, $\delta_2$ and $\delta_3$
in the ranges specified of the statement of the present Lemma, we deduce that there exists
a positive constant $\alpha(\Omega)$, independent of $(\bw, \bvarphi)$, such that
\begin{equation}\label{Ellipticity_A_phi}
\bA_{\bphi}(\vbx, \vbx)\ \geq\ \alpha(\Omega)\,\|\vbx\|_{\mathbb{H}}^{2}\qquad \forall\, \vbx\in \mathbb{H}\,.
\end{equation}
Next, by combining \eqref{Bounded-B_w} and \eqref{Ellipticity_A_phi}, we find that
\begin{equation*}
(\bA_{\bphi}+\bB_{\bw})(\vbx, \vbx)\ \geq\ \Big\{\alpha(\Omega)-c_1(\Omega)(2+\kappa_3^{2})^{1/2}\|\bw\|_{1,\Omega}\Big\}\|\vbx\|_{\mathbb{H}}^2\qquad \forall\, \vbx\in \mathbb{H}\,,
\end{equation*}
from which, we deduce that
\begin{equation}\label{Ellipticity-A_phi+B_w}
(\bA_{\bphi}+\bB_{\bw})(\vbx, \vbx)\;\geq \; \frac{\alpha(\Omega)}{2}\|\vbx\|_{\mathbb{H}}^2\qquad  \forall\, \vbx\in \mathbb{H}\,,
\end{equation}
provided that $\|\bw\|_{1,\Omega} \leq r_0$, with 
\begin{equation}\label{radius}
r_0\ :=\ \frac{\alpha(\Omega)}{2c_1(\Omega)(2+\kappa_3^2)^{1/2}}\,,
\end{equation}
which confirms the ellipticity of $\bA_{\bphi}+\bB_{\bw}$. On the other hand, by applying
the Cauchy-Schwarz inequality, we deduce
that $F_{\bphi}\in \mathbb{H}^{\prime}$ (cf. \eqref{F}) with
\begin{equation}\label{Bounded-F_phi}
\|F_{\bphi}\|\ \leq\ \sqrt{2}\,(2+\kappa_{2}^2)^{1/2}\,|\alp|\,|\bg|\,\|\bphi\|_{0,\Omega}\,.
\end{equation}
Consequently, a straightforward application of the Lax-Milgram lemma implies that there exists
a unique solution $\vbx\in \mathbb{H}$ of \eqref{Defining_operator_A}. Finally, from \eqref{Ellipticity-A_phi+B_w}
and \eqref{Bounded-F_phi}, and performing simple algebraic manipulations, we derive \eqref{Estimate-operatorA},
with $\displaystyle \CA:=\frac{2\,\sqrt{2}\,(2+\kappa_2^2)^{1/2}}{\alpha(\Omega)}>0$, independent of $(\bw, \bphi)$.
\end{proof}

\begin{rem}
At this point, we remark that for computational purposes, the constant $\alpha(\Omega)$
defined in Lemma {\rm \ref{Lemma_Flow}}, can be maximized by choosing the parameters
$\delta_1$, $\delta_2$, $\delta_3$, $\kappa_1$, $\kappa_2$ and $\kappa_3$ at the middle
points of their feasible ranges. Thus, adequate choices for these parameters are 
$\delta_1 := \frac{1}{2\nu_2}$, $\delta_2 := \frac{1}{\gamma}$ and $\delta_3 := \frac{1}{2}$,
which establish that
\begin{equation}\label{Parameters-maximized}
\kappa_1\ =\ \frac{\nu_1}{2}\,,\qquad \kappa_2\ =\ \frac{1}{\gamma}\qquad\text{and}\qquad \kappa_3\ =\ \frac{\nu_1}{2\nu_{2}^2}\,.
\end{equation}
\end{rem}

\begin{rem}\label{Equinorm-property-trace}
In order to establish the solvability of the problem \eqref{Defining_operator_B}, associated
to the operator $\B$, we first point out that $\bvarphi_{\D} \in \mathbf{H}^{1/2}(\GD)$ can be
continuously extended in the trace sense to $\mathbf{H}^{1/2}(\Gamma)$. Indeed, the set 
$$\mathcal{H}(\bvarphi_{\D})\ =\ (\gamma_0|_{\GD})^{-1}(\{\bvarphi_{\D}\})\ =\ \big\{\bvarphi\in \mathbf{H}^{1}(\Omega):\quad  \gamma_{0}(\bvarphi)|_{\GD}\;=\;\bvarphi_{\D}\big\}$$
is a closed and convex since the usual trace operator $\gamma_{0}:\mathbf{H}^{1}(\Omega)\to \mathbf{H}^{1/2}(\Gamma)$ is linear and continuous. Then, from The Best Approximation Theorem \cite[Theorem 7]{Daya Reddy} there exists a unique $\widetilde{\bvarphi}:=E(\bvarphi_{\D})\in \mathbf{H}^{1}(\Omega)$, where $E$ denotes the extension of $\bvarphi_{\D}$, such that $\gamma_{0}(\widetilde{\bvarphi})|_{\GD}\;=\;\bvarphi_{\mathrm{D}}$ with
\begin{equation}\label{equinorm}
\|\widetilde{\bvarphi}\|_{1,\Omega}\ =\ \inf_{\bvarphi\in \mathcal{H}(\bvarphi_{\D})}\|\bvarphi\|_{1, \Omega}\ =\ \text{\rm dist}({\bf 0},\mathcal{H}(\bvarphi_{\D}))\ =\ \|\bvarphi_{\D}\|_{1/2,\GD}.
\end{equation}
In this way, the suitable extension of $\bvarphi_{\D}\in \mathbf{H}^{1/2}(\GD)$ is not but the corresponding one of $\gamma_0(\widetilde{\bvarphi})\in \mathbf{H}^{1/2}(\Gamma)$ to $\mathbf{H}^{1}(\Omega)$.
\end{rem}

\begin{lem}\label{Lemma-Transport}
For each $(\bw, \bphi)\in \mathbf{H}$, problem \eqref{Defining_operator_B} has a unique solution $\bvarphi\in \mathbf{H}^{1}(\Omega)$, with $\bvarphi|_{\GD}=\bvarphi_{\D}$. Moreover, there exists a constant $\CB>0$ independent of $(\bw, \bphi)$, such that 
\begin{equation}\label{Estimate-operatorB}
\|\B(\bw, \bphi)\|_{1, \Omega}\;=\;\|\bvarphi\|_{1, \Omega}\;\leq\; \CB\big\{\|\bw\|_{1,\Omega}\|\bphi\|_{1,\Omega}\, +\, \|\bvarphi_{\D}\|_{1/2,\GD}\big\}\,.
\end{equation}
\end{lem}

\begin{proof} 
We begin by nothing, according to the aforementioned in Remark \ref{Equinorm-property-trace},
that given $\bvarphi_{\D}\in \mathbf{H}^{1/2}(\GD)$ there exists $\bvarphi_1\in \mathbf{H}^{1}(\Omega)$
such that
\begin{equation}\label{eqn:wp-C}
\bvarphi_1|_{\GD}\ =\ \bvarphi_{\D}\qquad\text{and}\qquad \|\bvarphi_1\|_{1,\Omega}\ =\ \|\bvarphi_{\D}\|_{1/2,\GD}\,.
\end{equation}
Then, we consider the auxiliary linear problem: Find $\bvarphi_0\in\mathbf{H}^{1}_{\GD}(\Omega)$
such that
\begin{equation}\label{alternative_problem}
\bC(\bvarphi_0, \bpsi)\ =\ G_{\bw,\bphi}(\bpsi)\, -\, \bC(\bvarphi_1, \bpsi)\qquad \forall\, \bpsi\in \mathbf{H}^{1}_{\GD}(\Omega)
\end{equation}
where $\bC$ and $G_{\bw,\bphi}$ are defined in \eqref{C} and \eqref{G}, respectively.
In addition, from \eqref{C} and the Cauchy-Schwarz's inequality, we deduce that
\begin{equation}\label{Bounded-C}
|\bC(\bchi, \bpsi)|\ \leq\ \|\mathbb{K}\|_{\infty,\Omega}\|\bchi\|_{1, \Omega}\|\bpsi\|_{1, \Omega}\qquad \forall\,\bchi, \bpsi \in \mathbf{H}^1(\Omega)\,,
\end{equation}
which, in particular, confirms the boundedness of $\bC$. Then, from \eqref{C}, using
\eqref{eq:def-pos-k} and the Poncair\'e inequality \eqref{Poincare-inequality}, we get
\begin{equation}\label{Ellipticity-C-H1}
\bC (\bpsi, \bpsi)\ \geq\ \widetilde{\alpha}\,\|\bpsi\|_{1, \Omega}^{2} \qquad \forall\,\bpsi \in \mathbf{H}^{1}_{\GD}(\Omega)\,,
\end{equation}
which proves that $\bC$ is $\mathbf{H}^{1}_{\GD}(\Omega)$-elliptic with constant
\begin{equation}\label{eq:elip-const-trans}
\widetilde{\alpha}\ :=\ k_0\,c_{\text{gp}}^{-1}\,.
\end{equation}
Next, applying the Cauchy-Schwarz's inequality, the boundedness of the continuous injection
$\boldsymbol{i}:\mathbf{H}^{1}(\Omega)\to \mathbf{L}^{4}(\Omega)$, with constant $c_1(\Omega):=\|\boldsymbol{i}\|^2$,
\eqref{Bounded-C}, and the second identity in \eqref{eqn:wp-C}, we easily deduce that
\begin{equation}\label{Bounded-F}
|G_{\bw,\bphi}(\bpsi)-\bC(\bvarphi_1, \bpsi)|\ \leq\ \Big\{c_1(\Omega)\|\bw\|_{1,\Omega}\|\bphi\|_{1,\Omega}+\|\mathbb{K}\|_{\infty, \Omega}\|\bvarphi_{\mathrm{D}}\|_{1/2,\GD}\Big\}\|\bpsi\|_{1, \Omega}\quad \forall\, \bpsi \in \mathbf{H}^{1}_{\GD}(\Omega)\,,
\end{equation}
which establishes the boundedness of the right-hand side of \eqref{alternative_problem}.
Consequently, a direct application of the Lax-Milgram lemma implies that there exists a unique $\bvarphi_{0}\in \mathbf{H}^{1}_{\GD}(\Omega)$ that satisfies \eqref{alternative_problem} with 
$$\|\bvarphi_{0}\|_{1,\Omega}\ \leq\ \frac{1}{\widetilde{\alpha}}\,\big\{c_1(\Omega)\|\bw\|_{1,\Omega}\|\bphi\|_{1,\Omega}\, +\, \|\mathbb{K}\|_{\infty, \Omega}\|\bvarphi_{\mathrm{D}}\|_{1/2,\GD}\big\}\,.$$
On the other hand, we now set $\bvarphi := \bvarphi_0 + \bvarphi_1$, which satisfies that
$\bvarphi|_{\GD} = \bvarphi_0|_{\GD} + \bvarphi_1|_{\GD} = \bvarphi_{\D}$.
In addition, it is easy to check that $\bvarphi$ is a solution of problem \eqref{Defining_operator_B},
where the uniqueness comes from \eqref{Ellipticity-C-H1}.
Finally, $\bvarphi$ verifies the estimate \eqref{Estimate-operatorB}, with
constant $\CB := \max\{\widetilde{\alpha}^{-1}c_1(\Omega), \widetilde{\alpha}^{-1}\|\mathbb{K}\|_{\infty, \Omega}+1\}$. 
\end{proof}

\medskip
We end this section by introducing a further regularity hypotheses on the problem definig $\A$, which will be employed to derive a Lipschitz-continuity property for the operator $\L$. More precisely, we assume that for each $(\bw, \bphi)\in \mathbf{H}$ with $\|\bw\|_{1, \Omega}+\|\bphi\|_{1, \Omega}\leq r$, $r>0$ given, there holds $(\br, \btau, \bv)=\A(\bw, \bphi)\in \mathbb{L}^{2}_{\tr}(\Omega)\cap \mathbb{H}^{\varepsilon}(\Omega)\times \mathbb{H}_{0}(\bdiv;\Omega)\cap \mathbb{H}^{\varepsilon}(\Omega)\times \mathbf{H}^{1+\varepsilon}(\Omega)$,  for some $\varepsilon\in (0,1)$ (when $n=2$) or $\varepsilon\in [\frac{1}{2}, 1)$ (when $n=3$), with 
\begin{equation}\label{Regularity-estimate}
\|\br\|_{\varepsilon, \Omega}\, +\, \|\btau\|_{\varepsilon, \Omega}\, +\, \|\bv\|_{\varepsilon, \Omega}\ \leq\ \widehat{C}(r)\,|\alp|\,|\bg|\,\|\bphi\|_{0,\Omega}\,,
\end{equation}
where $\widehat{C}(r)$ is a positive constant independent of $(\bw, \bphi)$ but depending on the upper bound $r$ of its norm. The reason of the estipulated ranges for $\varepsilon$ will be clarified  in the forthcomming analysis (specifically in the proofs of Lemmas \ref{Lipchtz-continuity-A} and \ref{Lipchtz-continuity-L} below). Also, we pay attention to the fact the while the estimate \eqref{Regularity-estimate} will be employed only to bound $\|\br\|_{\varepsilon, \Omega}$, we have stated it including the terms $\|\btau\|_{\varepsilon, \Omega}$ and $\|\bv\|_{1+\varepsilon, \Omega}$ as well, since due to the first and fourth equations of \eqref{eqn:model-1}, the regularities of $\br, \btau$ and $\bv$ will most likely be connected.

\subsection{Solvability analysis of the fixed point equation}\label{sec:WD-FPA2}

We begin by emphasizing that the well-posedness of the uncoupled problems \eqref{Defining_operator_A} and \eqref{Defining_operator_B} confirms that the operators $\A$, $\B$, and $\L$ (cf. \eqref{Defining_operator_L}) are well defined, and hence now we can address the solvability analysis of the fixed point problem presented in \eqref{Fixed-point-equation}. To this end, we will verify below the hypotheses of the Banach fixed-point theorem.

\begin{lem}\label{Ball-mapping}
Given $r\in (0,r_0)$, with $r_0$ given by \eqref{radius}, we let
$$W_r\ :=\ \big\{(\bw, \bphi)\in \mathbf{H}:\quad \|(\bw, \bphi)\|_{\mathbf{H}}\, \leq\, r\big\}\,,$$
and assume that the data satisfy
\begin{equation}\label{Small-data-1}
c(r)\,|\alp|\,|\bg|\, +\, \CB\,\|\bvarphi_{\D}\|_{1/2, \GD}\ \leq\ r
\end{equation}
where $c(r) := r\CA(1+r\CB)$, with $\CA$ and $\CB$ given by \eqref{Estimate-operatorA} and
\eqref{Estimate-operatorB}, respectively. Then, there holds $\L(W_r)\subseteq W_r$.
\end{lem}

\begin{proof}
It follows similar as in \cite[Lemma 3.5]{colm-2016-2}
\end{proof}

\begin{lem}\label{Lipchtz-continuity-A}
Let $r\in (0, r_0)$ with $r_0$ given by \eqref{radius}. Then, there exists a constant $\LA>0$, depending on the stabilization parameters $\kappa_2, \kappa_3$, and the constants $L_{\nu}, c_1(\Omega), \alpha(\Omega), C_{\varepsilon}$ (cf. \eqref{eqn:nu-properties}, \eqref{eqn:c2}, \eqref{Ellipticity_A_phi} and \eqref{epsilon-estimate}, respectively), such that for all
$(\bw, \bphi), (\widetilde{\bw}, \widetilde{\bphi})\in \mathbf{H}$,
with $\|\bw\|_{1,\Omega}, \|\widetilde{\bw}\|_{1,\Omega}\,\leq \;r$, there holds
\begin{equation}\label{LC-Estimate-A}
\begin{array}{l}
\|\A(\bw, \bphi)-\A(\widetilde{\bw}, \widetilde{\bphi})\|_{\mathbb{H}}\ \leq\
\LA\Big\{\|\A_3(\bw, \bphi)\|_{1, \Omega}\,\|\bw-\widetilde{\bw}\|_{1,\Omega}\\[1ex]
\qquad\qquad\qquad\qquad +\ \|\A_1(\bw, \bphi)\|_{\varepsilon, \Omega}\,\|\bphi-\widetilde{\bphi}\|_{0,n/\varepsilon,\Omega}\ +\ |\alp|\,|\bg|\,\|\bphi-\widetilde{\bphi}\|_{1, \Omega}\Big\}\,.
\end{array}
\end{equation}
\end{lem}

\begin{proof}
Given $(\bw, \bphi), (\widetilde{\bw}, \widetilde{\bphi})$ as stated, we let $\vbx = (\bt, \bsig, \bu)=\A(\bw, \bphi)$ and  $\vbz = (\widetilde{\bt}, \widetilde{\bsig}, \widetilde{\bu})=\A(\widetilde{\bw}, \widetilde{\bphi})$, which according to the definition of operator $\A$ (cf. \eqref{Defining_operator_A}), means that:
\begin{equation*}
\bA_{\bphi}(\vbx, \vby)+\bB_{\bw}(\vbx, \vby)\ =\ F_{\bphi}(\vby)\qquad\text{and}\qquad \bA_{\widetilde{\bphi}}(\vbz, \vby)+\bB_{\widetilde{\bw}}(\vbz, \vby)\ =\ F_{\widetilde{\bphi}}(\vby)\qquad \forall\,\vby\in\mathbb{H}\,.
\end{equation*}
Then, subtracting both identities, replacing $\vbz = \vbz-\vbx+\vbx$,
and using the bilinearity of $\bA_{\bphi}+\bB_{\bw}$ for any $\bphi$ and $\bw$,
it follows from \eqref{Defining_operator_A} that:
\begin{equation}\label{Bilinearity-Porperty-A_phi+B_w}
\displaystyle (\bA_{\widetilde{\bphi}}+\bB_{\widetilde{\bw}})(\vbx-\vbz, \vby)\ =\ (F_{\bphi}-F_{\widetilde{\bphi}})(\vby)\, +\, (\bA_{\widetilde{\bphi}}-\bA_{\bphi})(\vbx, \vby)\, +\, \bB_{\widetilde{\bw}-\bw}(\vbx,\vby)\qquad\forall\,\vby\in\mathbb{H}\,.
\end{equation}
Moreover, applying the ellipticity of $\bA_{\widetilde{\bphi}}+\bB_{\widetilde{\bw}}$
(cf. \eqref{Ellipticity-A_phi+B_w}), and then employing \eqref{Bilinearity-Porperty-A_phi+B_w}
with $\vby := \vbx-\vbz$, we find that
\begin{eqnarray}
\frac{\alpha(\Omega)}{2}\|\vbx-\vbz\|_{\mathbb{H}}^2 & \leq & (\bA_{\widetilde{\bphi}}+\bB_{\widetilde{\bw}})(\vbx-\vbz, \vbx-\vbz)\nonumber\\[2ex]
& = & (F_{\bphi}-F_{\widetilde{\bphi}})(\vbx-\vbz)\, +\, (\bA_{\widetilde{\bphi}}-\bA_{\bphi})(\vbx, \vbx- \vbz)\, +\, \bB_{\widetilde{\bw}-\bw}(\vbx, \vbx-\vbz)\,.\label{Estimate-from-ellipticity}
\end{eqnarray}
Then, for the first and third terms on the right-hand side of \eqref{Estimate-from-ellipticity},
we employ the Cauchy-Schwarz and H\"{o}lder inequalities, together with the continuous injection
$\mathbf{i}:\mathbf{H}^{1}(\Omega)\to \mathbf{L}^{4}(\Omega)$, similar as in \eqref{Bounded-F_phi}
and \eqref{Bounded-B_w}, in order to obtain that
\begin{equation}\label{Estimate_F_phi1-F_phi2}
\begin{array}{rcl}
|(F_{\bphi}-F_{\widetilde{\bphi}})(\vbx-\vbz)| & = & \displaystyle\Bigg|\int_{\Omega}(\alp (\bphi-\widetilde{\bphi}))\bg\cdot\big\{(\bu-\widetilde{\bu})-\kappa_2\,\bdiv(\bsig-\widetilde{\bsig})\big\}\Bigg|\\[3ex]
& \leq & \sqrt{2}\,(2+\kappa_{2}^{2})^{1/2}\,|\alp|\,|\bg|\,\|\bphi-\widetilde{\bphi}\|_{1,\Omega}\,\|\vbx-\vbz\|_{\mathbb{H}}
\end{array}
\end{equation}
and
\begin{equation}\label{B_w2-B_w1}
\begin{array}{rcl}
|\bB_{\bw-\widetilde{\bw}}(\vbx, \vbx-\vbz)| & = & \displaystyle \Bigg|\int_{\Omega}(\bu \otimes (\widetilde{\bw}-\bw))^{\tt{d}}:\big\{\kappa_3(\bsig-\widetilde{\bsig})^{\tt{d}}-(\bt-\widetilde{\bt})\big\}\Bigg|\\[3ex]
& \leq & c_1(\Omega)\,(2+\kappa_3^2)^{1/2}\,\|\bu\|_{1, \Omega}\,\|\bw-\widetilde{\bw}\|_{1, \Omega}\,\|\vbx-\vbz\|_{\mathbb{H}}
\end{array}
\end{equation}
On the other hand, for the second term in the right-hand side of \eqref{Estimate-from-ellipticity},
we apply the Lipschitz continuity property for $\nu$ given in \eqref{eqn:nu-properties}, the Cauchy-Schwarz and H\"{o}lder inequalities,
and the definition of $\bA_{\bvarphi}$ (cf. \eqref{A}), to obtain that
\begin{equation}\label{A_phi2-A_phi_1}
\begin{array}{l}
\displaystyle |(\bA_{\widetilde{\bphi}}-\bA_{\bphi})(\vbx, \vbx- \vbz)|\;=\;\Bigg|\int_{\Omega}2(\nu(\bphi)-\nu(\widetilde{\bphi}))\bt:\{(\bt-\widetilde{\bt})-\kappa_3(\bsig-\widetilde{\bsig})^{\tt{d}}\}\Bigg|\\[4ex]
\hspace{1cm} \quad \quad\;\leq \; 2L_{\nu}(2+\kappa_3^2)^{1/2}\|\bphi-\widetilde{\bphi}\|_{0,2q,\Omega}\|\bt\|_{0,2p,\Omega}\;\|\vbx-\vbz\|_{\mathbb{H}},
\end{array}
\end{equation}
with $p,q\in (1, +\infty)$ such that $\displaystyle \frac{1}{p}+\frac{1}{q}=1.$ We now proceed as in \cite[Lemma 3.9]{alvarez-2015}. In fact, given the further $\varepsilon$-regularity assumed in \eqref{Regularity-estimate}, we reall that the Sobolev embedding theorem (see e.g \cite[Theorem 4.12]{adams-2003}) establishes the continuous injection $i_{\varepsilon}:\mathrm{H}^{\varepsilon}(\Omega)\to \mathrm{L}^{2p}(\Omega)$ with boundedness constant $C_{\varepsilon}$, where
\begin{equation*}
2p\ =\ \left\{\begin{array}{cl}
\displaystyle \frac{2}{1-\varepsilon}  & \text{if } n=2\\[2ex]
\displaystyle \frac{6}{3-2\varepsilon} & \text{if } n=3
\end{array}\right.
\end{equation*}
and $2q=\frac{n}{\varepsilon}$, and therefore, there holds
\begin{equation}\label{epsilon-estimate}
\|\bt\|_{0,2p,\Omega}\;\leq \; C_{\varepsilon}\, \|\bt\|_{\varepsilon, \Omega}\qquad \forall\,\bt\in \mathbb{H}^{\varepsilon}(\Omega)\,.
\end{equation}
In this way, denoting
$$\LA\ :=\ \frac{2}{\alpha(\Omega)}\max\Big\{c_1(\Omega)\,(2+\kappa_3^2)^{1/2},\ 2\,L_{\nu}(2+\kappa_3^2)^{1/2}\,C_{\varepsilon},\ \sqrt{2}\,(2+\kappa_{2}^{2})^{1/2}\Big\}$$
from the inequalities \eqref{Estimate-from-ellipticity}, \eqref{Estimate_F_phi1-F_phi2}, \eqref{B_w2-B_w1}, \eqref{A_phi2-A_phi_1}, and recalling that $\bt=\A_{1}(\bw, \bphi)$ and $\bu=\A_3(\bw, \bphi)$, yields \eqref{LC-Estimate-A} and concludes the proof.
\end{proof}

\begin{lem}\label{Lipchtz-continuity-B}
There exists $\LB>0$, depending on $\widetilde{\alpha}$ and $c_1(\Omega)$ (cf. \eqref{Ellipticity-C-H1} and \eqref{Bounded-F}, respectively) such that for all $(\bw, \bphi), (\widetilde{\bw}, \widetilde{\bphi})\in W_r$ (cf. Lemma {\rm\ref{Ball-mapping}})there holds 
\begin{equation}\label{LC-Estimate-B}
\|\B(\bw, \bphi)-\B(\widetilde{\bw}, \widetilde{\bphi})\|_{1,\Omega}\ \leq\ \LB\Big\{\|\bw-\widetilde{\bw}\|_{1, \Omega}\, +\, \|\bphi-\widetilde{\bphi}\|_{1, \Omega}\Big\}\,.
\end{equation}
\end{lem}

\begin{proof}
Given $(\bw, \bphi), (\widetilde{\bw}, \widetilde{\bphi})\in W_r$, we let $\bvarphi, \widetilde{\bvarphi}\in \mathbf{H}^{1}(\Omega)$ be the corresponding solutions of \eqref{Defining_operator_B}, that is $\bvarphi:=\B(\bw, \bphi)$ and $\widetilde{\bvarphi}:=\B(\widetilde{\bw}, \widetilde{\bphi})$. Then, since $\bvarphi|_{\GD}=\widetilde{\bvarphi} |_{\GD}= \bvarphi_{\D}$, we realize that $\bvarphi-\widetilde{\bvarphi}$ belongs to $\mathbf{H}^{1}_{\GD}(\Omega)$. In this way, applying the ellipticity of $\bC$ (cf. \eqref{Ellipticity-C-H1}), using
\eqref{Defining_operator_B} and \eqref{G}, adding and subtracting $\widetilde{\bphi}$, and employing the 
H\"{o}lder inequality, the continuous injection $\mathbf{i}:\mathbf{H}^{1}(\Omega)\to \mathbf{L}^{4}(\Omega)$,
and the definition of $W_r$ (cf. Lemma \ref{Ball-mapping}), we readily deduce that
\begin{eqnarray}
\widetilde{\alpha}\; \|\bvarphi-\widetilde{\bvarphi}\|_{1, \Omega}^{2} & \leq & \bC(\bvarphi, \bvarphi-\widetilde{\bvarphi})-\bC(\widetilde{\bvarphi}, \bvarphi-\widetilde{\bvarphi})\nonumber\\[1ex] 
 & = & G_{\bw,\bphi}(\bvarphi-\widetilde{\bvarphi})-G_{\widetilde{\bw},\widetilde{\bphi}}(\bvarphi-\widetilde{\bvarphi})\label{Estimate-from-ellipticity-C}\\[1ex]
& = & -\int_{\Omega}\Big\{\big(\nabla(\bphi-\widetilde{\bphi})\big)\bw + (\nabla\widetilde{\bphi})(\bw-\widetilde{\bw})\Big\}\cdot(\bvarphi-\widetilde{\bvarphi})\nonumber\\[1ex]
& \leq & c_1(\Omega)\Big\{\|\bw\|_{1,\Omega}|\bphi-\widetilde{\bphi}|_{1,\Omega} + |\widetilde{\bphi}|_{1,\Omega}\|\bw-\widetilde{\bw}\|_{1,\Omega}\Big\}\|\bvarphi-\widetilde{\bvarphi}\|_{1,\Omega}\nonumber\\[1ex]
& \leq & rc_1(\Omega)\,\Big\{\|\bphi-\widetilde{\bphi}\|_{1, \Omega} + \|\bw-\widetilde{\bw}\|_{1, \Omega}\Big\}\|\bvarphi-\widetilde{\bvarphi}\|_{1, \Omega}\nonumber
\end{eqnarray}
which give \eqref{LC-Estimate-B} with $\LB := \frac{rc_1(\Omega)}{\widetilde{\alpha}}$.
\end{proof}

\begin{lem}\label{Lipchtz-continuity-L}
Let $r$ and $W_r$ as in Lemma {\rm\ref{Ball-mapping}}. Then, there exists a positive constant $L_{\,\tt{Lip}}$, depending on $r$,
$|\alp|$, $|\bg|$, and the constants $\CA$, $\LA$, and $\LB$ (\eqref{Estimate-operatorA}, \eqref{LC-Estimate-A}, and \eqref{LC-Estimate-B}, respectively), such that for all $(\bw, \bphi), (\widetilde{\bw}, \widetilde{\bphi})\in W_r$ there holds
\begin{equation}\label{LC-Estimate-L}
\begin{array}{l}
\|\L(\bw, \bphi)-\L(\widetilde{\bw}, \widetilde{\bphi})\|_{\mathbf{H}}\ \leq\ L_{\,\tt{Lip}}\,\|(\bw, \bphi)-(\widetilde{\bw}, \widetilde{\bphi})\|_{\mathbf{H}}\,.
\end{array}
\end{equation}
\end{lem}

\begin{proof}
Given $r\in (0,r_0)$ and $(\bw, \bphi), (\widetilde{\bw}, \widetilde{\bphi})\in W_r$, we first observe,
according to the definition of $\L$ (cf. \eqref{Defining_operator_L}), and the Lipschitz-continuity of
$\B$ (cf. \eqref{LC-Estimate-B}), that
\begin{eqnarray*}
\|\L(\bw, \bphi)-\L(\widetilde{\bw}, \widetilde{\bphi})\|_{\mathbf{H}} & = & \|\A_{3}(\bw, \bphi)-\A_{3}(\widetilde{\bw}, \widetilde{\bphi})\|_{1, \Omega}\, +\, \|\B(\A_3(\bw, \bphi), \bphi) - \B(\A_3(\widetilde{\bw}, \widetilde{\bphi}), \widetilde{\bphi})\|_{1, \Omega}\\[2ex]
& \leq & (1 + \LB)\,\|\A_{3}(\bw, \bphi)-\A_{3}(\widetilde{\bw}, \widetilde{\bphi})\|_{1, \Omega}\, +\, \LB\,\|\bphi-\widetilde{\bphi}\|_{1, \Omega}
\end{eqnarray*}
from which, employing the Lipschitz-continuity of $\A$ (cf. \eqref{Lipchtz-continuity-A}), yields
\begin{equation}\label{previous-estimate}
\begin{array}{l}
\|\L(\bw, \bphi)-\L(\widetilde{\bw}, \widetilde{\bphi})\|_{\mathbf{H}}\ \leq\ 
(1 + \LB)\LA\Big\{\|\A_3(\bw, \bphi)\|_{1, \Omega}\,\|\bw-\widetilde{\bw}\|_{1,\Omega}\\[2ex]
\qquad\qquad +\ \|\A_1(\bw,\bphi)\|_{\varepsilon,\Omega}\,\|\bphi-\widetilde{\bphi}\|_{0,n/\varepsilon,\Omega}\, +\, |\alp|\,|\bg|\,\|\bphi-\widetilde{\bphi}\|_{1, \Omega}\Big\}\, +\, \LB\,\|\bphi-\widetilde{\bphi}\|_{1, \Omega}
\end{array}
\end{equation}
Then, applying the bound \eqref{Estimate-operatorA} to estimate the term $\|\A_3(\bw, \bphi)\|_{1, \Omega}$, employing the continuous injection of $\mathbf{H}^{1}(\Omega)$ into $\mathbf{L}^{n/\varepsilon}(\Omega)$ with boundedness constant $\widetilde{C}_{\varepsilon}$, using the estimate \eqref{Regularity-estimate} to estimate the term $\|\A_1(\bw, \bphi)\|_{\varepsilon, \Omega}$, noting that $\|\bphi\|_{1,\Omega}\leq r$, and performing some algebraic manipulations, we get from \eqref{previous-estimate} that
\begin{eqnarray*}
\|\L(\bw, \bphi)-\L(\widetilde{\bw}, \widetilde{\bphi})\|_{\mathbf{H}} & \leq & (1+\LB)\LA\,|\alp|\,|\bg|\Big\{r\CA\|\bw-\widetilde{\bw}\|_{1,\Omega}\, +\, (r\widetilde{C}_{\varepsilon}\widehat{C}(r) + 1)\|\bphi-\widetilde{\bphi}\|_{1,\Omega}\Big\}\\
& & +\, \LB\|\bphi-\widetilde{\bphi}\|_{1,\Omega}\\[1ex]
& \leq & L_{\,\tt{Lip}}\|(\bw, \bphi)-(\widetilde{\bw}, \widetilde{\bphi})\|_{\mathbf{H}}
\end{eqnarray*}
In this way, \eqref{LC-Estimate-L} follows from the foregoing inequality by defining
\begin{equation}\label{def:LLip}
L_{\,\tt{Lip}}\ :=\ r\CA(1+\LB)\LA\,|\alp|\,|\bg|\, +\, (1+\LB)\LA\,|\alp|\,|\bg|(r\widetilde{C}_{\varepsilon}\widehat{C}(r) + 1)\, +\, \LB
\end{equation}
in order to complete the proof.
\end{proof}

\begin{theorem}\label{Main-result}
Suppose that the parameters $\kappa_1, \kappa_2$ and $\kappa_3$ satisfy the conditions
required by Lemma {\rm\ref{Lemma_Flow}}. Let $r$ and $W_r$ as in Lemma {\rm\ref{Ball-mapping}},
and assume that the data satisfy \eqref{Small-data-1} and allow to have
\begin{equation}\label{Small-data-2}
L_{\,\tt{Lip}}\ <\ 1\,,
\end{equation} 
with $L_{\,\tt{Lip}}$ defined in \eqref{def:LLip}.
Then, problem \eqref{augmented-mixed-primal} has a unique solution $(\vbx, \bvarphi)\in \mathbb{H}\times \mathbf{H}^{1}(\Omega)$ such that $\bvarphi|_{\GD}=\bvarphi_{\D}$, with  $(\bu, \bvarphi)\in W_r$ (cf. Lemma {\rm\ref{Ball-mapping}}), and there holds
\begin{equation}\label{C-D-D}
\|\vbx\|_{\mathbb{H}}\ \leq\ \CA\,r\,|\alp|\,|\bg|\,,
\end{equation}
and
$$\|\bvarphi\|_{1,\Omega}\ \leq\ \CB\,\big\{r\,\|\bu\|_{1,\Omega}\, +\, \|\bvarphi_{\D}\|_{1/2,\GD}\big\}\,.$$
\end{theorem}

\begin{proof}
It follows as a combination of Lemmas \ref{Ball-mapping} and \ref{Lipchtz-continuity-L}, the assumption \eqref{Small-data-2}, the Banach fixed-point theorem, and the a priori estimates \eqref{Estimate-operatorA} and \eqref{Estimate-operatorB}. We omit further details. 
\end{proof}

\section{Galerkin scheme}\label{FEM-specific-spaces}

In this section we introduce and analyze the Galerkin scheme of the semi-augmented 
mixed-primal problem \eqref{augmented-mixed-primal}.
To this end, we now let $\T$ be a regular triangulation of
$\Omega$ by triangles $K$ (in $\mathrm{R}^2$) or  tetrahedra $K$ ( in $\mathrm{R}^3$) of diameter $h_K$,
and define the meshsize $h:=\max\{h_K:\; K\in \T\}$.
In addition, given an integer $k \ge 0$, for each $K\in\T$ we let $\mathrm{P}_k(K)$ 
be the space of polynomial functions on $K$ of degree $\le\, k$, and define the corresponding
local Raviart-Thomas space of order $k$  as
\[
\mathbf{RT}_k(K) \,:=\, \mathbf{P}_k(K) \,\oplus\, \mathrm{P}_k(K)\,\bx\,,
\]
where, according to the notations described in Section \ref{Notations},
$\mathbf P_k(K) \,=\, [\mathrm P_k(K)]^n$, and $\bx$ is the generic vector
in $\mathrm{R}^n$.  Then, we consider piecewise polynomials of degree $\leq k$ for approximating entries of the strain rate $\bt$, the global Raviart-Thomas space of order $k$ to approximate rows of the pseudostress $\bsig$, and the Lagrange space given by the continuous piecewise polynomial vectors of degree $\leq k+1$ for the velocity $\bu$, respectively, that is 
\begin{eqnarray}
\mathbb{H}^{\bt}_{h} &:=&
\Big\{ \br_h\in \mathbb{L}^{2}_{\tt{tr}}(\Omega):\quad \br_{h}|_{K} \in \mathbb{P}_k(K)\quad \forall\, K \in \T\Big\}\,,\label{FE-1}\\[2ex]
\mathbf{H}^{\bsig}_{h} &:=&\Big\{ \btau_h\in  \mathbb{H}_0(\bdiv;\Omega) : \quad
\mathbf c^{\tt t}\, \btau_h|_K \in \mathbf{RT}_k(K) \quad \forall\,\mathbf c \in \mathrm{R}^n\,, 
\quad \forall\, K \in \T \Big\}\,,\label{FE-2}\\[2ex]
\mathbf{H}^{\bu}_{h} &:=& \Big\{ \bv_h\in \mathbf C(\overline{\Omega}) : \quad 
\bv_h|_K \in \mathbf{P}_{k+1}(K) \quad \forall\,  K\in \T, \quad \bv_h={\bf 0}\quad \mbox{on}\quad \Gamma\Big\}\,\label{FE-3}.
\end{eqnarray}
For the unknown $\bvarphi$ containing the tempeture and the concentration into its coordinates, we let $\mathbf{H}_{h}^{\bvarphi}\subset \mathbf{H}^{1}(\Omega)$ denote the Lagrange space of degree $\leq k+1$ with respect to $\T$ (similar to $\mathbf{H}^{\bu}_{h}$), and set 
\begin{equation}\label{FE-4}
\mathbf{H}^{\bvarphi}_{h, \GD}\;:=\;\Big\{\bpsi_h\in \mathbf{H}^{\bvarphi}_{h}:\quad \bpsi_h|_{\GD}=\mathbf{0}\Big\}
\end{equation}  
to be the analogous space with homogeneous Dirichlet boundary conditions. We define
\begin{equation}\label{def:phiDh}
\textcolor{black}{\bvarphi_{\D,h}\ :=\ \mathcal{I}_h^{SZ}(E(\bvarphi_{\D}))\big|_{\GD}}
\end{equation}
to be the approximate Dirichlet boundary data, where $\mathcal{I}_{h}^{SZ}:\mathbf{H}^{1}(\Omega)\to \mathbf{H}^{\bvarphi}_{h}$ denotes the Scott-Zhang interpolant operator of degree $k+1$, which satisfies the following stability and approximation properties, respectively, (see \cite[Lemma 1.130]{Ern-2004}).
\begin{lem}\label{Scott-Zhang}
Let $p$ and $\ell$ satisfy $1\leq p<+\infty$ and $\ell\geq 1$ if $p=1$,  and $1/p<\ell$ otherwise. Then,  there exists a positive constant $c$,  independent of $h$,  such that the following properties hold:
\begin{itemize}
\item[$(i)$]   For all $0\leq m\leq \mbox{min}\{1,\ell\}$,
\begin{equation}\label{stability} 
\forall h, \;\;\forall v\in W^{\ell,p}(\Omega), \quad  \quad \|\mathcal{I}_{h}^{SZ}(v)\|_{\,m,p,\Omega} \leq\  c\,\|v\|_{\,\ell,p,\Omega}.
\end{equation}
\item[$(ii)$]  Provided $\ell\leq k+1$,  for all $0\leq m\leq \ell$,
\begin{equation}\label{Approximation}
\forall h, \;\;\forall K\in \T, \; \forall v\in W^{\ell,p}(\Delta_K), \quad \quad \|v-\mathcal{I}_{h}^{SZ}(v)\|_{\,m, p, K} \leq \ c\,h^{\ell-m}_{K}\,|v|_{\,\ell,p,\Delta_K} 
\end{equation}
where $\Delta_K$ denotes the set of elements in $\T$, sharing at least one vertex with $K$ .
\end{itemize}
\end{lem}
Hence, $\mathbf{\bvarphi}_{\D,h}$ belongs to the discrete trace space on $\GD$ given by
$$\mathbf{H}^{1/2}_{h}(\GD)\;:=\;\Big\{\bpsi_{\D,h}\in \mathbf{C}(\GD):\quad \bpsi_{\D, h}|_{e}\in \mathbf{P}_{k+1}(e)\quad \forall \,e\in \mathcal{E}_{\GD}\Big\},$$
where $\mathcal{E}_{\GD}$ stands for the set of edges/faces on $\GD$.

\medskip
In this way, defining $\mathbb{H}_{h}:=\mathbb{H}^{\bt}_{h}\times \mathbf{H}^{\bsig}_{h} \times \mathbf{H}^{\bu}_{h}$ and denoting $\vbx_h := (\bt_h, \bsig_h, \bu_h)$, the underlying Galerkin scheme given by the discrete counterpart 
of \eqref{augmented-mixed-primal}, reads: Find 
$(\vbx_h,\bvarphi_h)\in\mathbb{H}_{h}\times\mathbf{H}^{\bvarphi}_{h}$, with $\bvarphi_{h}\big|_{\GD}=\bvarphi_{\D,h}$,   such that
\begin{subequations}\label{augmented-mixed-primal-h}
\begin{eqnarray}
\label{augmented-mixed-primal-1-h} \bA_{\bvarphi_h}(\vbx_h,\vby_h)+\bB_{\bu_h}(\vbx_h,\vby_h) & = & F_{\bvarphi_h}(\vby_h)\qquad\quad\;\; \forall\,\vby_h\in\mathbb{H}_h, \\[2ex]
\label{augmented-mixed-primal-2-h} \bC(\bvarphi_h,\bpsi_h) & = & G_{\bu_h,\bvarphi_h}(\bpsi_h)\qquad \forall\,\bpsi_h\in\mathbf{H}^{\bvarphi}_{h,\GD}.
\end{eqnarray}
\end{subequations}
Throughout the rest of this section we adopt the discrete analogue 
of the fixed point strategy introduced in Section \ref{sec:FPA}. Indeed, denoting $\mathbf{H}_h:=\mathbf{H}^{\bu}_{h}\times \mathbf{H}^{\bvarphi}_{h}$, we define  the operator $\A_h:\mathbf{H}_h \rightarrow \mathbb{H}_h$ by
$$\A_h(\bw_h,\bphi_h)\ =\ (\A_{1,h}(\bw_h,\bphi_h), \A_{2,h}(\bw_h,\bphi_h), \A_{3,h}(\bw_h, \bphi_h))\ :=\ \vbx_h\qquad  \forall\, (\bw_h, \bphi_h)\in \mathbf{H}_h$$
where $\vbx_h=(\bt_h, \bsig_h, \bu_h)\in \mathbb{H}_h$  is the unique solution of the discrete problem \eqref{augmented-mixed-primal-1-h} with $(\bw_h, \bphi_h)$ instead pf $(\bu_h, \bvarphi_h)$, that is
\begin{equation}\label{Defining_operator_A-h}
\bA_{\bphi_h}(\vbx_h, \vby_h)+\bB_{\bw_h}(\vbx_h, \vby_h)\ =\ F_{\bphi_h}(\vby_h)\qquad \forall\,\vby_h\in\mathbb{H}\,,
\end{equation}
where the bilinear forms $\bA_{\bphi_h}$, $\bB_{\bw_h}$ and the functional $F_{\bphi_h}$  are those corresponding to \eqref{A}, \eqref{B}, and \eqref{F}, respectively, with $\bw=\bw_h$ and $\bphi=\bphi_h$.

\medskip
In addition, we introduce the operator $\B_h:\mathbf{H}_h \rightarrow \mathbf{H}^{\bvarphi}_{h}$ defined as
$$\B_h(\bw_h,\bphi_h)\ :=\ \bvarphi_h\qquad \forall\, (\bw_h,\bphi_h)\in \mathbf{H}_h\,,$$
where $\bvarphi_h$ is the unique solution of the discrete problem \eqref{augmented-mixed-primal-2-h} with $(\bw_h, \bphi_h)$ instead pf $(\bu_h, \bvarphi_h)$, that is
\begin{equation}\label{Defining_operator_B-h}
\bC(\bvarphi_h,\bpsi_h)\ =\ G_{\bw_h,\bphi_h}(\bpsi_h)\qquad\forall\,\bpsi_h\in\mathbf{H}^{\bvarphi}_{h, \GD}\,,
\end{equation}
where the bilinear form $\bC$, and the functional $G_{\bw_h,\bphi_h}$ are defined as in
\eqref{C} and \eqref{G}, respectively, with $\bw=\bw_h$ and $\bphi=\bphi_h$.

\medskip
Therefore, by introducing the operator $\L_h: \mathbf{H}_h \rightarrow \mathbf{H}_h$ as
\begin{equation*}\label{Defining_operator_L-h}
\L_h(\bw_h,\bphi_h)\ :=\ (\A_{3,h}(\bw_h, \bphi_h),\, \B_h(\A_{3,h}(\bw_h,\bphi_h), \bphi_h))\qquad \forall\, (\bw_h, \phi_h)\in \mathbf{H}_h
\end{equation*}
we realize that solving \eqref{augmented-mixed-primal-h} is equivalent to seeking a fixed point of $\L_h$, that is: Find
$(\bu_h,\bvarphi_h)\in \mathbf{H}_h$ such that
\begin{equation}\label{Fixed-point-equation-h}
\L_h(\bu_h,\bvarphi_h)\ =\ (\bu_h,\bvarphi_h)\,.
\end{equation}
Certainly, all the above makes sense if we guarantee that the discrete problems
\eqref{Defining_operator_A-h} and \eqref{Defining_operator_B-h} are well-posed.
This is precisely the purpose of the next section.

\subsection{Well-posedness of the uncoupled problems}\label{sec:WD-FPA-h}

In this section, we establish the well-posedness of both \eqref{Defining_operator_A-h} and \eqref{Defining_operator_B-h}, thus confirming that the operators $\A_h$, $\B_h$, and hence $\L_h$, are well-defined. We begin with the corresponding result for $\A_h$, which actually follows almost verbatim to that of its continuous counterpart $\A$ (see Lemma \ref{Lemma_Flow}), and the proof can be omitted.
\begin{lem}\label{Lemma_Flow-h}
Assume that $\kappa_1\in\big(0, 2\delta_3\big(2\nu_1-\frac{\kappa_3\nu_2}{\delta_1}\big)\big)$, $\kappa_2\in(0,2\delta_2)$, $\kappa_3\in\big(0,\frac{2\nu_1\delta_1}{\nu_2}\big)$ with $\delta_1\in\big(0,\frac{1}{\nu_2}\big)$, $\delta_2\in\big(0, \frac{2}{\gamma}\big)$, $\delta_3\in (0,1)$. Then, there exists $r_0>0$ (cf. \eqref{radius}) such that for each $r\in (0,r_0)$, problem \eqref{Defining_operator_A-h} has a unique solution $\A_h(\bw_h, \bphi_h):=\vbx_h\in \mathbb{H}_h$, for each $(\bw_h, \bphi_h)\in \mathbf{H}_h$ with $\|\bw_h\|_{1, \Omega}\leq r$. Moreover, there exists  $\CA$, independent of $(\bw_h, \bphi_h)$, such that
\begin{equation*}\label{Estimate-operatorA-h}
\|\A_h(\bw_h, \bphi_h)\|_{\mathbb{H}_h}\ =\ \|\vbx_h\|_{\mathbb{H}_h}\ \leq\ \CA\,|\alp|\,|\bg|\,\|\bphi_h\|_{0,\Omega}\qquad \forall\,(\bw_h, \bphi_h)\in \mathbf{H}_h\,.
\end{equation*}
\end{lem}
We now provide  the discrete version of Lemma {\rm\ref{Lemma-Transport}}.

\begin{lem}\label{Lemma-Transport-h}
For each $(\bw_h, \bphi_h)\in \mathbf{H}_h$, problem \eqref{Defining_operator_B-h} has a unique solution $\bvarphi_h\in \mathbf{H}_{h}^{\bvarphi}$, with $\bvarphi_h|_{\Gamma_{\D}}=\bvarphi_{\D,h}$. Moreover, there exists a constant $\CBt>0$ independent of $(\bw_h, \bphi_h)$, such that 
\begin{equation}\label{Estimate-operatorB-h}
\|\B_h(\bw_h, \bphi_h)\|_{1, \Omega}\ =\ \|\bvarphi_h\|_{1, \Omega}\ \leq\ \CBt\big\{\|\bw_h\|_{1,\Omega}\|\bphi_h\|_{1,\Omega}\, +\, |\alp|\,|\bg|\,+\,\|\bvarphi_{\D}\|_{1/2,\GD}\big\}\,.
\end{equation}
\end{lem}

\begin{proof} 
Let $\textcolor{black}{\bvarphi_{1,h}:=\mathcal{I}_{h}^{SZ}(\bvarphi_1)}\in \mathbf{H}_{h}^{\bvarphi}$
which satisfies $\bvarphi_{1,h}|_{\GD} = \bvarphi_{\D,h}$. Then, similar to Lemma \ref{Lemma-Transport},
we consider the auxiliary discrete problem: Find $\bvarphi_{0,h}\in\mathbf{H}^{\bvarphi}_{h,\GD}$ such that
\begin{equation}\label{alternative_problem-h}
\bC(\bvarphi_{0,h}, \bpsi_h)\ =\ \widetilde{G}_{\bw_h,\bphi_h}(\bpsi_h)\qquad\forall\,\bpsi_h\in \mathbf{H}^{\bvarphi}_{h,\GD}\,,
\end{equation}
where 
\begin{equation*}\label{G_tilde-h}
\widetilde{G}_{\bw_h,\bphi_h}(\bpsi_h)\ :=\ -\int_{\Omega}( \nabla \bphi_h)\bw_h\cdot\bpsi_h\, -\, \int_{\Omega}\mathbb{K}\nabla \bvarphi_{1,h}:\nabla \bpsi_h \qquad \forall\, \bpsi_h \in \mathbf{H}^{\bvarphi}_{h,\GD}\,.
\end{equation*}
Next, the boundedness and the ellipticity of $\bC$ are obtained exactly as in the proof
of Lemma \ref{Lemma-Transport} with the same ellipticity constant $\widetilde{\alpha}$
given by \eqref{eq:elip-const-trans}.
On the other hand, reasoning as in the proof of Lemma \ref{Lemma-Transport}, \textcolor{black}{and applying the stability property of $\mathcal{I}_{h}^{SZ}$ (cf. \eqref{stability}) and the equinormic property of $\bvarphi_1$ (cf.  \eqref{eqn:wp-C})},  we easily deduce that
\begin{equation*}\label{Bounded-F-h}
|\widetilde{G}_{\bw_h, \bphi_h} (\bpsi_h)|\ \leq\ \Big\{c_1(\Omega)\|\bw_h\|_{1,\Omega}\|\bphi_h\|_{1,\Omega}\, +\textcolor{black}{\, c\,\|\mathbb{K}\|_{\infty, \Omega}\,\|\bvarphi_{\D}\|_{1/2,\GD}}\Big\}\|\bpsi_h\|_{1, \Omega},
\end{equation*}
$\forall\, \bpsi \in \mathbf{H}^{\bvarphi}_{h,\GD}$,  which says that $\widetilde{G}_{\bw_h, \bphi_h}\in [\mathbf{H}^{\bvarphi}_{h,\GD}]^{\prime}$
and 
$$\|\widetilde{G}_{\bw_h, \bphi_h}\| \leq c_1(\Omega)\|\bw_h\|_{1,\Omega}\|\bphi_h\|_{1,\Omega}\, +\textcolor{black}{\, c\,\|\mathbb{K}\|_{\infty, \Omega}\,\|\bvarphi_{\D}\|_{1/2,\GD}}.$$
Therefore, a direct application of the Lax-Milgram lemma
implies that there exists a unique $\bvarphi_{0,h}\in \mathbf{H}^{\bvarphi}_{h,\GD}$ that satisfies
\eqref{alternative_problem-h}  with 
$$\|\bvarphi_{0,h}\|_{1,\Omega}\ \leq\ \frac{1}{\widetilde{\alpha}}\,\Big\{c_1(\Omega)\|\bw_h\|_{1,\Omega}\|\bphi_h\|_{1,\Omega}\, +\textcolor{black}{\, c\,\|\mathbb{K}\|_{\infty, \Omega}\,\|\bvarphi_{\D}\|_{1/2,\GD}}\Big\}.$$
Then,  $\bvarphi_h\;:=\;\bvarphi_{0,h}+\bvarphi_{1,h}\in\mathbf{H}^{\bvarphi}_{h}$,
which in fact satisfies that $\bvarphi_{h}|_{\GD}=\bvarphi_{\D,h}$,  is the
unique solution of \eqref{Defining_operator_B-h}. In addition, the estimate \eqref{Estimate-operatorB-h}
holds with $\textcolor{black}{\CBt:=\max\{\widetilde{\alpha}^{-1}\,c_1(\Omega), \widetilde{c}\,(\widetilde{\alpha}^{-1}\,\|\mathbb{K}\|_{\infty, \Omega}+1)\}}$.
\end{proof}

\subsection{Solvability analysis of the fixed point equation}\label{sec:WD-FPA2-h}

In this section we establish the solvability of the fixed point problem \eqref{Fixed-point-equation-h}
by applying the Brouwer fixed-point theorem \cite[Theorem 9.9-2]{ciarlet-2013}. To this end, we begin
with the discrete version of lemma \ref{Ball-mapping}.
\begin{lem}\label{Ball-mapping-h}
Given $r\in (0,r_0)$, with $r_0$ given by \eqref{radius}, we let
$$W_r^h\ :=\ \big\{(\bw_h, \bphi_h)\in \mathbf{H}_h:\quad \|(\bw_h, \bphi_h)\|_{\mathbf{H}}\,\leq\, r\big\}\,,$$
and assume that the data satisfy
\begin{equation}\label{Small-data-1-h}
c(r)\,|\alp|\,|\bg|\, +\, \CBt\,\|\bvarphi_{\D}\|_{1/2, \GD}\ \leq\ r
\end{equation}
where $\textcolor{black}{c(r):=r\CA(1 + r\CBt)}$, with $\CA$ and $\CBt$ as in \eqref{Estimate-operatorA} and
\eqref{Estimate-operatorB-h}, respectively. Then, there holds $\L_h(W_r^h)\subseteq W_r^h$.
\end{lem}
\noindent In order to provide the discrete analogue of Lemma \ref{Lipchtz-continuity-A}, we notice in advance that, instead of the regularity assumptions employed in the continuous case, which are not applicable in the present case, we simple utilize an $\mathrm{L}^4$-$\mathrm{L}^{4}$-$\mathrm{L}^{2}$ argument.

\begin{lem}\label{Lipchtz-continuity-A-h}
Let $r\in (0, r_0)$ with $r_0$ given by \eqref{radius}. Then, there exists a constant $\LAt>0$, independent of $r$, 
such that for all $(\bw_h, \bphi_h), (\widetilde{\bw}_h, \widetilde{\bphi}_h)\in W_r^h$, with $\|\bw_h\|_{1,\Omega}, \|\widetilde{\bw}_h\|_{1,\Omega}\,\leq \;r$, there holds
\begin{equation}\label{LC-Estimate-A-h}
\begin{array}{l}
\|\A_h(\bw_h, \bphi_h)-\A_h(\widetilde{\bw}_h, \widetilde{\bphi}_h)\|_{\mathbb{H}}\ \leq\ \LAt\Big\{\|\A_{3,h}(\bw_h, \bphi_h)\|_{1, \Omega}\;\|\bw_h-\widetilde{\bw}_h\|_{1,\Omega} \\[1ex]
\qquad \qquad \qquad \qquad +\ \|\A_{1,h}(\bw_h, \bphi_h)\|_{0,4, \Omega}\,\|\bphi-\widetilde{\bphi}\|_{0,4,\Omega}\, +\, |\alp|\,|\bg|\,\|\bphi_h-\widetilde{\bphi}_h\|_{1, \Omega}\Big\}\,.
\end{array}
\end{equation}
\end{lem}

\begin{proof} It proceeds exactly as in the proof of Lemma \ref{Lipchtz-continuity-A}, except for the derivation of the discrete analogue of \eqref{A_phi2-A_phi_1}, where, instead of choosing the values of $p, q$ determined by the regularity parameter $\varepsilon$, it suffices to take $p=q=2$, thus obtaining
\begin{equation*}\label{A_phi2-A_phi_1-h}
\begin{array}{l}
\displaystyle |(\bA_{\widetilde{\bphi}_{h}}-\bA_{\bphi_{h}})(\vbx_{h}, \vbx_{h}- \vbz_{h})|\ \leq\
2\,L_{\nu}(2+\kappa_3^2)^{1/2}\,\|\bt_{h}\|_{0,4,\Omega}\,\|\bphi_{h}-\widetilde{\bphi}_{h}\,\|_{0,4,\Omega}\,\|\vbx_{h}-\vbz_{h}\|_{\mathbb{H}}\,,
\end{array}
\end{equation*}
for all $(\bw_{h}, \bphi_{h}), (\widetilde{\bw}_{h}, \widetilde{\bphi}_{h})$, with  $\vbx_{h}=(\bt_{h}, \bsig_{h}, \bu_{h}):=\A_h(\bw_{h}, \bphi_{h})\in \mathbb{H}_h$ and  $\vbz_{h}=(\widetilde{\bt}_{h}, \widetilde{\bsig}_{h}, \widetilde{\bu}_{h}):=\A_h(\widetilde{\bw}_{h}, \widetilde{\bphi}_{h})\in \mathbb{H}_h$. Thus, since the elements of $\mathbb{H}_{h}^{\,\bt}$ are piecewise polynomials, we can guarantee
that $\|\bt_h\|_{0,4,\Omega} < \infty$ for each $\bt_h\in \mathbb{H}_{h}^{\,\bt}$.
The proof concludes with $\LAt := \frac{2}{\alpha(\Omega)}\max\{c_1(\Omega)\,(2+\kappa_3^2)^{1/2},$ $2\,L_{\nu}(2+\kappa_3^2)^{1/2}, \sqrt{2}\,(2+\kappa_{2}^{2})^{1/2}\}$.
\end{proof}

The discrete version of Lemma \ref{Lipchtz-continuity-B} is given as follows.
\begin{lem}\label{Lipchtz-continuity-B-h}
Let $\LB>0$ as in Lemma {\rm\ref{Lipchtz-continuity-B}}. Then,
for all $(\bw_h, \bphi_h), (\widetilde{\bw}_{h}, \widetilde{\bphi}_{h})\in W_r^h$ there holds 
\begin{equation}\label{LC-Estimate-B-h}
\begin{array}{l}
\|\B_h(\bw_{h}, \bphi_{h})-\B_h(\widetilde{\bw}_{h}, \widetilde{\bphi}_{h})\|_{1,\Omega}\ \leq\ \LB\Big\{\|\bw_{h}-\widetilde{\bw}_{h}\|_{1, \Omega}\, +\, \|\bphi_{h}-\widetilde{\bphi}_{h}\|_{1, \Omega}\Big\}\,.
\end{array}
\end{equation}
\end{lem}
\begin{proof}
It corresponds to an adaptation of the proof of Lemma \ref{Lipchtz-continuity-B} to the discrete context.
\end{proof}

Now, combining Lemmas \ref{Lipchtz-continuity-A-h} and \ref{Lipchtz-continuity-B-h}, and employing the continuous injection of $\mathbf{H}^{1}(\Omega)$ into $\mathbf{L}^{4}(\Omega)$, we can prove the discrete version of Lemma \ref{Lipchtz-continuity-L}.
\begin{lem}\label{Lipchtz-continuity-L-h}
Let $r$ and $W_r^h$ as in Lemma \ref{Ball-mapping-h}. Then, there exits a positive constant $C$, depending only on $\LAt$ and $\LB$ (cf. \eqref{LC-Estimate-A-h} and \eqref{LC-Estimate-B-h}, respectively), such that for all $(\bw_h, \bphi_h), (\widetilde{\bw}_{h}, \widetilde{\bphi}_h)\in W_r^h$ there holds
\begin{eqnarray*}
\|\L_h(\bw_h, \bphi_h)-\L_h(\widetilde{\bw}_h, \widetilde{\bphi}_h)\|_{\mathbf{H}} & \leq & C\,\Big\{\|\A_{3,h}(\bw_h,\bphi_{h})\|_{1,\Omega}\, +\, c_1(\Omega)\,\|\A_{1,h}(\bw_h,\bphi_{h})\|_{0,4,\Omega}\nonumber\\[1ex]
&& \phantom{C\,\Big\{} +\, |\alp|\,|\bg|\, +\, \LB\Big\}\,\|(\bw_h, \bphi_h)-(\widetilde{\bw}_h, \widetilde{\bphi}_h)\|_{\mathbf{H}}\,.\label{LC-Estimate-L-h}
\end{eqnarray*}
More precisely, there holds $C = \max\{(1+\LB)\LAt, 1\}$.
\end{lem}

Consequently, since the foregoing lemma confirms the continuity of $\L_h$, by a straightforward
application of Brouwer fixed point theorem (cf. \cite[Theorem 9.9-2]{ciarlet-2013}) on the convex
and compact set $W_r^h\subseteq \mathbf{H}_h$, we can provide the main result of this section.
\begin{theorem}\label{Main-result-h}
Suppose that the parameters $\kappa_1, \kappa_2$ and $\kappa_3$ satisfy the conditions required
by Lemma {\rm\ref{Lemma_Flow-h}}. Let $r$ and $W_r^h$ as in Lemma {\rm\ref{Ball-mapping-h}}, and
assume that the data satisfy \eqref{Small-data-1-h}.
Then, problem \eqref{augmented-mixed-primal-h} has at least one solution $(\vbx_h, \bvarphi_h)\in \mathbb{H}_h\times \mathbf{H}_{h}^{\bvarphi}$ such that $\bvarphi_h|_{\GD}=\bvarphi_{\D,h}$, with  $(\bu_h, \bvarphi_h)\in W_r^h$, and there holds
\begin{equation*}\label{C-D-D-h}
\|\vbx_h\|_{\mathbb{H}}\ \leq\ \CA\,r\,|\alp|\,|\bg|\,,
\end{equation*}
and
$$\|\bvarphi_h\|_{1,\Omega}\ \leq\ \CBt\,\big\{r\,\|\bu_h\|_{1, \Omega}\,+\, \textcolor{black}{\|\bvarphi_{\D}\|_{1/2, \GD}}\big\}\,.$$
\end{theorem}
\section{A priori error analysis}\label{section5}

We now aim to derive the a priori error estimates for the Galerkin scheme given by
\eqref{augmented-mixed-primal-h}. To this end, given $((\bt, \bsig, \bu),\bvarphi):=(\vbx, \bvarphi)\in \mathbb{H}\times \mathbf{H}^{1}(\Omega)$,
with  $(\bu, \bvarphi)\in W_r$ and $((\bt_h, \bsig_h, \bu_h), \bvarphi_h) := (\vbx_h,\bvarphi_h)\in \mathbb{H}_h\times \mathbf{H}^{\bvarphi}_h$, with  $(\bu_h,\bvarphi_h)\in W_r^h$ solutions of \eqref{augmented-mixed-primal}
and \eqref{augmented-mixed-primal-h}, respectively, we first observe that the above problems can be rewritten as two pairs of corresponding continuous and discrete formulations, namely
\begin{subequations}\label{augmented-mixed-Flow}
\begin{eqnarray}
\label{augmented-mixed}   \bA_{\bvarphi}(\vbx,\vby)+\bB_{\bu}(\vbx,\vby) & = & F_{\bvarphi}(\vby)\qquad\quad \forall\,\vby\in\mathbb{H}\,,\\[1ex]
\label{augmented-mixed-h} \bA_{\bvarphi_h}(\vbx_h,\vby_h)+\bB_{\bu_h}(\vbx_h,\vby_h) & = & F_{\bvarphi_h}(\vby_h)\qquad \forall\,\vby_h\in\mathbb{H}_h\,,
\end{eqnarray}
\end{subequations}
and 
\begin{subequations}\label{primal-Transport}
\begin{eqnarray}
\label{primal}   \bC(\bvarphi,\bpsi) & = & G_{\bu,\bvarphi}(\bpsi)\qquad\quad\, \forall\,\bpsi\in\mathrm{H}^1_{\GD}(\Omega)\,,\\[1ex]
\label{primal-h} \bC(\bvarphi_h,\bpsi_h) & = & G_{\bu_h,\bvarphi_h}(\bpsi_h)\qquad \forall\,\bpsi_h\in\mathbf{H}^{\bvarphi}_{h, \GD}\,.
\end{eqnarray}
\end{subequations}
Our goal is to obtain an upper bound for the error $\displaystyle \|(\vbx, \bvarphi)-(\vbx_h, \bvarphi_h)\|_{\mathbb{H}\times\mathbf{H}^{1}(\Omega)}$. For this purpose, we first recall from \cite[Theorem 11.1]{Thomas-1991} an abstract result that corresponds to the standar Strang Lemma for elliptic variational problems, which will be straightforwardly applied to the pair \eqref{augmented-mixed}-\eqref{augmented-mixed-h}.

\begin{lem}\label{Strang}
Let $H$ be a Hilbert space, $F\in H^{\prime}$, and $A : H\times H\to\text{\rm R}$ be a
bounded and $H$-elliptic bilinear form. In addition, let $\{H_h\}_{h>0}$ be a sequence
of finite dimensional subspaces of $H$, and for each $h>0$ consider a bounded bilinear
form $A_h : H_h \times  H_h\to\text{\rm R}$ and a functional $F_h\in H_h^{\prime}$. 
Assume that the family $\{A_h\}_{h>0}$ is uniformly elliptic, that is, there exists a
constant $\beta>0$, independent of $h$ such that
$$A_{h}(v_h,v_h)\ \geq\ \beta\,\|v_h\|^2_H\qquad \forall\, v_h\in H_h\,, \quad \forall\, h>0\,.$$
In turn, let $u\in H$ and $u_h\in H_h$ such that
$$A(u,v)\ =\ F(v) \qquad\forall\, v\in H \qquad\text{and}\qquad A_h(u_h,v_h)\ =\ F_h(v_h) \qquad \forall\, v_h \in H_h\,.$$
Then for each $h>0$ there holds
\begin{eqnarray*}
\|u-u_h\|_H & \leq & C\,\Bigg\{\sup_{\substack{ {w_h\in H_h}\\ {w_h\neq \bf0} }}
\frac{\big|\,F(w_h) - F_{h}(w_h)\,\big|}{\|w_h\|_H}\, +\, \inf_{\substack{ {v_h \in H_h}\\ {v_h \neq\bf0} }}\Bigg(\|u-v_h\|_{H}\\[1ex]
&& \phantom{C\,\Bigg\{} +\, \sup_{\substack{ {w_h \in H_h}\\ {w_h\neq\bf0} }}\frac{\big|\,A(v_h,w_h)- A_h(v_h,w_h)\, \big|}{\|w_h\|_H}\Bigg)\Bigg\}\,,
\end{eqnarray*}
where $C := \beta^{-1}\,\max\{1, \|A\|\}$.
\end{lem}

In what follows, as usual, we denote
\[\mathrm{dist}\big(\vbx,\mathbb{H}_h\big)\ :=\ \inf_{\vby_h\in \mathbb{H}_h}\|\vbx - \vby_h\|_\mathbb{H}\qquad\text{and}\qquad
  \mathrm{dist}\big(\bvarphi,\mathbf{H}^{\varphi}_{h}\big)\ :=\ \inf_{\bpsi_h \in \mathbf{H}^{\bvarphi}_{h}}\|\bvarphi - \bpsi_h\|_{1,\Omega}\,.
\]
We now derive the a preliminary estimate for the error $\|\vbx-\vbx_h\|_{\mathbb{H}}=\|(\bt, \bsig, \bu)-(\bt_h, \bsig_h, \bu_h)\|_{\mathbb{H}}$.
\begin{lem}\label{Estimate-Error-Flow}
There exists a constant $C_{\tt{ST}}>0$ independent of $h$, such that
\begin{eqnarray}
\|\vbx-\vbx_h\|_{\mathbb{H}} & \leq & C_{\tt{ST}}\Big\{\mathrm{dist}\big(\vbx, \mathbb{H}_h\big)\, +\, |\alp|\,|\bg|\,\|\bvarphi-\bvarphi_h\|_{1, \Omega}\, +\, \|\bt\|_{\varepsilon, \Omega}\,\|\bvarphi-\bvarphi_h\|_{1, \Omega}\nonumber\\[1ex]
&& \phantom{C_{\tt{ST}}\Big\{} +\, \|\bu\|_{1, \Omega}\,\|\bu-\bu_h\|_{1, \Omega}\Big\}\,.\label{Estimate-error-F}
\end{eqnarray}
\end{lem}

\begin{proof}
From Lemma \ref{Lemma_Flow} we have that the bilinear forms $\bA_{\bvarphi}+\bB_{\bu}$ and
$\bA_{\bvarphi_h}+\bB_{\bu_h}$ are both bounded and uniformly elliptic, with ellipticity constant
$\frac{\alpha(\Omega)}{2}$ (cf. \eqref{Ellipticity-A_phi+B_w}). In turn, $F_{\bvarphi}$ and $F_{\bvarphi_h}$ are linear bounded functionals in $\mathbb{H}$ and $\mathbb{H}_h$ respectively. Then, a straightforward application of Lemma \ref{Strang} to the context given by \eqref{augmented-mixed}-\eqref{augmented-mixed-h}, gives
\begin{equation}\label{Eq-strang-flow}
\begin{array}{l}
\disp \|\vbx-\vbx_h\|_{\mathbb{H}}\ \leq\ C\,\Bigg\{\sup_{\substack{ {\vby_h\in \mathbb{H}_h}\\ {\vby_h\neq \bf0} }}
\frac{\big|\,F_{\bvarphi}(\vby_h) - F_{\bvarphi_h}(\vby_h)\,\big|}{\|\vby_h\|_{\mathbb{H}}}\\[5ex]
\disp\qquad\qquad +\ \inf_{\substack{ {\vbz_h \in \mathbb{H}_h}\\ {\vbz_h \neq\bf0}}}\Bigg(\|\vbx-\vbz_h\|_{\mathbb{H}} \,+\, \sup_{\substack{ {\vby_h \in \mathbb{H}_h}\\ {\vby_h\neq\bf0}}}
\frac{\big|\,(\bA_{\bvarphi} + \bB_{\bu})(\vbz_h,\vby_h)- (\bA_{\bvarphi_h}+\bB_{\bu_h})(\vbz_h,\vby_h)\, \big|}{\|\vby_h\|_{\mathbb{H}}}\Bigg)\Bigg\}\,,
\end{array}
\end{equation}
where $C:=\frac{2}{\alpha(\Omega)}\max\{1, \|\bA_{\bvarphi}+\bB_{\bu}\|\}$.
It is important to recall here, from \eqref{Bounded-A_phi+B_w}, that $\|\bA_{\bvarphi}+\bB_{\bu}\|$
depends only on $\kappa_1$, $\kappa_2$, $\kappa_3$, $\nu_2$, $\gamma$, $c_1(\Omega)$ and $\|\bu\|_{1,\Omega}$,
where $\|\bu\|_{1,\Omega}\leq r$. Furthermore, we now proceed to estimate each term appearing at the
right-hand side of \eqref{Eq-strang-flow}. In order to do that, we first apply the same arguments
employed to obtain \eqref{Estimate_F_phi1-F_phi2}, to find that 
\begin{equation}\label{Estimate-I}
|F_{\bvarphi}(\vby_h) - F_{\bvarphi_h}(\vby_h)|\ \leq\ \sqrt{2}\,(2+\kappa_2^2)^{1/2}\,|\alp|\,|\bg|\,\|\bvarphi-\bvarphi_h\|_{1,\Omega}\,\|\vby_h\|_{\mathbb{H}}\,.
\end{equation}
Next, in order to estimate the last supremum in \eqref{Eq-strang-flow},
we add and subtract $\vbx := (\bt,\bsig,\bu)$, we note that
\begin{eqnarray*}
&   & (\bA_{\bvarphi}+\bB_{\bu})(\vbz_h, \vby_h)-(\bA_{\bvarphi_h}+\bB_{\bu_h})(\vbz_h, \vby_h)\\[1ex]
& = & (\bA_{\bvarphi}+\bB_{\bu})(\vbx, \vby_h)-(\bA_{\bvarphi_h}+\bB_{\bu_h})(\vbx, \vby_h) - 
      (\bA_{\bvarphi}+\bB_{\bu})(\vbx-\vbz_h, \vby_h) + (\bA_{\bvarphi_h}+\bB_{\bu_h})(\vbx-\vbz_h, \vby_h)\\[1ex]
 & = & (\bA_{\bvarphi}-\bA_{\bvarphi_h})(\vbx, \vby_h)\, +\, \bB_{\bu-\bu_h}(\vbx,\vby_h) - 
      (\bA_{\bvarphi}+\bB_{\bu})(\vbx-\vbz_h, \vby_h) + (\bA_{\bvarphi_h}+\bB_{\bu_h})(\vbx-\vbz_h, \vby_h)
\end{eqnarray*}
where, applying the same approach used in \eqref{A_phi2-A_phi_1} and \eqref{B_w2-B_w1},
together with \eqref{epsilon-estimate} and the continuous embedding $\mathbf{H}^{1}(\Omega)\to \mathbf{L}^{n/\varepsilon}(\Omega)$
with constant $\widetilde{C}_{\varepsilon}$, and the boundedness of the bilinear forms $\bA_{\bvarphi}+\bB_{\bu}$ and $\bA_{\bvarphi_h}+\bB_{\bu_h}$, it follows that
\begin{equation}\label{Estimate-II}
\begin{array}{r}
|(\bA_{\bvarphi}+\bB_{\bu})(\vbz_h, \vby_h)-(\bA_{\bvarphi_h}+\bB_{\bu_h})(\vbz_h, \vby_h)|\ \leq\ \Big\{ 2\,L_{\nu}C_{\varepsilon}\widetilde{C}_{\varepsilon}(2+\kappa_3^2)^{1/2}\,\|\bt\|_{\varepsilon, \Omega}\,\|\bvarphi-\bvarphi_h\|_{1, \Omega}\\[2ex]
\qquad +\ c_1(\Omega)\,(2+\kappa_3^2)^{1/2}\,\|\bu\|_{1, \Omega}\,\|\bu-\bu_h\|_{1, \Omega}\\[2ex]
\qquad +\ \big(\|\bA_{\bvarphi}+\bB_{\bu}\| + \|\bA_{\bvarphi_h} + \bB_{\bu_h}\|\big)\|\vbx-\vbz_h\|_{\mathbb{H}}\Big\}\|\vby_h\|_{\mathbb{H}}\,.
\end{array}
\end{equation}
In this way, by replacing \eqref{Estimate-I} and \eqref{Estimate-II} back into \eqref{Eq-strang-flow},
we obtain \eqref{Estimate-error-F} with $C_{\tt{ST}}$ is a positive constant depending on $\alpha(\Omega)$, $L_{\nu}$, $C_{\varepsilon}$, $\widetilde{C}_{\varepsilon}$, $c_1(\Omega)$, $\kappa_1$, $\kappa_2$, $\kappa_3$, $\nu_2$, $\gamma$, and $r$.
\end{proof}

\medskip
The following result present a  estimate for the error $\|\bvarphi-\bvarphi_h\|_{1,\Omega}$.
\begin{lem}\label{Estimate-Error-Transport}
Assume that $r$ satisfy that
\begin{equation}\label{Data-constraint-1}
r\,\frac{c_1(\Omega)}{\widetilde{\alpha}}\ \leq\ \frac{1}{2}\,,
\end{equation}
where $\widetilde{\alpha}$ is defined in \eqref{eq:elip-const-trans},
and $c_1(\Omega):=\|\boldsymbol{i}\|^2$ is the boundedness constants
of the continuous injection $\boldsymbol{i}:\mathbf{H}^{1}(\Omega)\to \mathbf{L}^{4}(\Omega)$.
Then, there exists a constant $\widetilde{C}_{\tt{ST}}>0$, independent of $h$, such that 
\begin{equation}\label{Cea-Estimate-Transport}
\|\bvarphi-\bvarphi_h\|_{1, \Omega}\ \leq\ \widetilde{C}_{\tt{ST}}\,\|\bvarphi-\mathcal{I}^{SZ}_{h}(\bvarphi)\|_{1, \Omega}\, +\, \|\bu-\bu_h\|_{1, \Omega}\,.
\end{equation}
where $\mathcal{I}_{h}^{SZ}$ denotes the Scott-Zhang interpolant operator introduced in Section \ref{FEM-specific-spaces}.
\end{lem}

\begin{proof}
We proceed similarly as in the proof of \cite[Lemma 5.3]{alvarez-2016}. Indeed, by applying the triangle inequality we have that
\begin{equation}\label{Eq-1}
\|\bvarphi-\bvarphi_h\|_{1, \Omega}\ \leq\ \|\bvarphi-\textcolor{black}{\bpsi_{h}}\|_{1, \Omega}\, +\, \|\textcolor{black}{\bpsi_{h}}-\bvarphi_h\|_{1,\Omega},
\end{equation}
where $\textcolor{black}{\bpsi_{h}\;=\;\mathcal{I}^{SZ}_{h}(\bvarphi)} \in \mathbf{H}^{\bvarphi}_h$.
Moreover, nothing that $\bvarphi|_{\GD} = E(\bvarphi_{\D})|_{\GD} = \bvarphi_{\D}$, we can employ 
\cite[eq. (2.17)]{ScottZhang90} and \eqref{def:phiDh}, to obtain that $\bpsi_{h}|_{\GD}=\bvarphi_{\D,h}$.
Now, utilizing the ellipticity of bilinear the form $\bC(\cdot, \cdot)$ on $\mathbf{H}_{h, \GD}^{\bvarphi}$
(see the proof of Lemma \ref{Lemma-Transport-h}) with constant $\widetilde{\alpha}$, along with the fact that
$\bC(\bvarphi_h, \textcolor{black}{\bpsi_{h}}-\bvarphi_h)=G_{\bu_h, \bvarphi_h}(\textcolor{black}{\bpsi_{h}}-\bvarphi_h)$ and $\bC(\bvarphi, \textcolor{black}{\bpsi_{h}}-\bvarphi_h) = G_{\bu, \bvarphi}(\textcolor{black}{\bpsi_{h}}-\bvarphi_h)$ (see \eqref{primal-Transport}), we deduce that
\begin{eqnarray}
\widetilde{\alpha}\,\|\textcolor{black}{\bpsi_{h}}-\bvarphi_h\|_{1,\Omega}^2 & \leq & \bC(\bvarphi,\textcolor{black}{\bpsi_{h}}-\bvarphi_h) - \bC(\bvarphi-\textcolor{black}{\bpsi_{h}},\textcolor{black}{\bpsi_{h}}-\bvarphi_h) - \bC(\bvarphi_h,\textcolor{black}{\bpsi_{h}}-\bvarphi_h)\nonumber\\[1ex]
& \leq & |G_{\bu, \bvarphi}(\textcolor{black}{\bpsi_{h}}-\bvarphi_h)-G_{\bu_h, \bvarphi_h}(\textcolor{black}{\bpsi_{h}}-\bvarphi_h)|\, +\, |\bC(\bvarphi-\textcolor{black}{\bpsi_{h}}, \textcolor{black}{\bpsi_{h}}-\bvarphi_h)|\,.\label{Eq-2}
\end{eqnarray}
Next, we apply the estimate \eqref{Estimate-from-ellipticity-C} to bound the first term on the
right-hand side of \eqref{Eq-2}, whereas for second term we use the boundedness of $\bC(\cdot, \cdot)$
(cf. \eqref{Bounded-C}), and \eqref{Data-constraint-1}. Then, it follows that
\begin{eqnarray}
\|\textcolor{black}{\bpsi_{h}}-\bvarphi_h\|_{1, \Omega} & \leq & \frac{\|\mathbb{K}\|_{\infty,\Omega}}{\widetilde{\alpha}}\,\|\bvarphi-\textcolor{black}{\bpsi_{h}}\|_{1, \Omega}\ +\ \frac{rc_1(\Omega)}{\widetilde{\alpha}}\,\Big\{\|\bvarphi-\bvarphi_h\|_{1,\Omega}\, +\, \|\bu-\bu_h\|_{1,\Omega}\Big\}\nonumber\\[1ex]
& \leq & \frac{\|\mathbb{K}\|_{\infty,\Omega}}{\widetilde{\alpha}}\,\|\bvarphi-\textcolor{black}{\bpsi_{h}}\|_{1, \Omega}\ +\ \frac{1}{2}\,\Big\{\|\bvarphi-\bvarphi_h\|_{1,\Omega}\, +\, \|\bu-\bu_h\|_{1,\Omega}\Big\}\,.\label{Eq-3}
\end{eqnarray}
Then, by replacing \eqref{Eq-3} back into \eqref{Eq-1},  we conclude \eqref{Cea-Estimate-Transport}
with $\widetilde{C}_{\tt{ST}} := 2(1 + \widetilde{\alpha}^{-1}\|\mathbb{K}\|_{\infty, \Omega})$.
\end{proof}

\medskip
We now proceed to combine Lemmas \ref{Estimate-Error-Flow}  and \ref{Estimate-Error-Transport}
to derive the C\'ea estimate for the total error $$\|\vbx-\vbx_h\|_{\mathbb{H}}\, +\, \|\bvarphi-\bvarphi_h\|_{1, \Omega}\,.$$
In fact, by replacing the estimate  for $\|\bvarphi-\bvarphi_h\|_{1,\Omega}$ given by
\eqref{Cea-Estimate-Transport} into the right-hand side of \eqref{Estimate-error-F}, and
using the fact that $\|\bu\|_{1, \Omega}\leq r\,\CA|\alp||\bg|$ (cf. \eqref{C-D-D}) and
$\|\bt\|_{\varepsilon, \Omega}\leq \widehat{C}(r)|\alp||\bg|\|\bvarphi\|_{0,\Omega}$
(cf. \eqref{Regularity-estimate}) along with $\|\bvarphi\|_{1,\Omega} \leq r$, we find that
\begin{eqnarray*}
\|\vbx-\vbx_h\|_{\mathbb{H}} & \leq & C_{\tt{ST}}\,\mathrm{dist}\big(\vbx,\mathbb{H}_h\big)
\, +\, C_{\tt{ST}}\widetilde{C}_{\tt{ST}}\,\Big\{|\alp|\,|\bg| + \|\bt\|_{\varepsilon, \Omega}\Big\}\,\|\bvarphi-\mathcal{I}^{SZ}_{h}(\bvarphi)\|\nonumber\\
& & +\ C_{\tt{ST}}\,\Big\{|\alp|\,|\bg| + \|\bt\|_{\varepsilon, \Omega} + \|\bu\|_{1, \Omega}\Big\}\|\bu-\bu_{h}\|_{1, \Omega}\nonumber\\[1ex]
& \leq & C_{\tt{ST}}\,\mathrm{dist}\big(\vbx,\mathbb{H}_h\big)
\, +\, C_0\,\|\bvarphi-\mathcal{I}^{SZ}_{h}(\bvarphi)\|_{1, \Omega}\, +\, C_1\,|\alp|\,|\bg|\,\|\bu-\bu_{h}\|_{1, \Omega}\,,\label{Eq-4}
\end{eqnarray*}
where $C_0 := C_{\tt{ST}}\widetilde{C}_{\tt{ST}}(1+r\widehat{C}(r))\,|\alp|\,|\bg|$
and $C_1 := C_{\tt{ST}}(1 + r\widehat{C}(r) + r\CA)$.
In this way, assuming that the data $\alp$ and $\bg$ satisfy that
\begin{equation}\label{Data-constraint-2}
C_1\,|\alp|\,|\bg|\ \leq\ \frac{1}{2}\,,
\end{equation}
we can conclude that
\begin{equation}\label{Eq-5}
\|\vbx-\vbx_h\|_{\mathbb{H}}\ \leq\ 2\,C_{\tt{ST}}\,\mathrm{dist}\big(\vbx, \mathbb{H}_h\big)\, +\, 2\,C_0\,\|\bvarphi-\mathcal{I}^{SZ}_{h}(\bvarphi)\|_{1, \Omega}.
\end{equation}

Consequently, we now can establish the following main result.
\begin{theorem}\label{Cea}
Assume that $r$ and the data $\alp$ and $\bg$ are sufficiently small so that
\eqref{Data-constraint-1} and \eqref{Data-constraint-2} hold, respectively.
Then, there exists a positive constant $C^{*}$, independent of $h$, such that
\begin{equation}\label{Cea-Estimate}
\|\vbx-\vbx_h\|_{\mathbb{H}}\, +\, \|\bvarphi-\bvarphi_h\|_{1, \Omega}\ \leq\ C^{*}\,\Big\{\mathrm{dist}\big(\vbx, \mathbb{H}_h\big)\;+\;\|\bvarphi-\mathcal{I}^{SZ}_{h}(\bvarphi)\|_{1, \Omega}\Big\}\,.
\end{equation}
\end{theorem}

\begin{proof}
It follows straightforwardly from the C\'ea estimates \eqref{Eq-5} and \eqref{Cea-Estimate-Transport}.
\end{proof}

\medskip
In order to provided the result concerning to the theoretical rate of convergence
of \eqref{augmented-mixed-primal-h}, we recall from \cite{gatica-2014}, the
approximation properties of the specific finite element subspaces introduced in
Section \ref{FEM-specific-spaces}.

\noindent\noindent
($\mathbf{AP}_{h}^{\bt}$) There exists $c>0$, independent of $h$, such that for each $s\in (0, k+1]$, and for each $\bt\in \mathbb{H}^{s}(\Omega)\cap \mathbb{L}^{2}_{\tt{tr}}(\Omega)$, there holds 
\begin{equation*}\label{App-t}
\mathrm{dist}\big(\bt, \mathbb{H}^{\bt}_{h}\big)\;\leq \; c\,h^{s}\,\|\bt\|_{s,\Omega}\,.
\end{equation*}
($\mathbf{AP}_{h}^{\bsig}$) There exists $c>0$, independent of $h$, such that for each $s\in (0, k+1]$, and for each $\bsig\in \mathbb{H}^{s}(\Omega)\cap \mathbb{H}_{0}(\bdiv;\Omega)$, with $\bdiv\bsig \in \mathbf{H}^{s}(\Omega)$,  there holds 
\begin{equation*}\label{App-sig}
\mathrm{dist}\big(\bsig, \mathbb{H}^{\bsig}_{h}\big)\;\leq \; c\,h^{s}\,\big\{\|\bsig\|_{s, \Omega}\, +\, \|\bdiv \bsig\|_{s, \Omega}\big\}\,.
\end{equation*}
($\mathbf{AP}_{h}^{\bu}$) There exists $c>0$, independent of $h$, such that for each $s\in (0, k+1]$, and for each $\bu\in \mathbf{H}^{s+1}(\Omega)$,  there holds 
\begin{equation*}\label{App-u}
\mathrm{dist}\big(\bu, \mathbf{H}^{\bu}_{h}\big)\;\leq \; c\,h^{s}\,\|\bu\|_{s+1, \Omega}\,.
\end{equation*}
Finally,  thanks to the approximation property of $\mathcal{I}^{SZ}_{h}$ given in Lemma \ref{Scott-Zhang},  there exists $c>0$, independent of $h$, such that for each $s\in (0, k+1]$, and for each $\bvarphi \in \mathbf{H}^{s+1}(\Omega)$,  there holds 
\begin{equation*}\label{App-varphi}
\|\bvarphi-\mathcal{I}^{SZ}_{h}(\bvarphi)\|_{1, \Omega}\;\leq \; c\,h^{s}\,\|\bvarphi\|_{s+1, \Omega}\,.
\end{equation*}
Therefore,  from the C\'ea estimate \eqref{Cea-Estimate}, employing the aforementioned approximation properties, we can establish the following result.
\begin{theorem}\label{theoreticalRates}
In addition to the hypotheses of Theorems {\rm\ref{Main-result}}, {\rm\ref{Main-result-h}}
and {\rm\ref{Cea}}, assume that there exists $s>0$ such that $\bt\in \mathbb{H}^{s}(\Omega), \bsig\in \mathbb{H}^{s}(\Omega), \bdiv \,\bsig \in \mathbf{H}^{s}(\Omega), \bu\in \mathbf{H}^{1+s}(\Omega), \bvarphi\in \mathbf{H}^{1+s}(\Omega)$. Then, there exists a positive constant $C$, independent of $h$, such that with the finite element subspaces defined by \eqref{FE-1}, \eqref{FE-2}, \eqref{FE-3}, \eqref{FE-4}, there holds
\begin{equation*}\label{Final-estimate}
\|\vbx-\vbx_h\|_{\mathbb{H}}+\|\bvarphi-\bvarphi_h\|_{1, \Omega}\;\leq \; C\, h^{\min\{s, k+1\}}\Big\{\|\bt\|_{s, \Omega}+\|\bsig\|_{s,\Omega}+\|\bdiv \,\bsig\|_{s, \Omega}+\|\bu\|_{1+s, \Omega}+\|\bvarphi\|_{1+s, \Omega}\Big\}.
\end{equation*}
\end{theorem}

\section{Numerical results}\label{sec:Nu-Ex}

In this section we present three numerical experiments illustrating the performance
of our semi-augmented mixed-primal finite element scheme \eqref{augmented-mixed-primal-h},
and confirming the rates of convergence provided by Theorem \ref{theoreticalRates}.
More precisely, we take the stabilization parameters $\kappa_1$, $\kappa_2$ and $\kappa_3$
as in \eqref{Parameters-maximized}, which satisfies the assumption of Lemma \ref{Lemma_Flow-h}.
In addition, the zero integral mean condition for tensors in the space \eqref{FE-2}
is imposed via a real Lagrange multiplier. In turn, the nonlinear algebraic systems
obtained are solved by the fixed-point method with a tolerance of $10^{-6}$,
along with the Newton method for approximate the solution of \eqref{augmented-mixed-primal-1-h}
in each fixed-point's iteration. We take as initial guess the solution of a similar
linear problem (in particular, satisfying the boundary conditions for $\bu_h$ and $\bvarphi_h$).
The numerical results presented below were obtained using a \Cpp code, where the corresponding
linear systems arising from \eqref{augmented-mixed-primal-1-h} are solved using the BiCGSTAB
method, whereas for \eqref{augmented-mixed-primal-2-h} we employ the Conjugate Gradient
method as the main solver. Finally, in all experiments we let $\bg = (0,-1)^{\mathrm{t}}$
be the gravitational force, and utilizing structure triangulations of the corresponding
domain in 2D. Furthermore, for the first two examples we consider polynomial degrees
$k\in\{0,1,2\}$, whereas we only use $k = 0$ in the last example.

We now introduce some additional notation. The individual errors are denoted by:
$$\mathtt{e}(\bt)\ :=\ \|\bt-\bt_h\|_{0,\Omega}\,,\qquad\mathtt{e}(\bsig)\ :=\ \|\bsig-\bsig_h\|_{\mathbf{div},\Omega}\,,\qquad\mathtt{e}(\bu)\ :=\ \|\bu-\bu_h\|_{1,\Omega}\,,$$
$$\mathtt{e}(\bvarphi)\ :=\ \|\bvarphi-\bvarphi_h\|_{1,\Omega}\qquad\text{and}\qquad \mathtt{e}(p)\ :=\ \|p-p_h\|_{1,\Omega}\,,$$
where, according to \eqref{eqn:pressure} and \eqref{eqn:def-c0}, $p_h$ can be computed as:
$$p_h\ =\ -\dfrac{1}{n}\Big[\tr(\bsig_h\, +\, (\bu_h\otimes\bu_h))\Big]\, +\, \dfrac{1}{n |\Omega|}\|\bu_h\|_{0,\Omega}^2\,.$$
On the other hand as is usual, we let $\mathtt{r}(\cdot)$ be the experimental rate of convergence given by
$$\mathtt{r}(\cdot)\ :=\ \frac{\mathrm{log}(\mathtt{e}(\cdot)\,/\,\mathtt{e}^{\prime}(\cdot))}{\mathrm{log}(h\,/\,h^{\prime})}\,,$$
where $\mathtt{e}$ and $\mathtt{e}^{\prime}$ denote errors computed on two consecutive meshes of sizes $h$ and $h^{\prime}$, respectively.
In addition, $N$ stands for the total number of degrees of freedom (unknowns) of \eqref{augmented-mixed-primal-h},
that is,
\begin{eqnarray*}
N & := & 4\times\{\text{number of nodes in }\T\}\, +\, \big\{2(k+1) + 4k\big\}\times\{\text{number of edges in } \T\}\\[1ex]
& & +\ \big\{(k+1)(k+2) + 2k(k+1) + 2k(k-1)\big\}\times\{\text{number of elements in } \T\}\ +\ 1\,.
\end{eqnarray*}

\noindent {\bf Example 1.} 
We first consider the square $\Omega = (0,1)^2$, and set $\GD = \left\{(s,0),(s,1)\in\text{\rm R}^2\,:\, 0 \leq s\leq 1\right\}$,\,
$\GN = \Gamma\setminus\GD$,\, $\gamma = 0.1$,\, $\nu(\bx) = (x_1^2 + x_2^2 + 1)^{-1}$
(here $\nu_1=1$ and $\nu_2 = 2$),\, $\alp = (0.5,1.5)^{\mathrm{t}}$,\, the thermal conductivity
$\mathbb{K} = 2\,\mathbb{I}$, and adequately manufacture the data so that the exact solution is
given by the smooth functions
$$\bu(\bx)\ =\ \left(\begin{array}{r}
-\sin^2(2\pi x_1)\sin(4\pi x_2)\\[1ex]
 \sin(4\pi x_1)\sin^2(2\pi x_2)
\end{array}\right)\,,\qquad p(\bx)\ =\ \cos(x_1)\cos(x_2)\, -\, \sin^2(1)\,,$$
and
$$\bvarphi(\bx)\ =\ \left(\begin{array}{c}
x_1x_2\\[1ex] \exp(x_1 + x_2)
\end{array}\right)\,,$$
for all $\bx := (x_1,x_2)^{\mathrm{t}}\in\Omega$. In Table \ref{tab:exa01}, we summarize the
convergence history of the finite element scheme \eqref{augmented-mixed-primal-h} as applied
to Example 1. We notice there that the rate of convergence $O(h^{k+1})$ predicted by Theorem
\ref{theoreticalRates} is attained by all the unknowns.

\begin{table}[h!t]\centering\footnotesize
\begin{tabular}{|c|c|c|c@{\hspace{1.5ex}}c|c@{\hspace{1.5ex}}c|c@{\hspace{1.5ex}}c|c@{\hspace{1.5ex}}c|c@{\hspace{1.5ex}}c|}\hline
$k$ & $h$ & $N$ & $\mathtt{e}(\bt)$ & $\mathtt{r}(\bt)$ & $\mathtt{e}(\bsig)$ & $\mathtt{r}(\bsig)$ & $\mathtt{e}(\bu)$ & $\mathtt{r}(\bu)$ & $\mathtt{e}(\bvarphi)$ & $\mathtt{r}(\bvarphi)$ & $\mathtt{e}(p)$ & $\mathtt{r}(p)$\\ \hline
  & 0.0404 & 17575 & 7.98e-01 & $--$ & 2.68e+00 & $--$ & 1.37e+00 & $--$ & 4.13e-02 & $--$ & 1.09e-01 & $--$\\
  & 0.0314 & 28895 & 6.05e-01 & 1.10 & 2.09e+00 & 1.00 & 1.05e+00 & 1.05 & 3.20e-02 & 1.02 & 8.53e-02 & 0.98\\
  & 0.0257 & 43015 & 4.87e-01 & 1.09 & 1.71e+00 & 1.00 & 8.51e-01 & 1.05 & 2.61e-02 & 1.01 & 6.99e-02 & 0.99\\
  & 0.0218 & 59935 & 4.07e-01 & 1.07 & 1.45e+00 & 1.00 & 7.16e-01 & 1.04 & 2.20e-02 & 1.01 & 5.92e-02 & 1.00\\
0 & 0.0189 & 79655 & 3.50e-01 & 1.06 & 1.25e+00 & 1.00 & 6.18e-01 & 1.03 & 1.91e-02 & 1.01 & 5.13e-02 & 1.00\\
  & 0.0129 & 170725 & 2.35e-01 & 1.04 & 8.54e-01 & 1.00 & 4.18e-01 & 1.02 & 1.30e-02 & 1.00 & 3.49e-02 & 1.00\\
  & 0.0094 & 316805 & 1.71e-01 & 1.02 & 6.27e-01 & 1.00 & 3.05e-01 & 1.01 & 9.52e-03 & 1.00 & 2.56e-02 & 1.00\\
  & 0.0071 & 562405 & 1.28e-01 & 1.01 & 4.70e-01 & 1.00 & 2.29e-01 & 1.01 & 7.14e-03 & 1.00 & 1.92e-02 & 1.00\\
  & 0.0057 & 878005 & 1.02e-01 & 1.00 & 3.76e-01 & 1.00 & 1.83e-01 & 1.01 & 5.71e-03 & 1.00 & 1.54e-02 & 1.00\\ \hline
  & 0.0404 & 59645 & 5.37e-02 & $--$ & 1.77e-01 & $--$ & 8.81e-02 & $--$ & 1.45e-04 & $--$ & 7.38e-03 & $--$\\
  & 0.0314 & 98285 & 3.20e-02 & 2.07 & 1.07e-01 & 1.99 & 5.28e-02 & 2.04 & 8.48e-05 & 2.12 & 4.53e-03 & 1.94\\
1 & 0.0257 & 146525 & 2.12e-02 & 2.05 & 7.18e-02 & 2.00 & 3.51e-02 & 2.03 & 5.61e-05 & 2.06 & 3.05e-03 & 1.96\\
  & 0.0218 & 204365 & 1.51e-02 & 2.04 & 5.14e-02 & 2.00 & 2.51e-02 & 2.02 & 4.00e-05 & 2.03 & 2.20e-03 & 1.97\\
  & 0.0189 & 271805 & 1.13e-02 & 2.03 & 3.86e-02 & 2.00 & 1.88e-02 & 2.02 & 3.00e-05 & 2.01 & 1.66e-03 & 1.98\\
  & 0.0129 & 583445 & 5.19e-03 & 2.02 & 1.80e-02 & 2.00 & 8.70e-03 & 2.01 & 1.39e-05 & 2.01 & 7.75e-04 & 1.98\\ \hline
  & 0.0404 & 126215 & 2.70e-03 & $--$ & 8.82e-03 & $--$ & 3.76e-03 & $--$ & 3.25e-07 & $--$ & 3.45e-04 & $--$\\
  & 0.0314 & 208175 & 1.26e-03 & 3.04 & 4.15e-03 & 2.99 & 1.76e-03 & 3.03 & 1.47e-07 & 3.17 & 1.63e-04 & 2.99\\
2 & 0.0257 & 310535 & 6.87e-04 & 3.03 & 2.28e-03 & 3.00 & 9.58e-04 & 3.02 & 7.96e-08 & 3.05 & 8.92e-05 & 2.99\\
  & 0.0218 & 433295 & 4.15e-04 & 3.01 & 1.36e-03 & 3.08 & 5.79e-04 & 3.02 & 4.79e-08 & 3.04 & 5.40e-05 & 3.01\\
  & 0.0189 & 576455 & 2.70e-04 & 3.00 & 8.83e-04 & 3.03 & 3.76e-04 & 3.01 & 3.10e-08 & 3.03 & 3.51e-05 & 3.00\\ \hline
\end{tabular}
\caption{History of convergence for Example 1.}\label{tab:exa01}
\end{table}

\noindent {\bf Example 2.} 
Next, we adapt \cite[Example 3]{cg-vem-stokes}, and consider the $L$-shaped domain
$\Omega = (-1,1)^2\,\setminus\,[0,1]^2$, and set $\GN = \left\{(s,0),(0,s)\in\text{\rm R}^2\,:\, 0 \leq s\leq 1\right\}$,\,
$\GD = \Gamma\setminus\bar{\Gamma}_N$,\, $\gamma = 10^{-3}$,\, $\nu(\bx) = 1 + \exp(-x_1^2)$
(once again $\nu_1=1$ and $\nu_2 = 2$),\, $\alp = (1,0.5)^{\mathrm{t}}$,\,
$\mathbb{K} = \left(\begin{smallmatrix}1 & \;0\\[0.5ex] 0 & \;2\end{smallmatrix}\right)$, and adequately manufacture the data so that the
exact solution is given by the smooth functions
$$\bu(\bx)\ =\ \left(\begin{array}{r}
 x_2^2\\[1ex]
-x_1^2
\end{array}\right)\,,\qquad p(\bx)\ =\ (x_1^2\, +\, x_2^2)^{1/3}\, -\, p_0\,,\qquad\text{and}\qquad\bvarphi(\bx)\ =\ \left(\begin{array}{c}
\exp(-x_1^2 - x_2^2)\\[1ex] \exp(-x_1x_2)
\end{array}\right)\,,$$
for all $\bx := (x_1,x_2)^{\mathrm{t}}\in\Omega$, where $p_0\in\mathrm{R}$ is such that
$\int_{\Omega}p = 0$ holds ($p_0\approx \text{8.211056552903396e-01}$). In addition,
we remark here that the partial derivatives of $p$, and hence, in particular $\bdiv(\bsig)$,
are singular at the origin. Indeed, according to the power $1/3$, there holds $\bsig\in\mathbb{H}^{5/3-\varepsilon}(\Omega)$
and $\bdiv(\bsig)\in\textbf{H}^{2/3-\varepsilon}(\Omega)$ for each $\varepsilon > 0$.
In fact, in Table \ref{tab:exa02} we present the corresponding convergence history of
Example 2, where, as predicted in advance, we note that the orders $O(h^{\min\{k+1,5/3\}})$
and $O(h^{2/3})$ are attained by $(\bt_h, \bu_h)$ and $\bsig_h$,
respectively.
Once again, the rate of convergence predicted by Theorem \ref{theoreticalRates} is attained
by all the unknowns, except for the variable $\bvarphi_h$ that preserves $O(h^{k+1})$.
The foregoing phenomenon could be a special feature of this example.
Furthermore, the results in Example 2 suggest that our approach should certainly be strengthened
with the further incorporation of an adaptive strategy based on a suitable a-posteriori error
estimates. This issue will also be addressed in a forthcoming paper.

\begin{table}[h!t]\centering\footnotesize
\begin{tabular}{|c|c|c|c@{\hspace{1.5ex}}c|c@{\hspace{1.5ex}}c|c@{\hspace{1.5ex}}c|c@{\hspace{1.5ex}}c|c@{\hspace{1.5ex}}c|}\hline
$k$ & $h$ & $N$ & $\mathtt{e}(\bt)$ & $\mathtt{r}(\bt)$ & $\mathtt{e}(\bsig)$ & $\mathtt{r}(\bsig)$ & $\mathtt{e}(\bu)$ & $\mathtt{r}(\bu)$ & $\mathtt{e}(\bvarphi)$ & $\mathtt{r}(\bvarphi)$ & $\mathtt{e}(p)$ & $\mathtt{r}(p)$\\ \hline
  & 0.0707 & 17285 & 5.60e-02 & $--$ & 2.77e-01 & $--$ & 7.07e-02 & $--$ & 1.02e-01 & $--$ & 4.85e-02 & $--$\\
  & 0.0566 & 26855 & 4.48e-02 & 1.00 & 2.23e-01 & 0.97 & 5.66e-02 & 1.00 & 8.15e-02 & 1.00 & 3.88e-02 & 1.00\\
  & 0.0471 & 38525 & 3.73e-02 & 1.00 & 1.87e-01 & 0.97 & 4.72e-02 & 1.00 & 6.80e-02 & 1.00 & 3.23e-02 & 1.00\\
0 & 0.0404 & 52295 & 3.20e-02 & 1.00 & 1.61e-01 & 0.96 & 4.04e-02 & 1.00 & 5.83e-02 & 1.00 & 2.77e-02 & 1.00\\
  & 0.0354 & 68165 & 2.80e-02 & 1.00 & 1.42e-01 & 0.96 & 3.54e-02 & 1.00 & 5.10e-02 & 1.00 & 2.42e-02 & 1.00\\
  & 0.0236 & 152645 & 1.87e-02 & 1.00 & 9.63e-02 & 0.95 & 2.36e-02 & 1.00 & 3.40e-02 & 1.00 & 1.61e-02 & 1.00\\
  & 0.0166 & 305495 & 1.32e-02 & 1.00 & 6.94e-02 & 0.94 & 1.66e-02 & 1.00 & 2.40e-02 & 1.00 & 1.14e-02 & 1.00\\
  & 0.0135 & 465575 & 1.07e-02 & 1.00 & 5.69e-02 & 0.93 & 1.35e-02 & 1.00 & 1.94e-02 & 1.00 & 9.22e-03 & 1.00\\ \hline
  & 0.0707 & 58565 & 2.34e-04 & $--$ & 2.67e-02 & $--$ & 1.35e-04 & $--$ & 1.45e-03 & $--$ & 7.06e-04 & $--$\\
  & 0.0566 & 91205 & 1.56e-04 & 1.82 & 2.30e-02 & 0.68 & 9.20e-05 & 1.71 & 9.30e-04 & 2.00 & 4.62e-04 & 1.90\\
1 & 0.0471 & 131045 & 1.12e-04 & 1.81 & 2.03e-02 & 0.67 & 6.75e-05 & 1.69 & 6.46e-04 & 2.00 & 3.27e-04 & 1.89\\
  & 0.0404 & 178085 & 8.47e-05 & 1.81 & 1.83e-02 & 0.67 & 5.21e-05 & 1.69 & 4.75e-04 & 2.00 & 2.44e-04 & 1.89\\
  & 0.0354 & 232325 & 6.66e-05 & 1.80 & 1.68e-02 & 0.67 & 4.16e-05 & 1.68 & 3.63e-04 & 2.00 & 1.90e-04 & 1.88\\
  & 0.0236 & 521285 & 3.23e-05 & 1.79 & 1.28e-02 & 0.67 & 2.11e-05 & 1.68 & 1.62e-04 & 2.00 & 8.93e-05 & 1.86\\ \hline
  & 0.0707 & 123845 & 3.00e-05 & $--$ & 1.52e-02 & $--$ & 3.29e-05 & $--$ & 1.64e-05 & $--$ & 8.29e-05 & $--$\\
  & 0.0566 & 193055 & 2.07e-05 & 1.67 & 1.31e-02 & 0.67 & 2.27e-05 & 1.67 & 8.39e-06 & 3.01 & 5.70e-05 & 1.68\\
2 & 0.0471 & 277565 & 1.52e-05 & 1.67 & 1.16e-02 & 0.67 & 1.67e-05 & 1.67 & 4.85e-06 & 3.01 & 4.20e-05 & 1.67\\
  & 0.0404 & 377375 & 1.18e-05 & 1.67 & 1.05e-02 & 0.67 & 1.29e-05 & 1.67 & 3.05e-06 & 3.00 & 3.25e-05 & 1.67\\
  & 0.0354 & 492485 & 9.44e-06 & 1.67 & 9.60e-03 & 0.67 & 1.04e-05 & 1.67 & 2.04e-06 & 3.00 & 2.60e-05 & 1.67\\ \hline
\end{tabular}
\caption{History of convergence for Example 2.}\label{tab:exa02}
\end{table}

\noindent {\bf Example 3.} 
Finally, we aim to illustrate the accuracy of our method by considering a case in which the exact solution
is unknown in the a time-dependent approach. More precisely, we add $\partial_t\bu$ and $\partial_t\bvarphi$
to the left-hand side of first and last equations of \eqref{eqn:model}, respectively, which, together with
the boundary conditions \eqref{eqn:bc}, we consider initial conditions
$$\bu(\cdot,0)\,=\,\bu_0\,\hspace{0.2cm}\mbox{in } \Omega
\qquad \mbox{and}\qquad
\bvarphi(\cdot,0)\, =\, \bvarphi_0\,\hspace{0.2cm}\mbox{in } \Omega\,.$$
We remark here that, in similar way of \cite{gss-IMA-2017}, the analysis presented along of this paper can
be extended to this time-dependent problem by employing backward Euler time stepping in order to obtain a
fully-discrete method. On the other hand, for this example we consider once again the unit square
$\Omega = (0,1)^2$, and set $\GD = \left\{(0,s),(1,s)\in\text{\rm R}^2\,:\, 0 \leq s\leq 1\right\}$,\,
$\GN = \Gamma\setminus\GD$,\, $\gamma = 10^{-3}$,\, $\nu(\bx) = 10^{-2}$,\,
$\alp = (1,10)^{\mathrm{t}}$,\, $\mathbb{K} = \left(\begin{smallmatrix}1 & 0\\[0.5ex] 0 & 10^{-1}\end{smallmatrix}\right)$.
The boundary condition is defined as
$$\bvarphi_{\D}(\bx,t)\ =\ \left\{\begin{array}{cl}
(1,1)^{\mathrm{t}}   & \text{if } x_1 = 0\\[1ex]
(-1,-1)^{\mathrm{t}} & \text{if } x_1 = 1
\end{array}\right.$$
for all $\bx := (x_1,x_2)^{\mathrm{t}}\in\Omega$,
whereas the initial conditions are given by
$$\bu_0(\bx)\ =\ \left(\begin{array}{r}
\sin^2(\pi x_1)\sin(2\pi x_2)\\[1ex] -\sin(2\pi x_1)\sin^2(\pi x_2)
\end{array}\right)\qquad\text{and}\qquad \bvarphi_0(\bx)\ =\ \left(\begin{array}{r}
\exp(x_1+x_2)\\[1ex] \exp(x_1-x_2)
\end{array}\right)\,.$$
In addition, for the time stepping technique we use $\triangle t = \frac{1}{50}$.

\begin{figure}[h!t]\centering
\scalebox{0.172}{\begin{minipage}{25cm}
\includegraphics{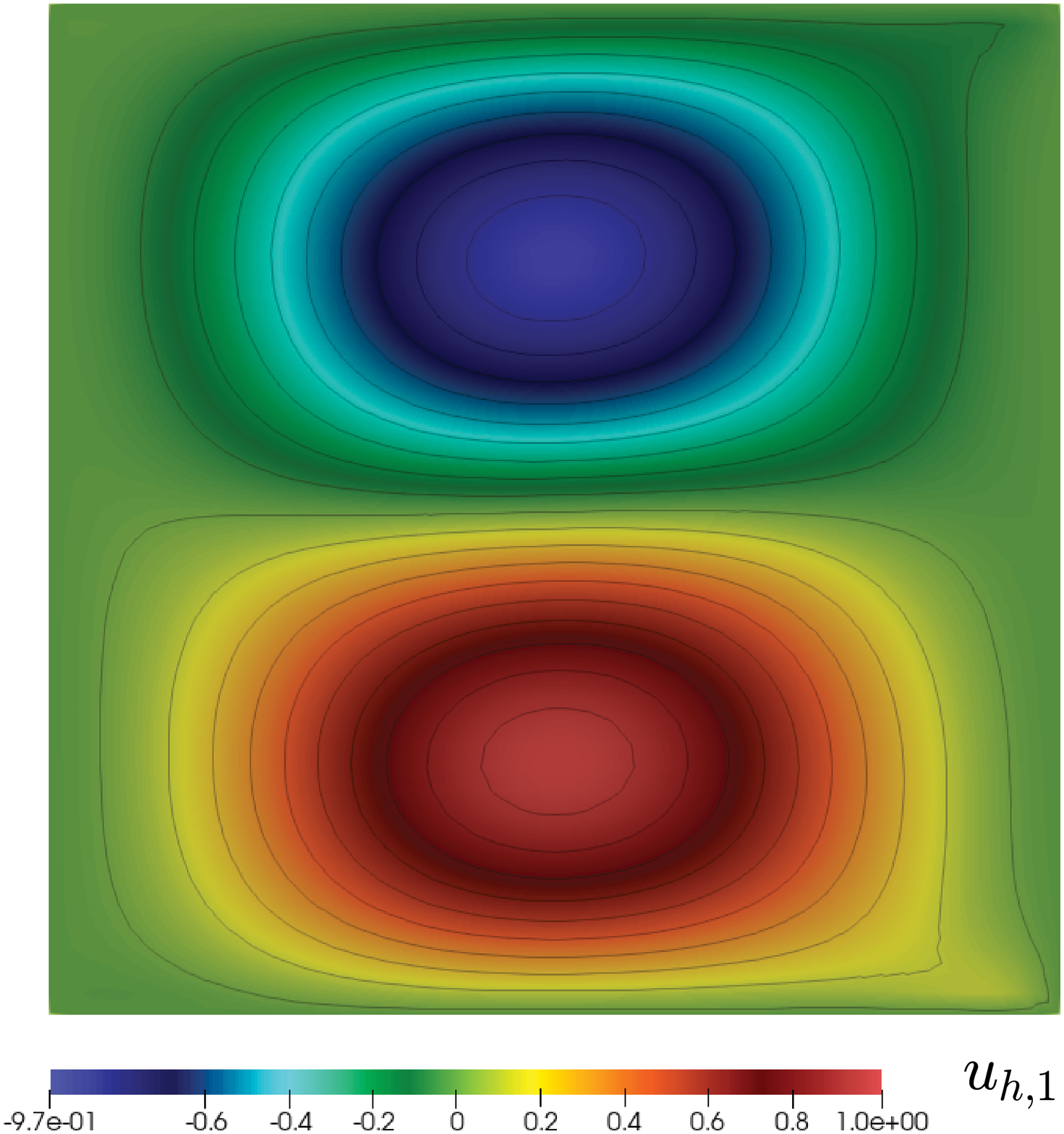}\\[2ex]
\includegraphics{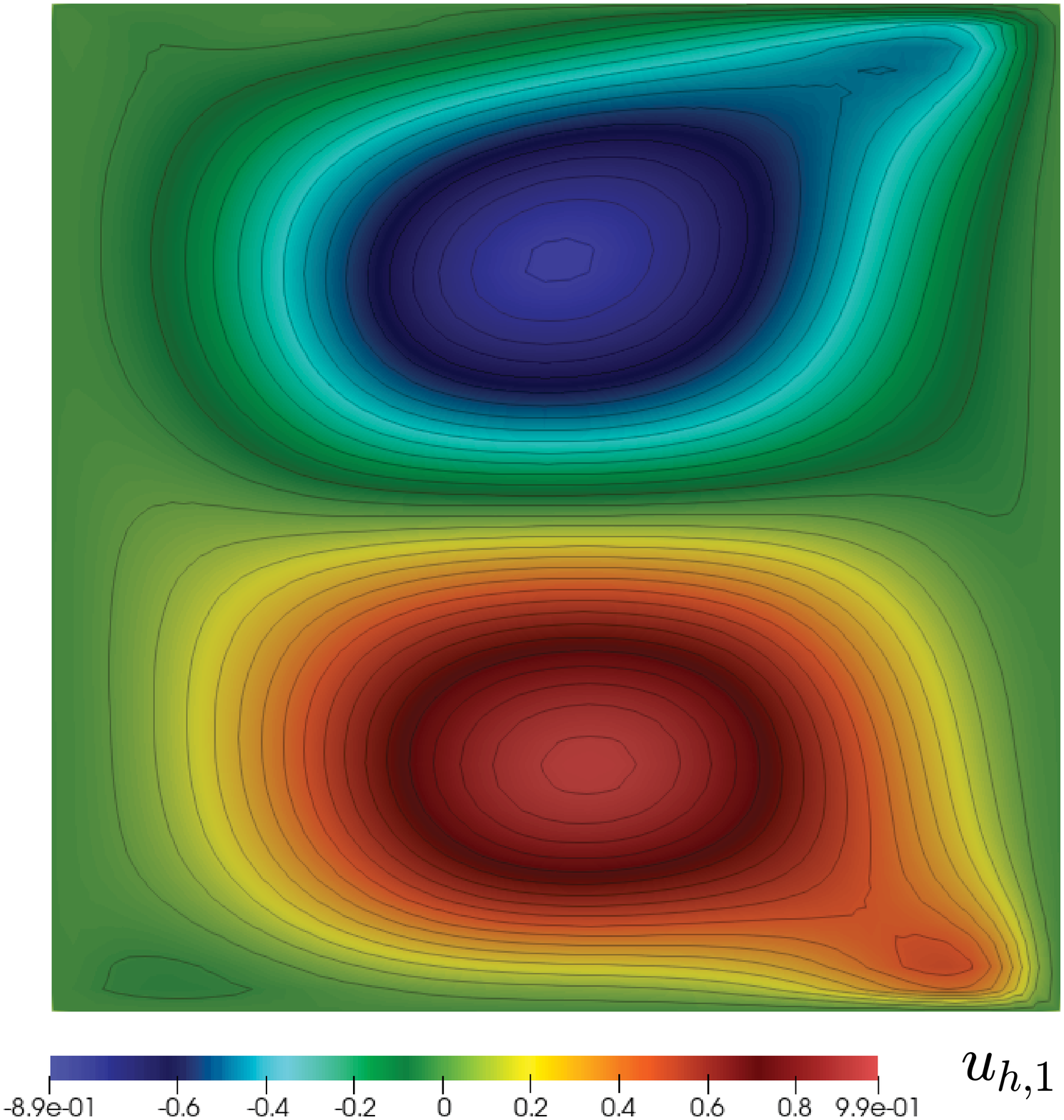}\\[2ex]
\includegraphics{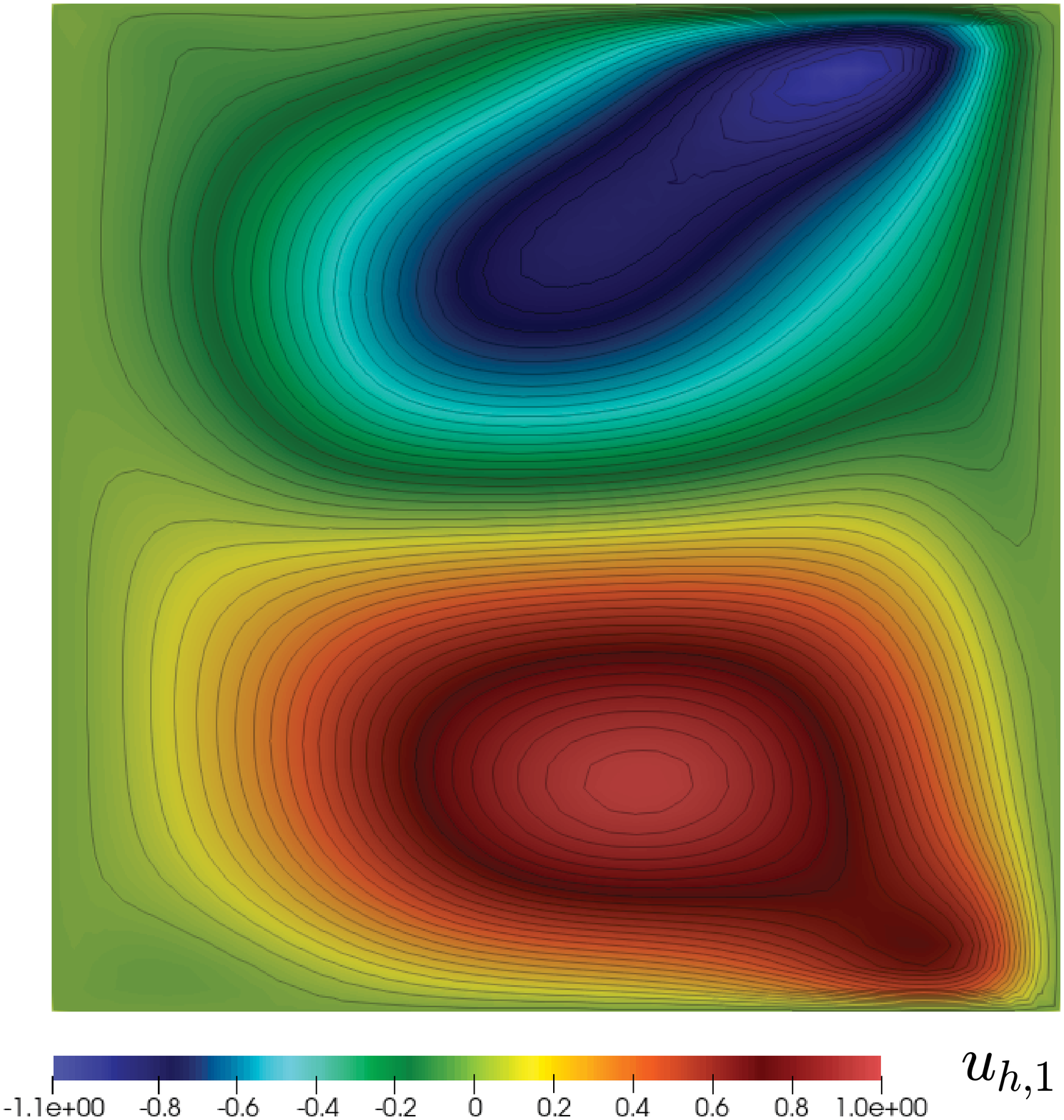}\\[2ex]
\includegraphics{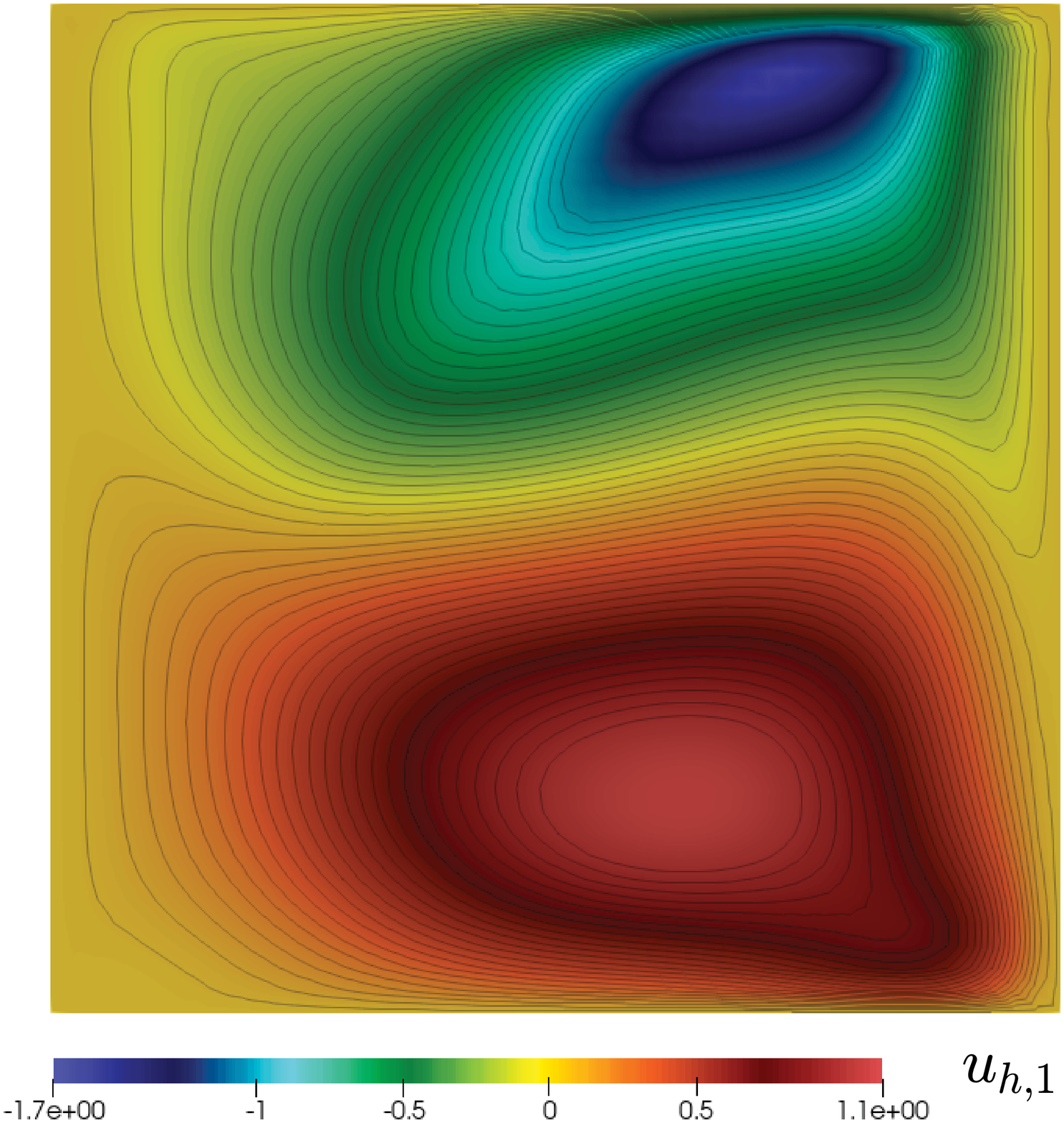}\\[2ex]
\includegraphics{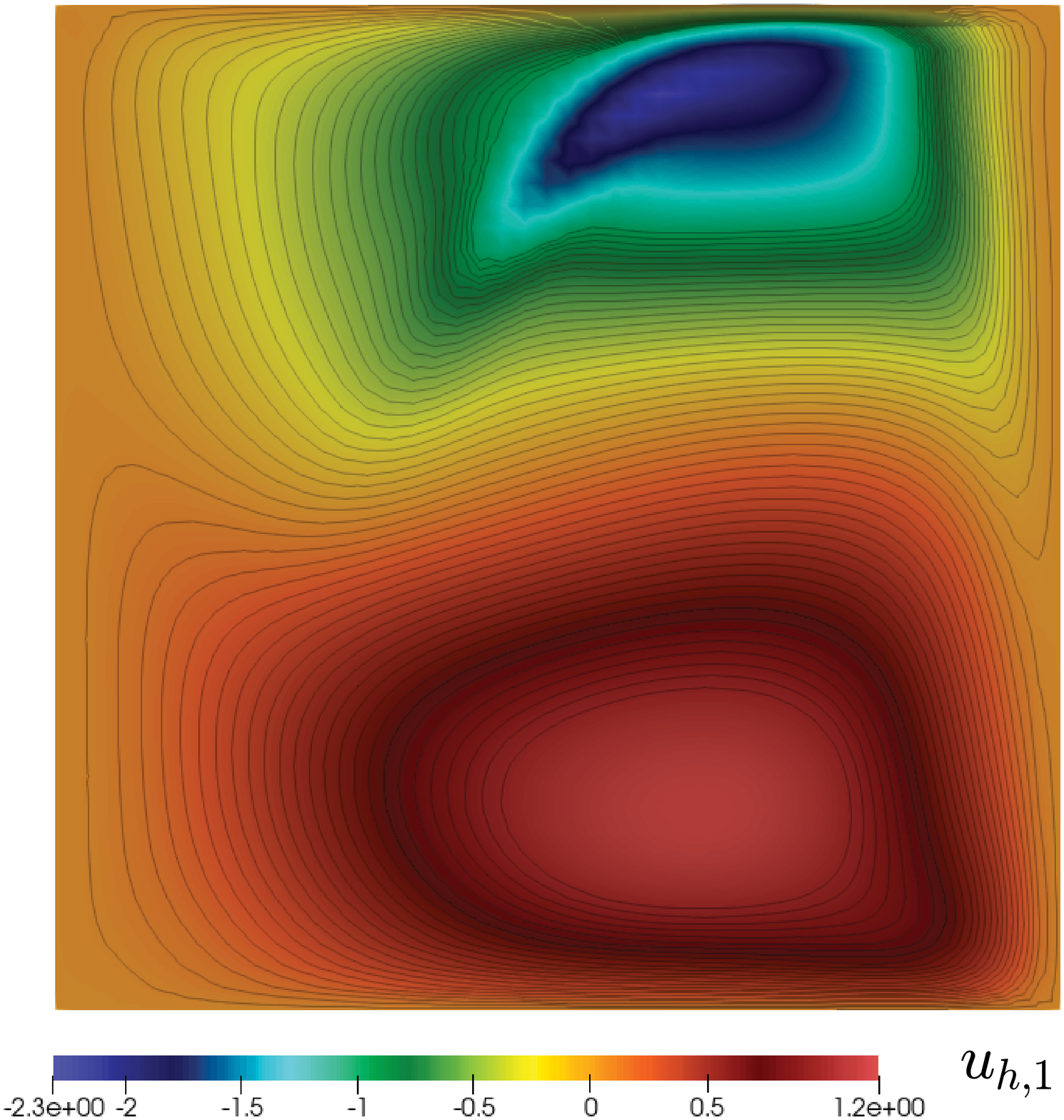}
\end{minipage}\hspace{6cm}\begin{minipage}{25cm}
\includegraphics{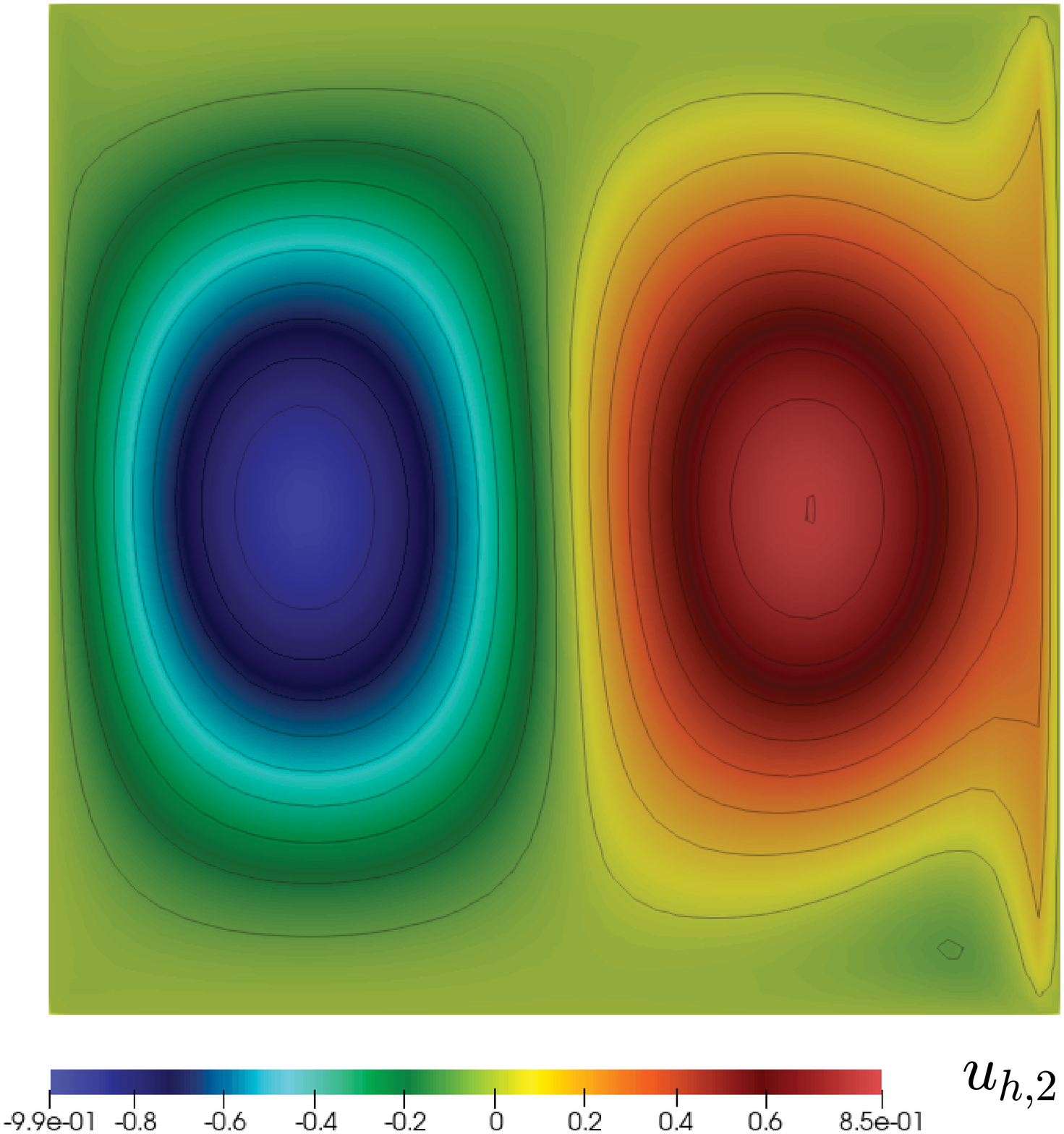}\\[2ex]
\includegraphics{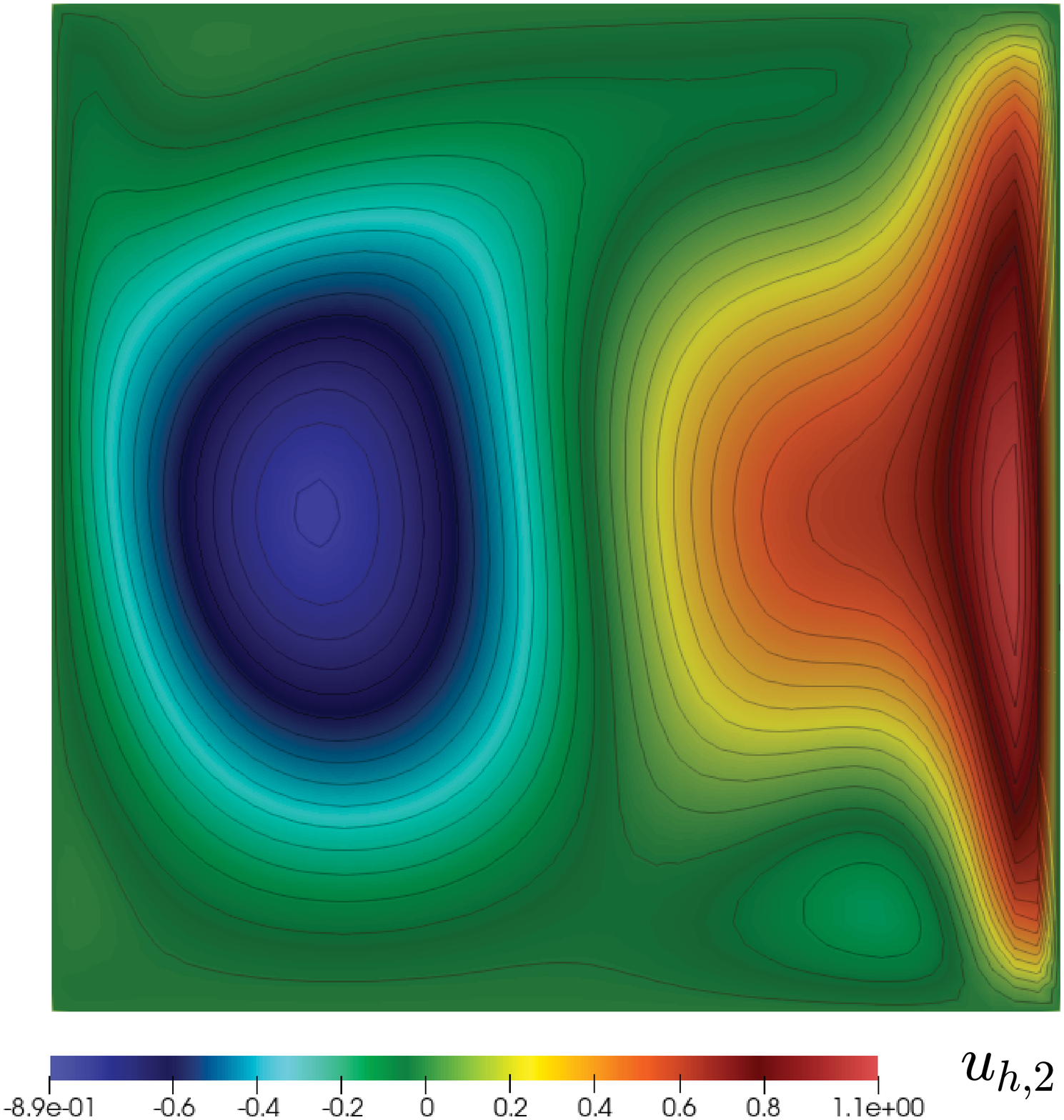}\\[2ex]
\includegraphics{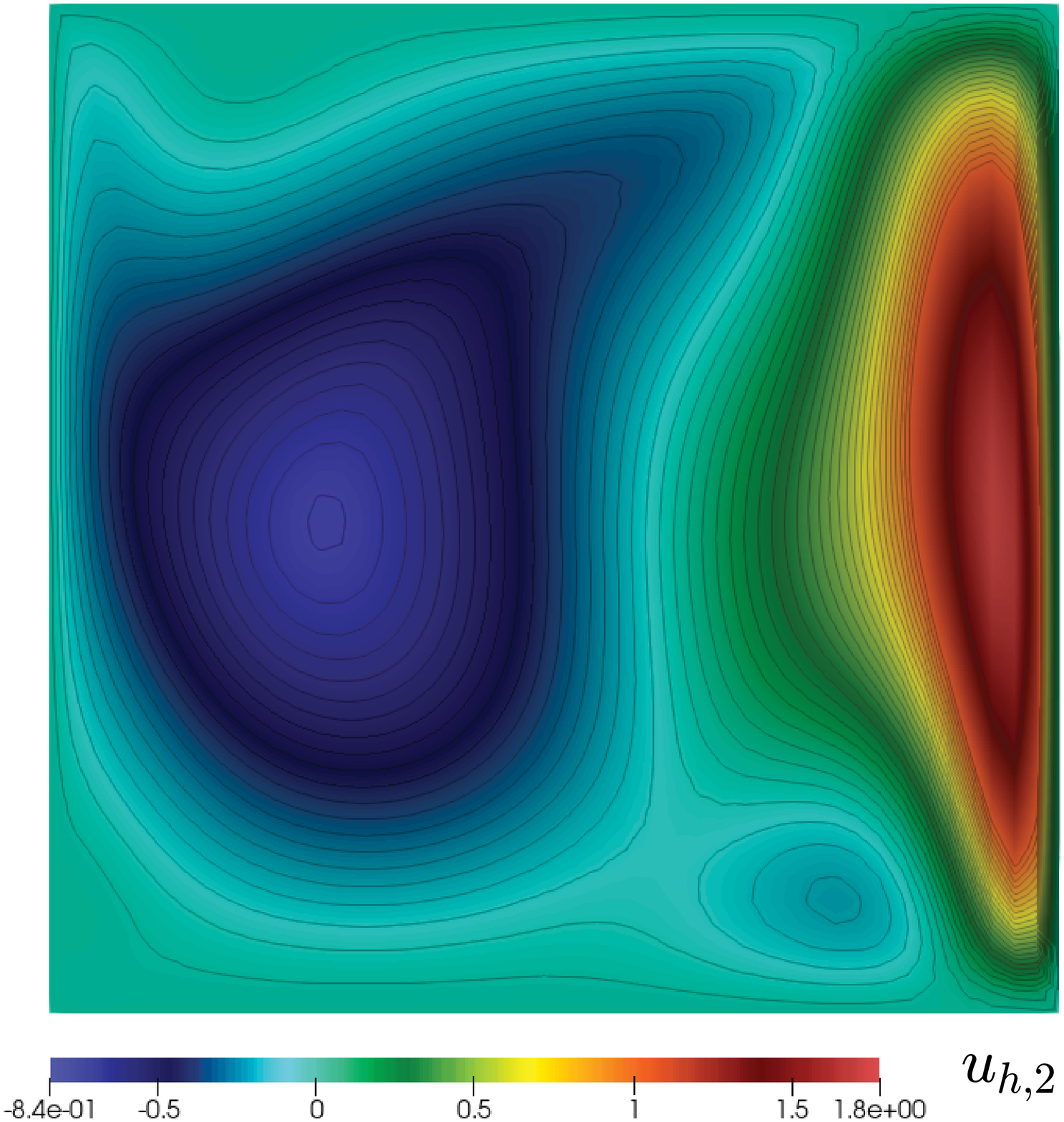}\\[2ex]
\includegraphics{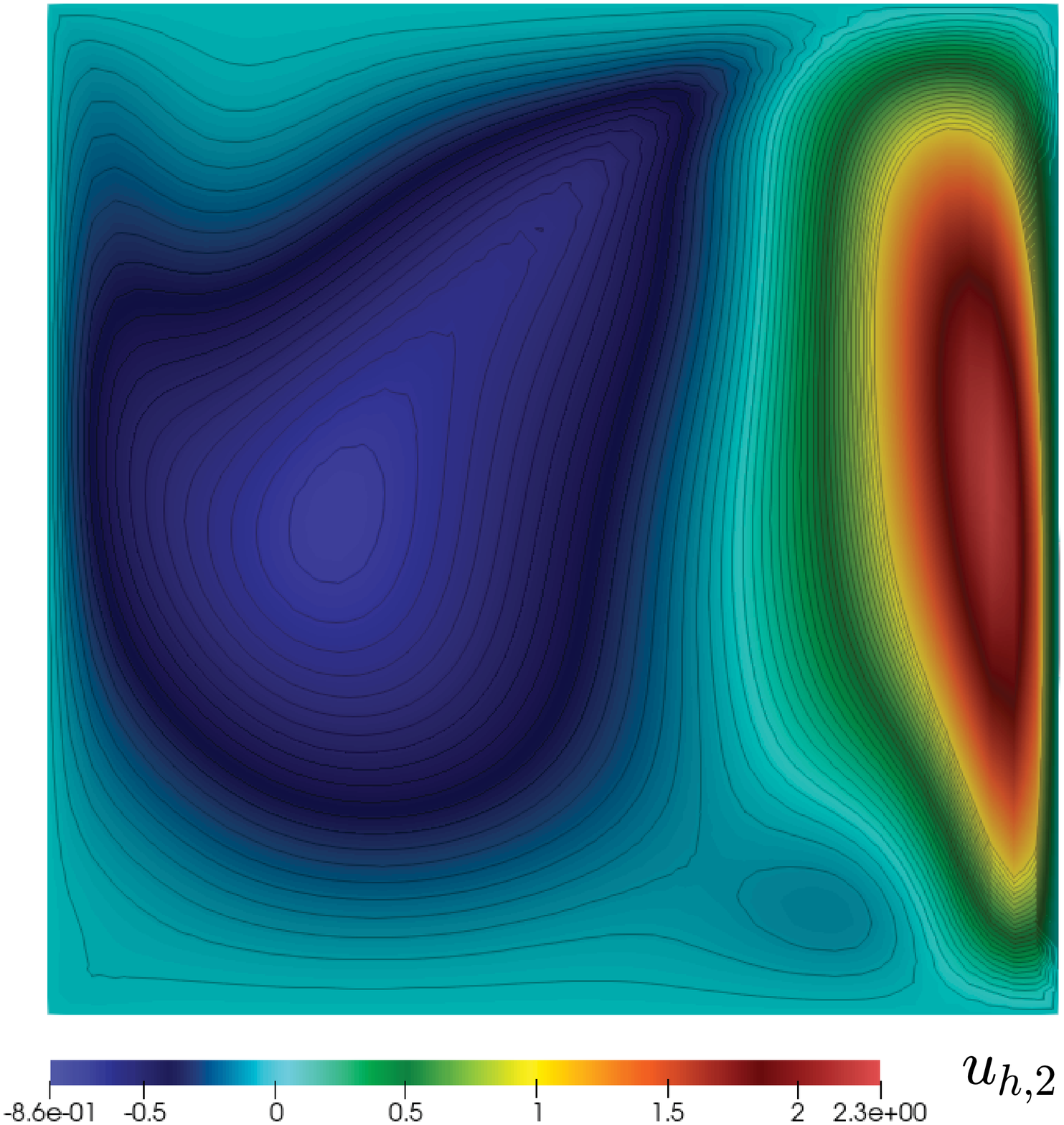}\\[2ex]
\includegraphics{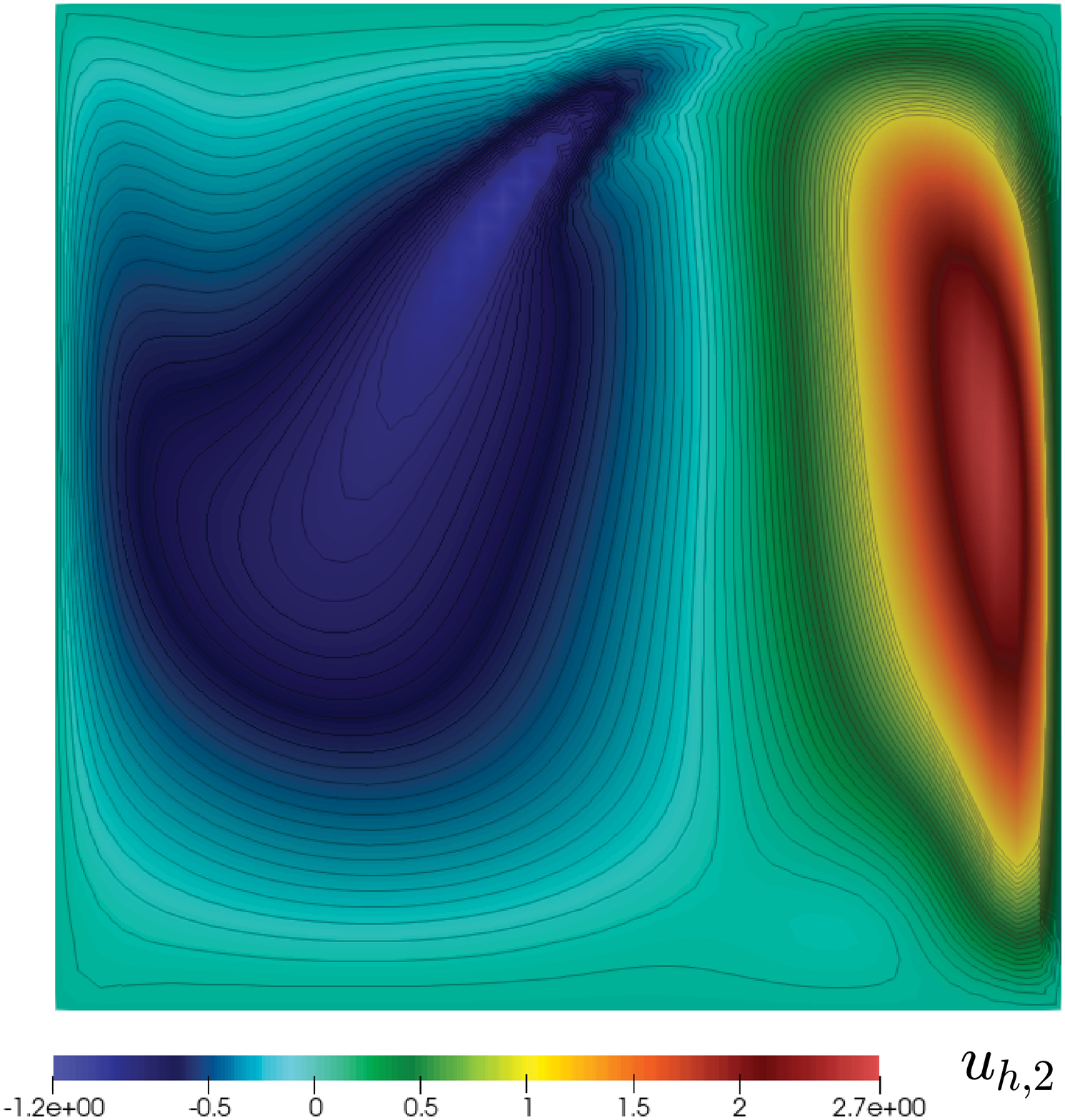}
\end{minipage}}
\caption{Example 3, velocity components $u_{h,1}$ (left) and $u_{h,2}$ (right),
obtained with a fully-discrete time-dependent mixed method with no manufactured
analytical solution using $k = 0$ and $N = 28895$ degrees of freedom. We plot
for $t\in\{t_1, t_5, t_{10}, t_{15}, t_{20}\}$, en each row.}\label{fig:exa03-1}
\end{figure}

\begin{figure}[h!t]\centering
\scalebox{0.172}{\begin{minipage}{25cm}
\includegraphics{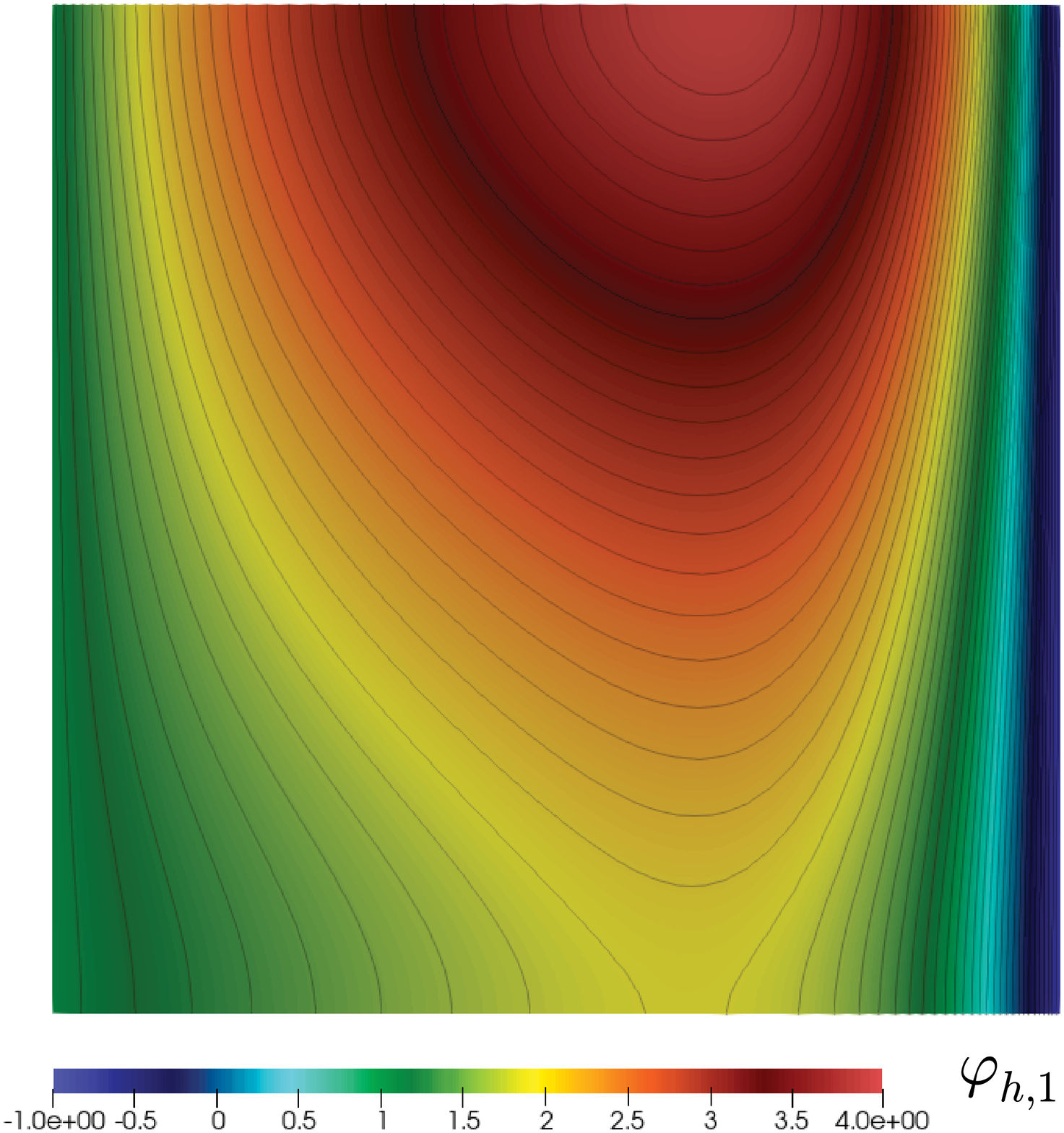}\\[2ex]
\includegraphics{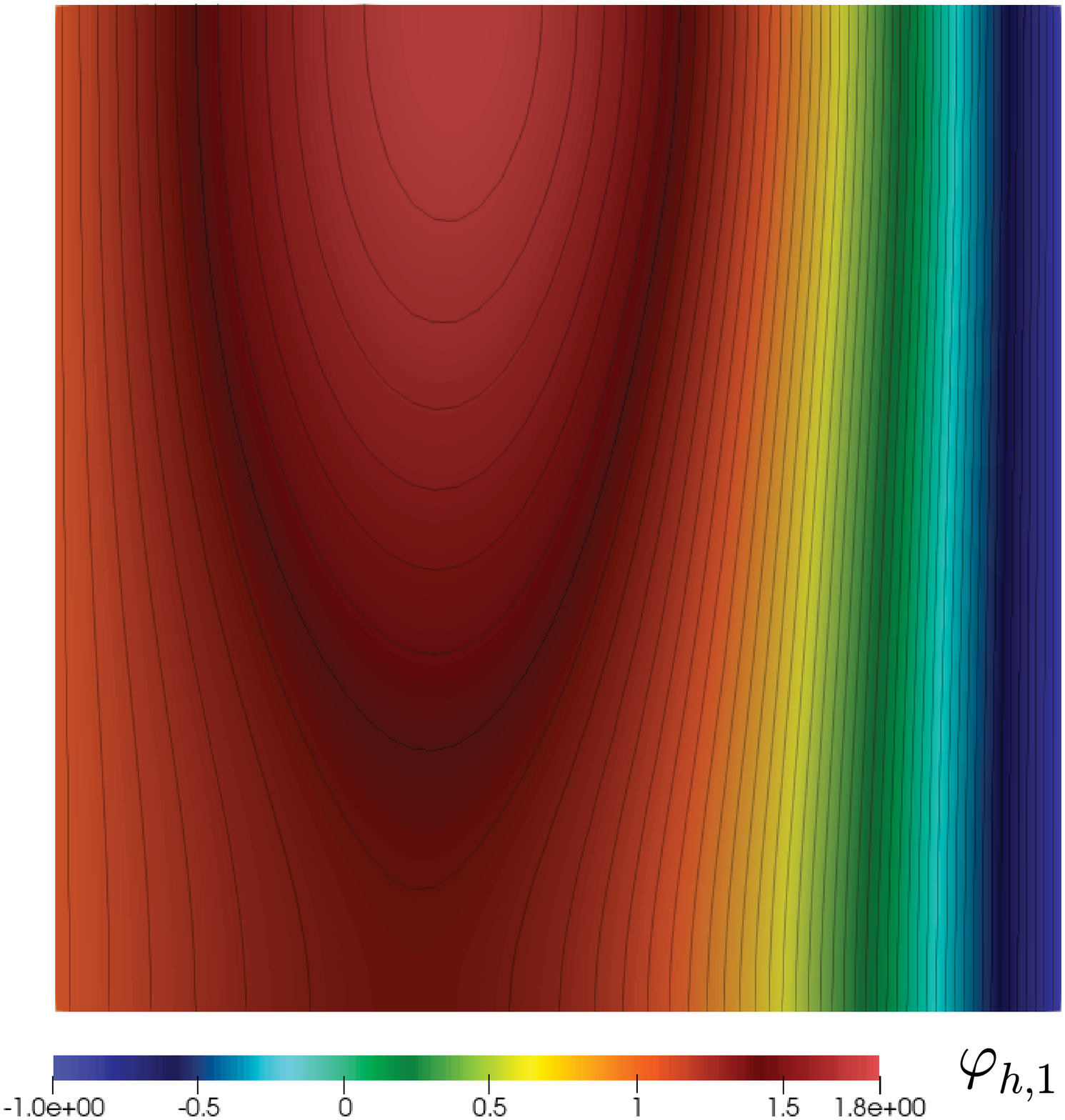}\\[2ex]
\includegraphics{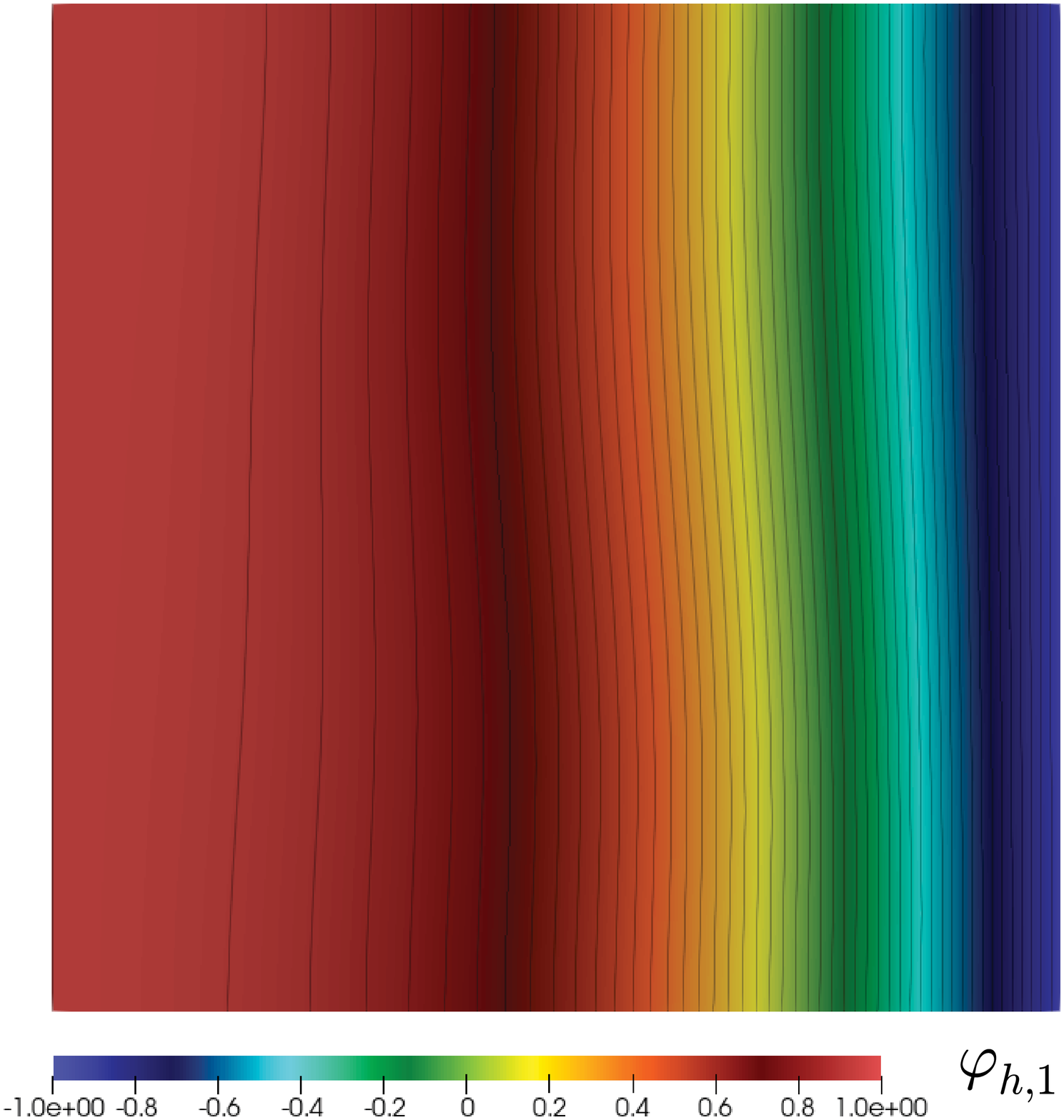}\\[2ex]
\includegraphics{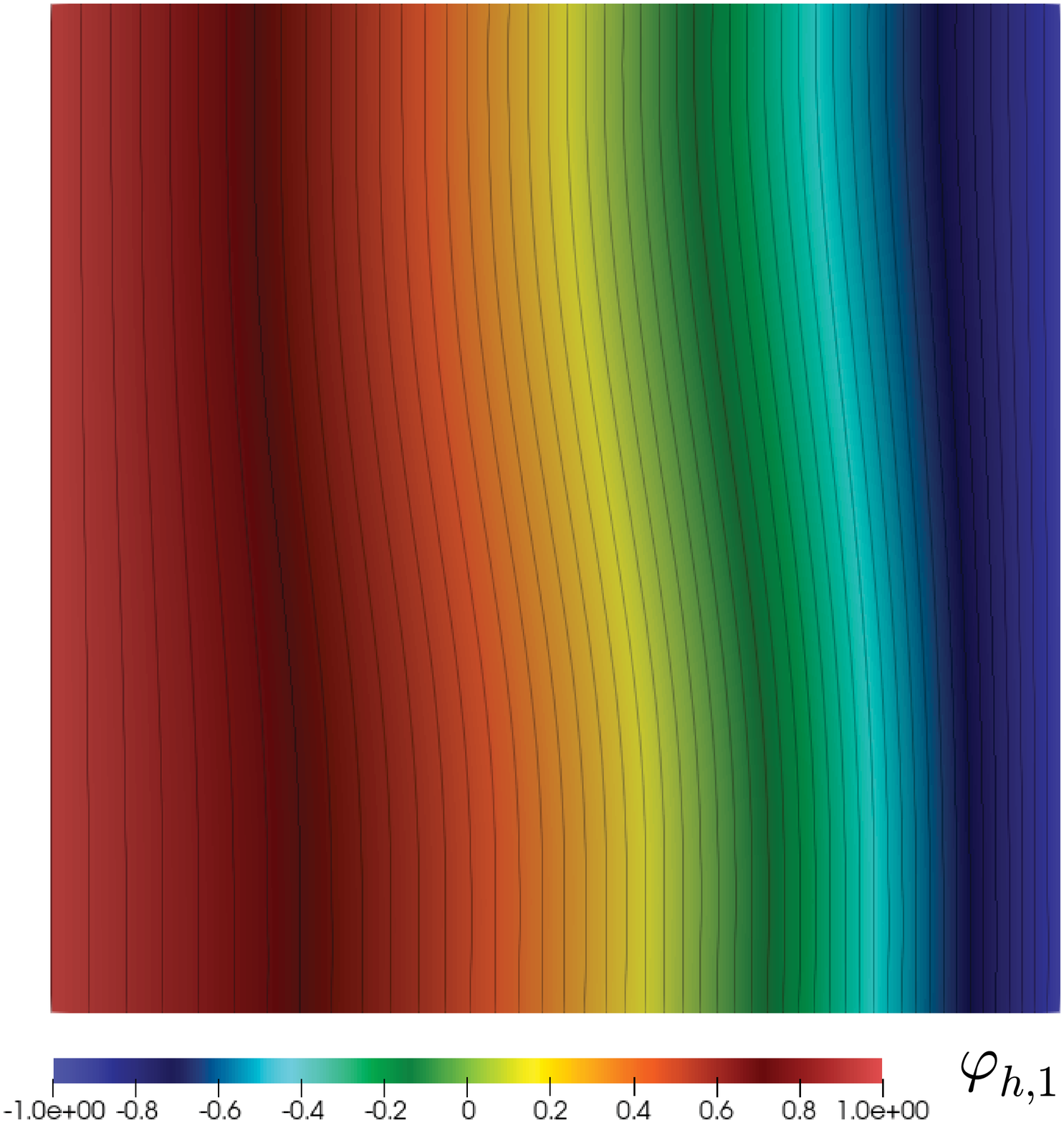}\\[2ex]
\includegraphics{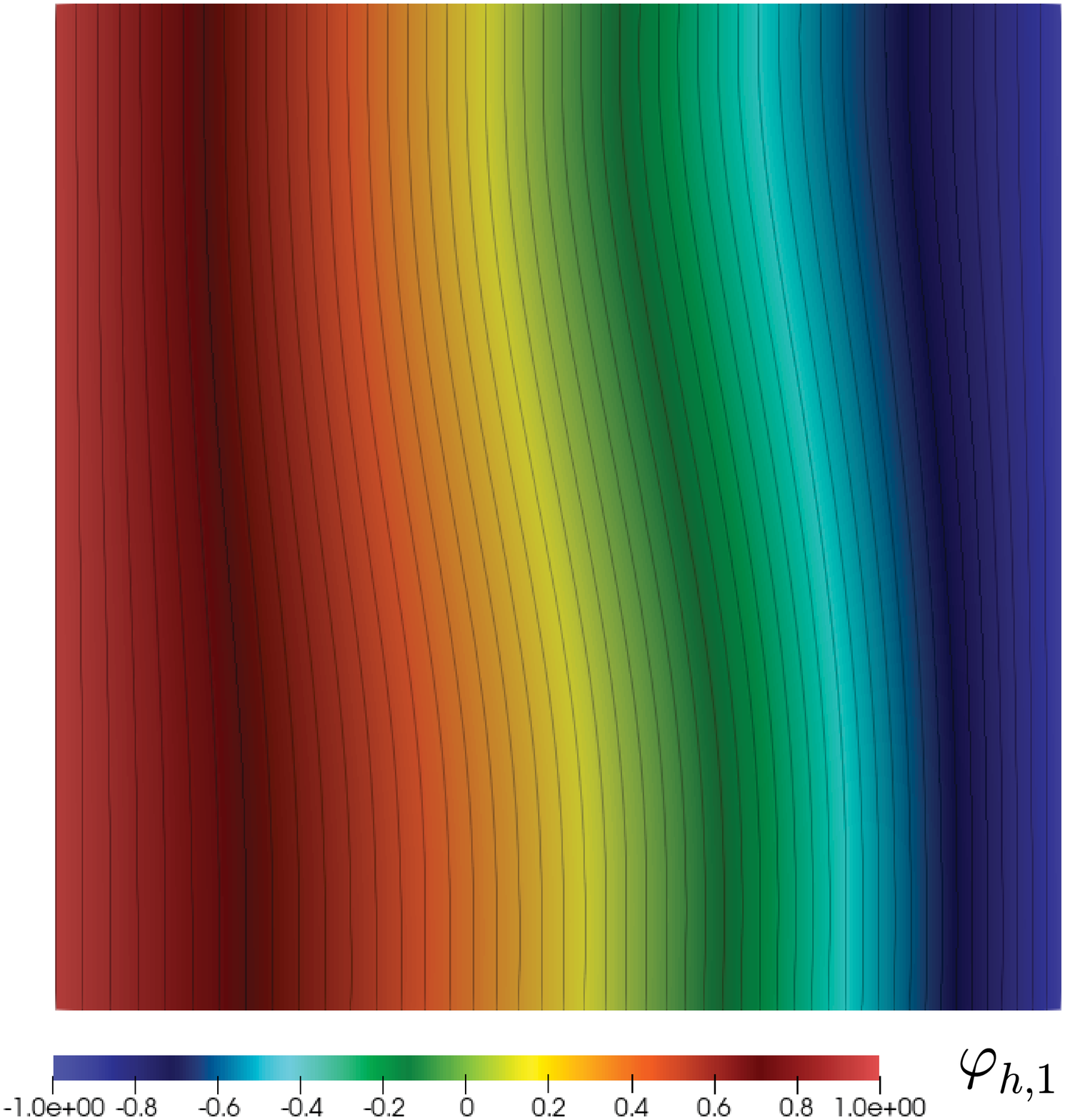}
\end{minipage}\hspace{6cm}\begin{minipage}{25cm}
\includegraphics{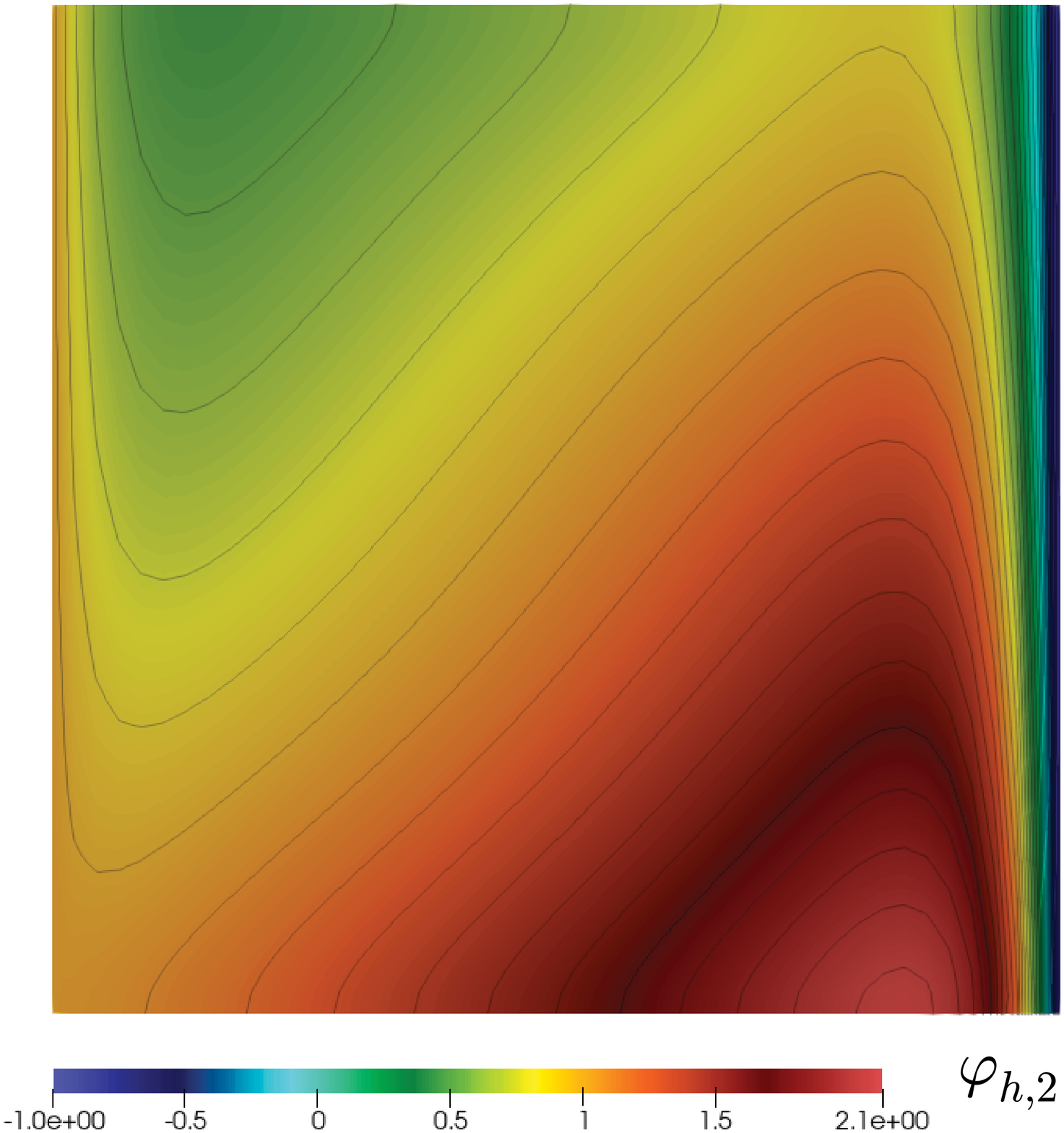}\\[2ex]
\includegraphics{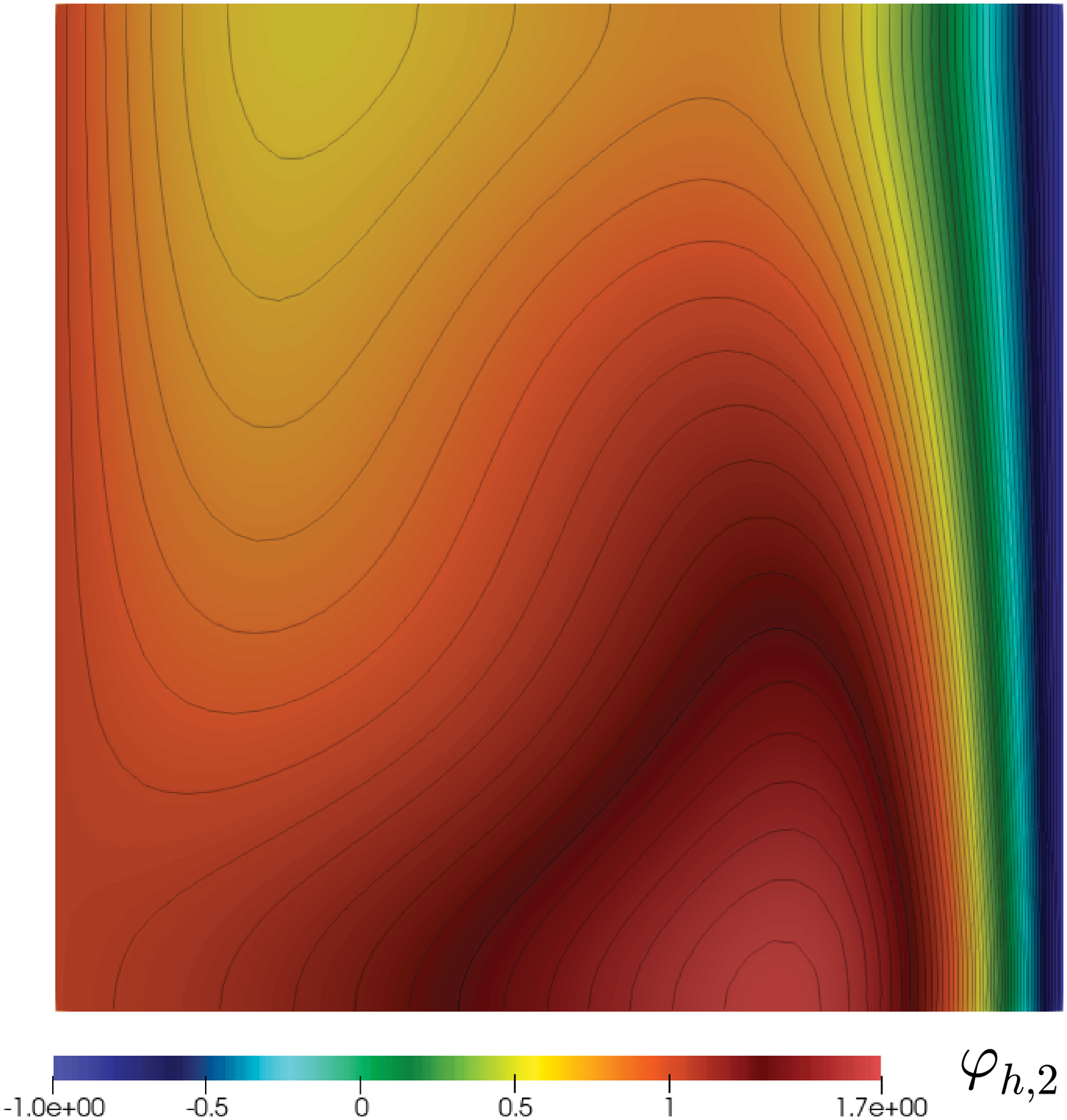}\\[2ex]
\includegraphics{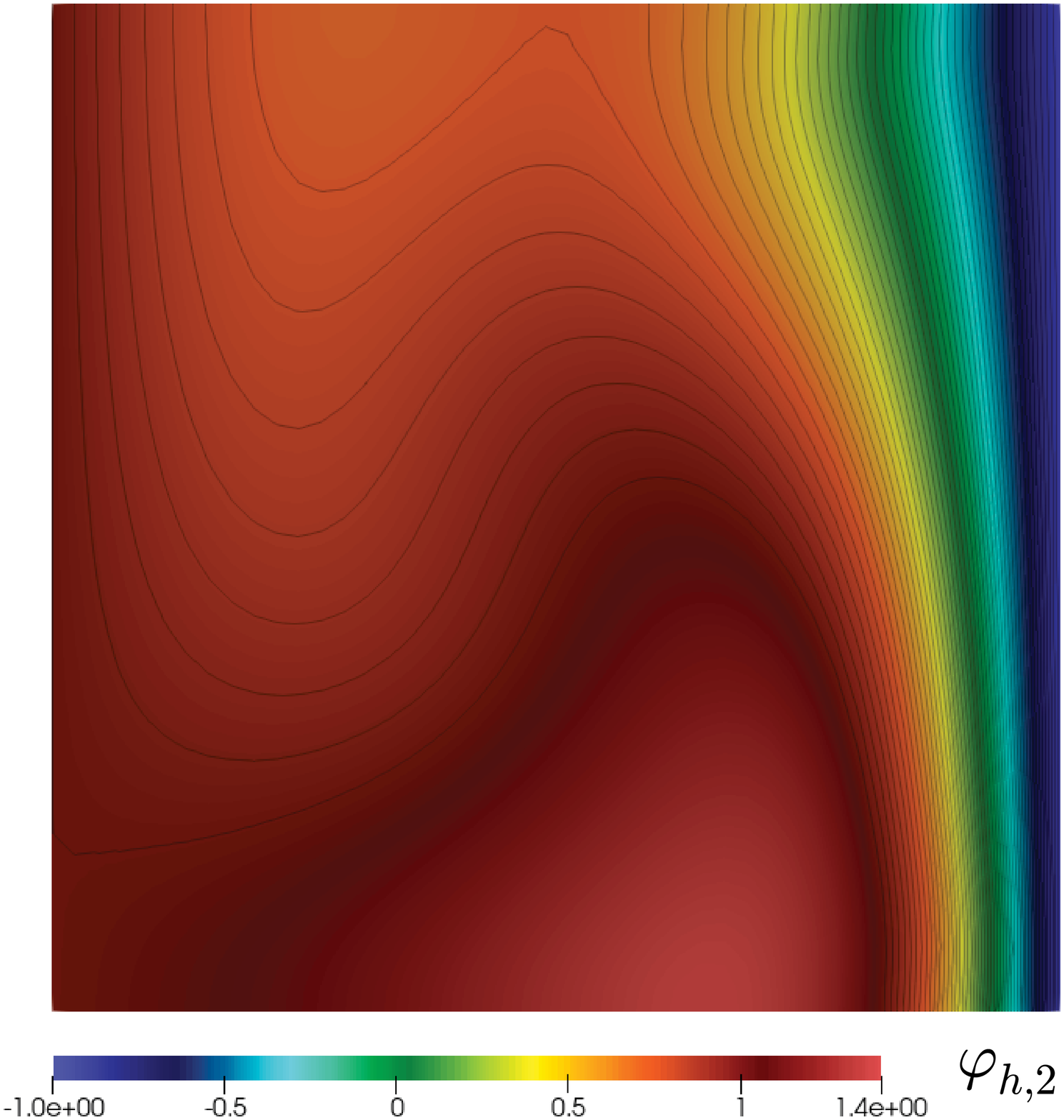}\\[2ex]
\includegraphics{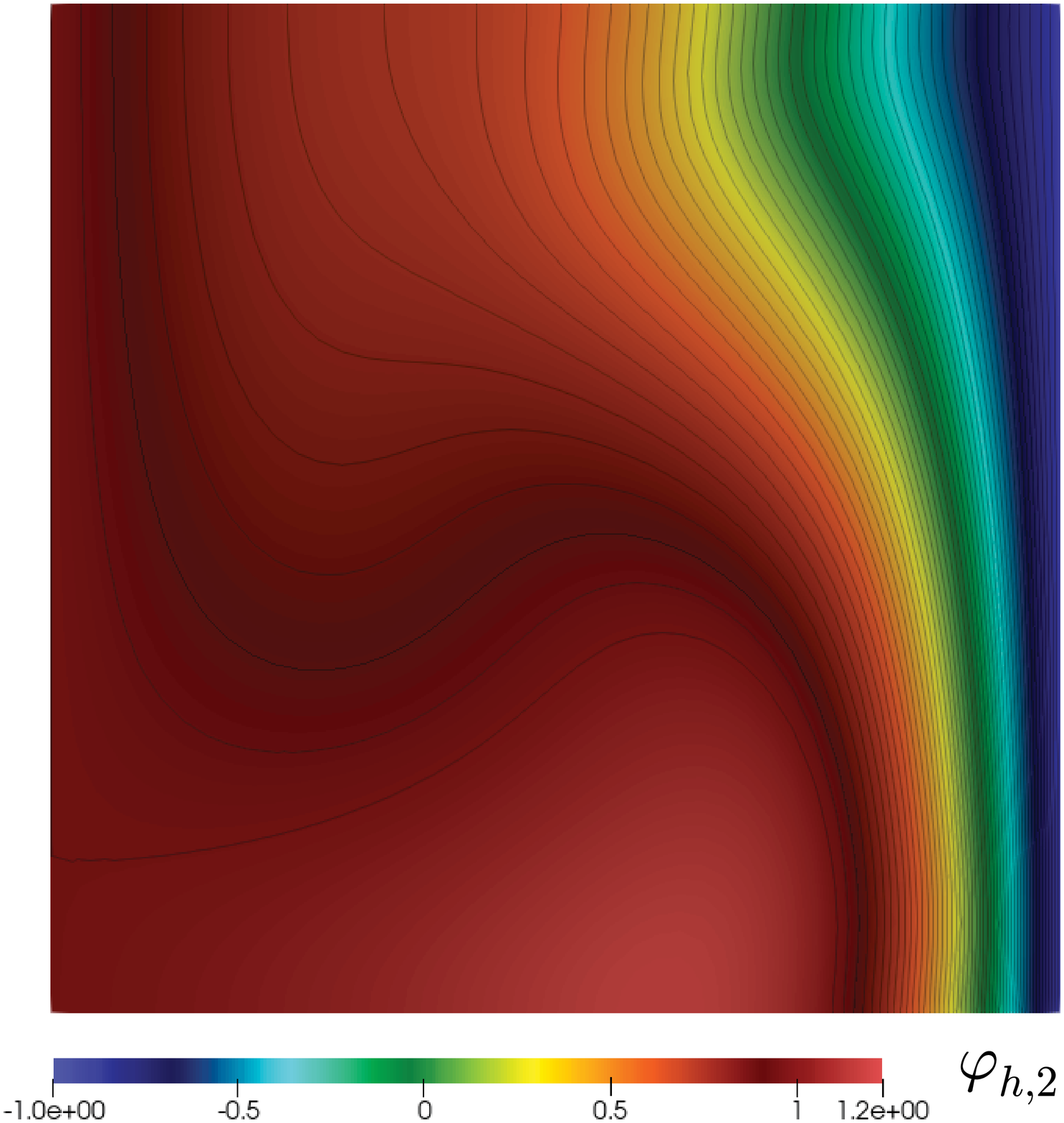}\\[2ex]
\includegraphics{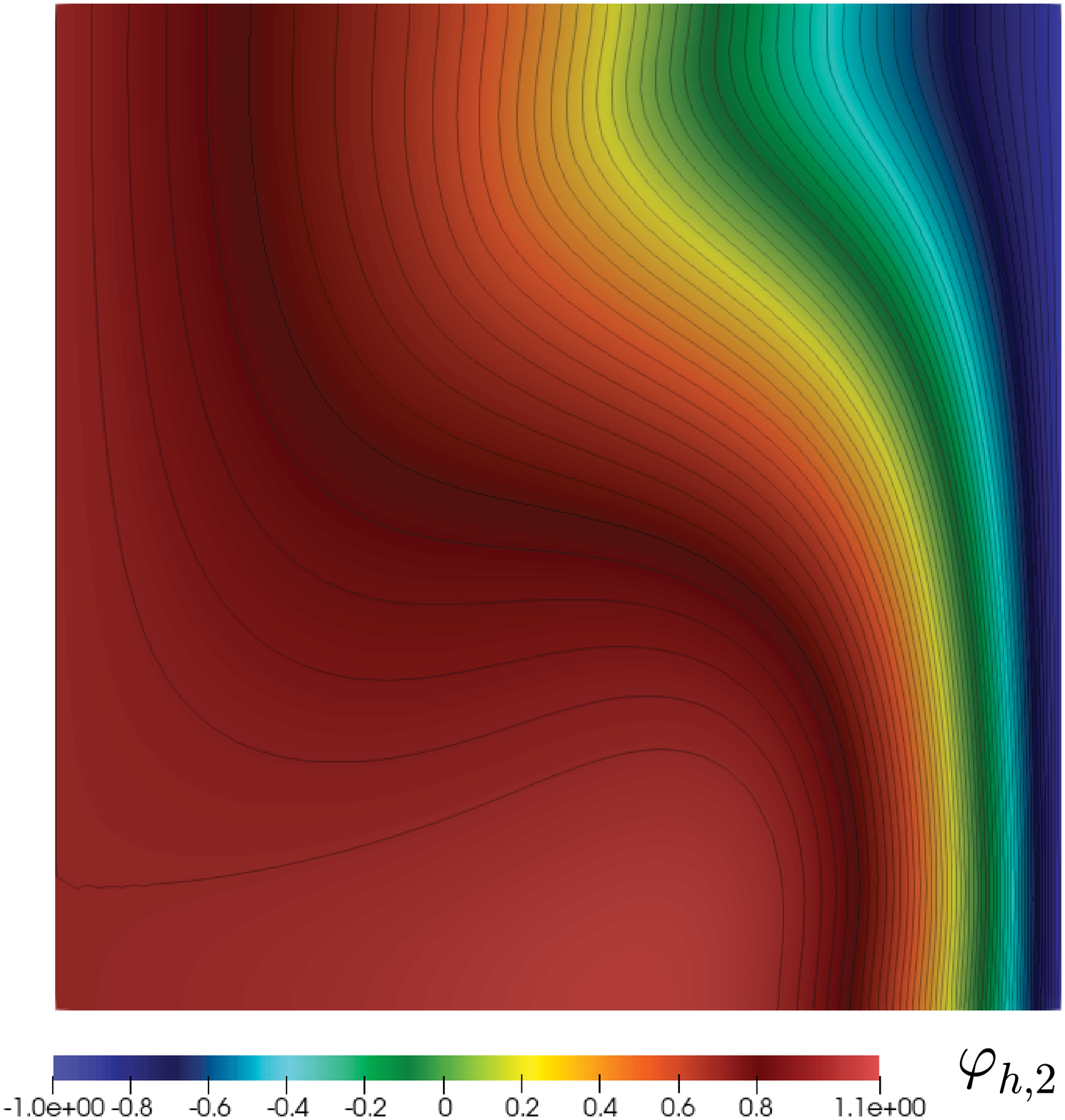}
\end{minipage}}
\caption{Example 3, temperature $\varphi_{h,1}$ (left) and concentration $\varphi_{h,2}$ (right),
obtained with a fully-discrete time-dependent mixed method with no manufactured
analytical solution using $k = 0$ and $N = 28895$ degrees of freedom. We plot
for $t\in\{t_1, t_5, t_{10}, t_{15}, t_{20}\}$, en each row.}\label{fig:exa03-2}
\end{figure}

In Table \ref{tab:exa03}, we summarize the convergence history, where
we can note that the rate of convergence $O(h^{k+1})$ predicted by Theorem
\ref{theoreticalRates} is attained by all the unknowns for $k=0$ and time
step $t_{25} = 0.5$.
We mention that the errors and the convergence rates are computed by
considering the discrete solution obtained with a finer mesh ($N = 28895$)
as the exact solution. Additionally, in Figures \ref{fig:exa03-1} and
\ref{fig:exa03-2}, we display the approximation of the velocity components,
temperature and concentration. All the figures presented there were obtained
with $N = 28895$ degrees of freedom (used as exact solution) in the time
step $t_{\ell} := \ell\cdot\triangle t$, with $\ell\in\{1,5,10,15,20\}$.

\begin{table}[h!t]\centering\footnotesize
\begin{tabular}{|c|c|c|c@{\hspace{1.5ex}}c|c@{\hspace{1.5ex}}c|c@{\hspace{1.5ex}}c|c@{\hspace{1.5ex}}c|c@{\hspace{1.5ex}}c|}\hline
$k$ & $h$ & $N$ & $\mathtt{e}(\bt)$ & $\mathtt{r}(\bt)$ & $\mathtt{e}(\bsig)$ & $\mathtt{r}(\bsig)$ & $\mathtt{e}(\bu)$ & $\mathtt{r}(\bu)$ & $\mathtt{e}(\bvarphi)$ & $\mathtt{r}(\bvarphi)$ & $\mathtt{e}(p)$ & $\mathtt{r}(p)$\\ \hline
  & 0.0707 & 5845 & 4.47e+01 & $--$ & 1.07e+01 & $--$ & 5.79e-01 & $--$ & 5.15e-01 & $--$ & 8.63e-01 & $--$\\
  & 0.0566 & 9055 & 3.57e+01 & 1.01 & 8.35e+00 & 1.10 & 4.54e-01 & 1.08 & 4.02e-01 & 1.12 & 6.86e-01 & 1.03\\
0 & 0.0471 & 12965 & 2.95e+01 & 1.05 & 6.80e+00 & 1.13 & 3.61e-01 & 1.26 & 3.21e-01 & 1.24 & 5.66e-01 & 1.05\\
  & 0.0404 & 17575 & 2.53e+01 & 0.99 & 5.72e+00 & 1.12 & 3.08e-01 & 1.03 & 2.70e-01 & 1.11 & 4.85e-01 & 1.01\\
  & 0.0354 & 22885 & 2.21e+01 & 1.03 & 4.93e+00 & 1.11 & 2.68e-01 & 1.02 & 2.34e-01 & 1.08 & 4.19e-01 & 1.09\\ \hline
\end{tabular}
\caption{History of convergence for Example 3 for $t = 0.5$.}\label{tab:exa03}
\end{table}

\clearpage


\begin{thebibliography}{99}
	

	\bibitem{abedi}
		{\sc J. Abedi, S. Aliabadi},
		{\it Simulation of incompressible flows with heat and mass transfer using parallel finite element method}.
		Electronic J. Diff. Eq., Conference 10 (2003), 1--11.
		




\bibitem{adams-2003}
{\sc R.A. Adams and J.J.F. Fournier},
{\it Sobolev Spaces}.
Second edition. Pure and Applied Mathematics (Amsterdan), 140. Elsevier/Academic Press, Amsterdam, 2003.


		\bibitem{dimitri} \textsc{G. Alekseev, D. Tereshko and V. Pukhnachev}. Boundary control problems for Oberbeck-Boussinesq model of heat and mass transfer. Advanced Topics in Mass Transfer (2011), 485--512.


\bibitem{almonacid-2018}
{\sc J.A. Almonacid, G.N. Gatica and R. Oyarz\'ua},
{\it A new mixed finite element method for the n-dimensional Boussinesq problem with temperature-dependent viscosity}.
Calcolo 55 (2018), no. 3, Art. 36, 42pp.



\bibitem{aggr19}
{\sc M. Alvarez, G.N. Gatica, B. G\'omez-Vargas and R. Ruiz-Baier},
{\it New mixed finite element methods for natural convection with phase-change in porous media}.
J. Sci. Comput. 80 (2019), no. 1, 141--174.


\bibitem{alvarez-2016}
{\sc M. Alvarez, G.N. Gatica and R. Ruiz-Baier},
{\it A mixed-primal element approximation of a sedimentation-consolidation system}.
M3AS: Math. Models Methods Appl. Sci., 26 (2016), no. 5, 867--900.


\bibitem{alvarez-2015}
{\sc M. Alvarez, G.N. Gatica and R. Ruiz-Baier},
{\it An augmented mixed-primal finite element method for a coupled flow-transport problem}.
ESAIM Math. Model. Numer. Anal. 49 (2015), no. 5, 1399--1427.




\bibitem{alvarez-2019}
{\sc M. Alvarez, G.N. Gatica, B. Gomez-Vargas, R. Ruiz-Baier},
{\it New Mixed Finite Element Nethods For Natural Convection with Phase-Change in Porous Media}.
Journal of Scientific Computing. vol. 80, 1, pp. 141-174, (2019).


\bibitem{brezzi-1991}
{\sc F. Brezzi and M. Fortin},
{\it Mixed and Hybrid Finite Element Methods}.
Springer Series in Computational Mathematics, 15. Springer-Verlag, New York, 1991.

\bibitem{bmr19}
{\sc R. B\"urger, P.E. M\'endez and R. Ruiz-Baier},
{\it On H(div)-conforming methods for double-diffusion equations in porous media.}  
SIAM J. Numer. Anal. 57 (2019), no. 3, 1318--1343.



\bibitem{cg-vem-stokes}
{\sc E. C\'aceres and G.N. Gatica},
{\it A mixed virtual element method for the pseudostress-velocity formulation of the Stokes problem}.
IMA J. Numer. Anal. 37 (2017), no. 1, 296--331.

\bibitem{camano-2016}
{\sc J. Cama\~no, G.N. Gatica, R. Oyarz\'ua and G. Tierra},
{\it An augmented mixed finite element method for the Navier-Stokes equations with variable viscosity}.
SIAM J. Numer. Anal. 54 (2016), no. 2, 1069--1092.

\bibitem{caucao-2017}
{\sc S. Caucao, G. N. Gatica and R. Oyarz\'ua},
{\it  Analysis of an augmented fully-mixed formulation for the coupling of the Stokes and heat equations}.
ESAIM Math. Model. Numer. Anal. 52 (2018), no. 5, 1947--1980.

\bibitem{Yuan2016} 	
{\sc Y.Y. Chen, B.W. Li and J.K. Zhang}, 
{\it Spectral collocation method for natural convection in a square porous cavity with local thermal equilibrium and non-equilibrium models}.
Int. J. Heat Mass Transfer 64 (2013), 35--49.

\bibitem{Cheng79}
{\sc P. Cheng},
{\it Heat transfer in geothermal systems}.
Adv. Heat Transfer 14 (1979), 1--105.

\bibitem{ciarlet-2013}
{\sc P.G. Ciarlet},
{\it Linear and Nonlinear Functional Analysis with Applications}.
Society for Industrial and Applied Mathematics, Philadelphia, PA, 2013.

\bibitem{cibik}
{\sc A. \c{C}ibik and S. Kaya},
\textit{Finite element analysis of a projection-based stabilization method for
	the Darcy-Brinkman equations in double-diffusive convection}.
Appl. Numer. Math. 64 (2013), 35--49.

\bibitem{CGMR-2020}
{\sc E. Colmenares, G.N. Gatica, S. Moraga and R. Ruiz-Baier},
{\it A fully-mixed finite element method for the steady state Oberbeck-Boussinesq system.} SMAI J. Comput. Math. 6 (2020),  125--157.


\bibitem{colm-2016}
{\sc E. Colmenares, G.N. Gatica and R. Oyarz\'ua},
{\it Fixed point strategies for mixed variational formulations of the stationary Boussinesq problem}.
C. R. Math. Acad. Sci. Paris 354 (2016), no. 1, 57--62.

\bibitem{colm-2016-2}
{\sc E. Colmenares, G.N. Gatica and R. Oyarz\'ua},
{\it  Analysis of an augmented mixed-primal formulation for the stationary Boussinesq problem}.
Numer. Methods Partial Differential Equations  32 (2016), no. 2, 445--478.

\bibitem{colm-2017}
{\sc E. Colmenares, G.N. Gatica and R. Oyarz\'ua},
{\it An augmented fully-mixed finite element method for the stationary Boussinesq problem}.
Calcolo 54 (2017), no. 1, 167--205.




\bibitem{dallmann16} 
{\sc H. Dallmann and D. Arndt}, {\it Stabilized finite element 
	methods for the Oberbeck-Boussinesq model}. 
J. Sci. Comput. 69 (2016), no. 1, 244--273.


\bibitem{davis-2004}
{\sc T.A. Davis},
{\it Algorithm 832: UMFPACK V4.3--an unsymmetric-pattern multifrontal method}.
ACM Trans. Math. Software 30 (2004), no. 2, 196--199.

\bibitem{Daya Reddy}
{\sc B. Daya Reddy}
{\it Introductory Functional Analysis: With Appliations to Boundary Value Problems and Finite Elements}.
Texts in Applied Mathematics,  Springer (1998), no. 27.

\bibitem{Ern-2004}
{\sc A. Ern and J.L. Guermond},
{\it Theory and Practice of Finite Elements}, Applied Mathematical Sciences,  Springer (2004), no. 159.

\bibitem{fgm-SISC-2008}
{\sc L.E. Figueroa, G.N. Gatica and A. M\'arquez},
{\it Augmented mixed finite element methods for the stationary 
	Stokes equations}. 
SIAM J. Sci. Comput. 31 (2008/09), no. 2, 1082--1119.

\bibitem{gatica-2014}
{\sc G.N. Gatica},
{\it  A Simple Introduction to the Mixed Finite Element Method. Theory and Applications}.
Springer Briefs in Mathematics. Springer, Cham, 2014.

\bibitem{gatica-2014F}
{\sc G.N. Gatica},
{\it Introducci\'on al An\'alisis Funcional. Teor\'ia y Aplicaciones}.
Editorial Reverte, Barcelona Bogot\'a Buenos Aires Caracas M\'exico, 2014.

\bibitem{ob-2} \textsc{B. Gebhart and L. Pera},
\textit{The nature of vertical natural convection flows resulting from the combined
buoyancy effects of thermal and mass diffusion}. 
Int. J. Heat Mass Transfer 14 (1971), no. 12,  2025--2050.

\bibitem{ob-3} \textsc{B. Gebhart and L. Pera},
\textit{Natural convection flows adjacent to horizontal surfaces resulting from the
combined buoyancy effects of thermal and mass diffusion}. 
Int. J. Heat Mass Transfer 15 (1972), no. 2,  269--278.

\bibitem{gmsh} 
{\sc C. Geuzaine and J.-F. Remacle}, {\it Gmsh: a three-dimensional
	finite element mesh generator with built-in pre- and post-processing facilities}.
Int. J. Numer. Methods Engrg. {79} (2009), no. 11, 1309--1331.


\bibitem{girault-1986}
{\sc V. Girault and P-A. Raviart},
{\it Finite Element Methods for Navier-Stokes Equations. Theory and
	Algorithms}.
Springer Series in Computational Mathematics, 5. Springer-Verlag, Berlin, 1986.

\bibitem{gss-IMA-2017}
{\sc J. Guzm\'an, F.A. Sequeira and C.W. Shu},
{\it H(div) conforming and DG methods for incompressible Euler's equations}. 
IMA J. Numer. Anal. 37 (2017), no. 4, 1733--1771.

\bibitem{ob-1} \textsc{M. Hammami, M. Mseddi and Mounir Baccar},
\textit{Numerical
Study of Coupled Heat and Mass Transfer in a Trapezoidal Cavity}. Eng. Appl. Comput. Fluid Mech. 1 (2007), no. 3,  216--226. 		


\bibitem{Karimi97} 	
{\sc M. Karimi-Fard, M.C. Charrier-Mojtabi and K. Vafai}, 
{\it Non-Darcian effects on double-diffusive convection within a porous medium}.
Numer. Heat Transfer, Part A: Appl. 31 (1997), no. 8, 837--852.

\bibitem{McLean_2000}
{\sc W. McLean}, 
{\it Strongly elliptic systems and boundary integral equations}.
Cambridge University Press, Cambridge, 2000.

\bibitem{Mojtabi2005} 
{\sc A. Mojtabi and M.C. Charrier-Mojtabi}, Double-diffusive convection in porous media. Handbook of Porous Media, Part III. Taylor and Francis, 2005.

\bibitem{moraga} 
\textsc{N. Moraga, C. Zambra, P. Torres and R Lemus-Mondaca},
\textit{Fluid  dinamics, heat and mass transfer modelling by finite volume method for agrofood processes}.
DYNA 78 (2011), 140--149.

\bibitem{Nield99} 
{\sc D.A. Nield and A. Bejan}, Convection in porous media. Second edition. Springer Science and \& Business Media. 
Springer-Verlag, New York, 1999. 

\bibitem{oyarzua-2014}
{\sc R. Oyarz\'ua, T. Qin and D. Sch\"otzau},
{\it An exactly divergence-free finite element method for a generalized Boussinesq problem}.
IMA J. Numer. Anal. 34 (2014), no. 3, 1104--1135.

\bibitem{quartero-1994}
{\sc A. Quarteroni and A. Valli}.
{\it Numerical Approximation of Partial Differential Equations}.
Springer Series in Computational Mathematics, 23. Springer-Verlag, Berlin, 1994.

\bibitem{raviart-1983}
{\sc P.-A. Raviart and J.-M. Thomas},
{\it Introduction \`a l'Analyse Num\'erique des  \'Equations aux D\'eriv\'ees Partielles}.
(French) [Introduction to the numerical analysis of partial differential equations] Collection Math\'ematiques Appliqu\'ees pour la Ma\^itrise. Masson, Paris, 1983.

\bibitem{Thomas-1991}
{\sc J.E. Roberts and J.-M. Thomas},
{\it Mixed and Hybrid Methods}.
Handbook of numerical analysis, vol. II, 523--639, Handb. Numer. Anal., II, North-Holland, Amsterdam, 1991.

\bibitem{ob-5} \textsc{S. C. Saha and A. Hossain},
\textit{Natural convection flow with combined buoyancy effects due to
thermal and mass diffusion in a thermally stratified media}. 
Nonlinear Analysis: Modelling and Control 9 (2004),  89--102.

\bibitem{ScottZhang90}
{\sc L.R. Scott and S. Zhang},
{\it Finite element interpolation of nonsmooth functions satisfying boundary conditions}.
Math. Comp. 54 (1990), no. 190, 483--493.

\end{thebibliography}
\end{document}